\documentclass[a4paper]{amsart}
\usepackage{amsmath,amsthm,amssymb,mathtools,mathrsfs,bm}
\usepackage{tikz-cd}
\usepackage{xcolor}
\usepackage[top=2.3cm,bottom=2.3cm,left=2.3cm,right=2.3cm]{geometry}
\usepackage{setspace}
\usepackage[expansion=false]{microtype}
\usepackage[shortlabels]{enumitem}
\usepackage{needspace}
\usepackage{hyperref}
\definecolor{urlcolor}{rgb}{0,0,0}
\definecolor{linkcolor}{rgb}{.7,.10,.2}
\definecolor{citecolor}{rgb}{.12,.54,.11}
\hypersetup{breaklinks=true,colorlinks=true,urlcolor=urlcolor,
linkcolor=linkcolor,citecolor=citecolor,
pdftitle={Tautological generation of BPS cohomology},
pdfauthor={Ben Davison and Lucien Hennecart}}
\usepackage{cleveref}
\numberwithin{equation}{section}
\newtheorem{theorem}{Theorem}[section]
\newtheorem{proposition}[theorem]{Proposition}
\newtheorem{lemma}[theorem]{Lemma}
\newtheorem{corollary}[theorem]{Corollary}
\theoremstyle{definition}
\newtheorem{definition}[theorem]{Definition}

\theoremstyle{remark}
\newtheorem{remark}[theorem]{Remark}

\newcommand{\BoQ}{\mathbf Q}
\newcommand{\BoC}{\mathbf C}
\newcommand{\BoN}{\mathbf N}
\newcommand{\BoZ}{\mathbf Z}
\newcommand{\dd}{\mathbf d}
\newcommand{\ee}{\mathbf e}

\newcommand{\FM}{\mathfrak M}
\newcommand{\CM}{\mathcal M}
\newcommand{\CA}{\mathcal A}
\newcommand{\Fg}{\mathfrak g}
\newcommand{\CB}{\mathcal B}
\newcommand{\CV}{\mathcal V}
\newcommand{\CL}{\mathcal L}
\newcommand{\CE}{\mathcal E}
\newcommand{\CF}{\mathcal F}
\newcommand{\HO}{\mathrm H}
\newcommand{\GL}{\mathrm{GL}}
\newcommand{\SL}{\mathrm{SL}}
\newcommand{\BG}{\mathrm B}
\newcommand{\Gm}{\mathbf G_{\mathrm m}}
\newcommand{\BM}{\mathrm{BM}}
\newcommand{\rd}{\mathrm{rd}}
\newcommand{\vir}{\mathrm{vir}}

\newcommand{\sgn}{\mathrm{sgn}}

\newcommand{\calBPS}{\mathcal{BPS}}
\newcommand{\JH}{\mathrm{JH}}
\newcommand{\Tr}{\operatorname{Tr}}
\newcommand{\tr}{\operatorname{tr}}
\newcommand{\Rep}{\operatorname{Rep}}
\newcommand{\End}{\operatorname{End}}
\newcommand{\Hom}{\operatorname{Hom}}
\newcommand{\Sym}{\operatorname{Sym}}

\newcommand{\Prim}{\operatorname{Prim}}
\newcommand{\Flag}{\operatorname{Flag}}
\newcommand{\rk}{\operatorname{rk}}
\newcommand{\im}{\operatorname{im}}
\newcommand{\ad}{\operatorname{ad}}
\newcommand{\vdim}{\operatorname{vdim}}

\newcommand{\Hilb}{\operatorname{Hilb}}

\newcommand{\id}{\mathrm{id}}
\newcommand{\Fr}{\mathfrak F}
\newcommand{\Zleg}{\mathfrak Z}
\newcommand{\Y}{\mathfrak Y}
\newcommand{\XX}{\mathfrak X}
\newcommand{\IC}{\operatorname{IC}}

\usepackage[dvipsnames]{xcolor}

\title[Tautological generation of BPS cohomology]{Tautological generation of preprojective BPS cohomology}

\author{Ben Davison}

\address{School of Mathematics, University of Edinburgh, Edinburgh, UK}
\email{Ben.Davison@ed.ac.uk}

\author{Lucien Hennecart}

\address{Laboratoire Ami\'enois de Math\'ematique Fondamentale et Appliqu\'ee, CNRS UMR 7352, Universit\'e de Picardie Jules Verne, 33 rue Saint Leu, 80000 Amiens, France}

\email{lucien.hennecart@u-picardie.fr}
\date{\today}

\begin{document}

\begin{abstract}
Let $Q$ be a quiver and let $\dd$ be a dimension vector. We prove that the affinized BPS cohomology of the stack of $\dd$-dimensional representations of the preprojective algebra $\Pi_Q$ is generated by its reduced virtual fundamental class under the action of tautological classes. This is a generalisation to singular quiver varieties of the Kirwan surjectivity statement of McGerty--Nevins for smooth framed quiver varieties. We prove that the statement holds for every King stability condition. To prove this, we attach legs to vertices of $Q$ and construct a correspondence with a smooth quiver variety for the extended quiver. This correspondence sends the reduced virtual fundamental class of the singular preprojective quiver stack to the Vandermonde determinant capped with the fundamental class of the smooth preprojective quiver stack for the leg-extended quiver. We construct the Weyl group action by the simultaneous Springer resolution of Dobrovolska--Ginzburg--Travkin, providing a comparison of the cohomology of fibres. This comparison preserves the tautological Chern classes. The sign-isotypic formula involving Kac polynomials of Hausel--Letellier--Rodriguez-Villegas, combined with Kirwan surjectivity for smooth framed quiver varieties of McGerty--Nevins, then implies the tautological generation theorem. We also describe the sign isotypic component of the cohomology of quiver varieties using the real-root subalgebras of the BPS Lie algebra. The tautological generation result can be extended via similar methods to deformed preprojective algebras of Crawley-Boevey--Holland at central parameters.
\end{abstract}
\maketitle
\enlargethispage{3pt}
\setcounter{tocdepth}{1}
\tableofcontents

\section{Introduction}
\subsection{Background}
The cohomology of a smooth moduli space is often generated by characteristic classes of its universal object. For smooth Nakajima quiver varieties, this is the Kirwan surjectivity theorem of McGerty--Nevins \cite{MN}. It states that the cohomology of any smooth Nakajima quiver variety is generated by the Chern classes of the tautological family of representations of the preprojective algebra it parametrizes. Going beyond the linear setting of quiver representations, the McGerty--Nevins theorem is analogue of the theorem of Markman, stating that cohomology rings of moduli spaces of sheaves on symplectic surfaces are similarly generated by tautological classes \cite{markman2007integral}. 

For singular moduli spaces, ordinary singular cohomology does not usually have the same relation with the geometry of the moduli problem and is in general harder to understand. In the case of preprojective algebras, the Borel--Moore homology of the moduli stack carries a cohomological Hall algebra structure \cite{schiffmann2013cherednik,yang2018cohomological}, and its BPS Lie algebra \cite{DM,DavBPS}, constructed from the corresponding tripled quiver with canonical cubic potential, provides a distinguished part of this homology. A first fact showing that the BPS cohomology is a meaningful subspace of the Borel--Moore homology of the stack of representations of the preprojective algebra is the identification of its Poincaré polynomial with the Kac polynomials of the quiver \cite{davison2018purity}, providing a new proof of the Kac positivity conjecture \cite[Conjecture~2]{kac2006root}, first proven in \cite{HLRV} and later in \cite{DGT}. Combining ideas from all these sources, in this paper we show that the BPS cohomology is tautologically generated, and so provides the language for generalising the McGerty--Nevins theorem past the smooth setting.

The methods of McGerty--Nevins also give the tautological generation for the pure part of the cohomology of smooth multiplicative quiver varieties \cite{MNmult}. They are multiplicative analogues of quiver varieties, first considered by Crawley-Boevey and Shaw in their work on the Deligne--Simpson problem \cite{CBS}.

The aim of this paper, see \Cref{thm:main}, is to identify the part of the Borel--Moore homology of the stack of representations of a preprojective algebra generated by tautological classes from a single geometric class: the \emph{reduced virtual fundamental class} (\Cref{sec:virtual}). The reduction removes the trivial one-dimensional part of the obstruction space corresponding to simultaneous scalar automorphisms. The resulting class can be nonzero even though the ordinary virtual fundamental class vanishes, and it need not agree with the ordinary fundamental class of the classical stack (which may not have the correct dimension or may not be irreducible). It does however coincide with the latter when the stack is irreducible, generically reduced, and its reduced virtual dimension and classical dimension coincide, see \Cref{lem:stablefund}. A hint that the BPS cohomology should be generated by a single class under tautological classes is that Kac polynomials are monic polynomials \cite[\S1.15, Proposition]{kac2006root}. The leading coefficient $1$ then corresponds to the reduced virtual fundamental class (\Cref{cor:lowest}).

The \emph{affinized BPS cohomology} admits an intrinsic description as the primitive subspace for a factorisation coproduct (\Cref{subsection:factorisation-coproduct}), see \cite{DHKSV} for a general definition, and \cite{JinThesis,Jin} for detailed versions that apply for categories of modules over preprojective algebras. The coproduct is constructed by separating the spectrum of the central endomorphism given by the action of $\omega$ in representations of $\Pi_Q[\omega]$ (\Cref{subsection:factorisation-coproduct}), which is isomorphic to the Jacobi algebra of the tripled quiver with canonical cubic potential (\Cref{equation:jacobi-algebra}). The scalar-microlocal description of the reduced virtual class then implies that every tautological multiple of this class is primitive (\Cref{prop:tautprimitive}). The tautological generation theorem (\Cref{thm:main}) amounts to proving the opposite inclusion, that is that every class in the affinized BPS cohomology is a tautological multiple of the reduced virtual fundamental class (\Cref{subsection:dimension-argument}). This is done by considering a flag correspondence between the singular preprojective stack and a smooth framed preprojective stack, \Cref{eq:flagcorr}, and by identifying the action of the flag map \eqref{equation:flag-map} on the affinized BPS cohomology.

The tautological generation result presented in this paper is a preprojective algebra version of the tautological generation of the affinized BPS cohomology for symplectic surfaces presented in \cite{DHKSV}. For symplectic surfaces, the tautological generation of the affinized BPS cohomology is a generalisation of Markman's theorem \cite{markman2007integral} regarding the tautological generation of the cohomology of smooth moduli spaces. This tautological generation provides a characterisation of the BPS cohomology and of its affinized version that do not require the consideration of the tripled quiver with the canonical cubic potential. Indeed, the reduced virtual fundamental class can be defined within the geometry of the preprojective algebra (using the moment map, \Cref{sec:virtual}) and so can the action of tautological classes. We expect such a phenomenon to be more general, see in particular \Cref{subsection:generation-2CY} for a question raised in the $2$-Calabi--Yau context. We highlight an interesting generalisation of the question raised in \Cref{subsection:generation-2CY} to $0$-shifted symplectic stacks but we do not take the risk of making a conjecture.

\subsection{Main results}
Let $Q=(Q_0,Q_1)$ be a finite quiver with set of vertices $Q_0$ and set of arrows $Q_1$, possibly with loops and multiple arrows. We write $\Pi_Q$ for its preprojective algebra (\Cref{subsection:quivers-moment-maps}) and for a dimension vector $\dd\in\BoN^{Q_0}$, we define the product of general linear groups
\[
 G_{\dd}=\prod_{i\in Q_0}\GL_{d_i}.
\]
For a King stability parameter $\theta$ satisfying $\theta\cdot\dd=0$ (see \cite{King}, \Cref{subsection:generic} for details), we let $\FM_{\Pi_Q,\dd}^{\theta}$ be the stack of $\theta$-semistable $\Pi_Q$-modules of dimension vector $\dd$, and let $\CM_{\Pi_Q,\dd}^{\theta}$ be its good moduli space (\Cref{subsection:quivers-moment-maps}), defined as a GIT quotient. Throughout this paper, we work with constructible sheaves with rational coefficients, although careful readers will quickly notice that all arguments and results extend to derived categories of mixed Hodge modules. In the framework of quiver varieties, where all mixed Hodge modules of relevance (in particular, the BPS sheaf) are pure and all mixed Hodge structures are pure Tate by \cite[Corollary~6.2]{DavInt}, working with constructible sheaves is sufficient to deduce immediately the result in the refined context of mixed Hodge modules/structures.

We denote by $\Fg_{\Pi_Q,\dd}^{\theta}$ the BPS cohomology (\Cref{subspace:BPS-subspace}) and by $\widehat{\Fg}_{\Pi_Q,\dd}^{\theta}$ its affinization, realised as the primitive subspace of the Borel--Moore homology of $\FM_{\Pi_Q,\dd}^{\theta}$ for a cocommutative coproduct (\Cref{subsection:factorisation-coproduct}). As graded vector spaces, there is an isomorphism (see \Cref{thm:PBW} and \Cref{prop:primitives})
\[
 \widehat{\Fg}_{\Pi_Q,\dd}^{\theta}
 \cong \Fg_{\Pi_Q,\dd}^{\theta}\otimes\BoQ[u],
 \qquad \deg u=2.
\]
The embedding of BPS cohomology and its affinized version in the Borel--Moore homology of the stack $\FM_{\Pi_Q,\dd}^{\theta}$ and the cohomological degree conventions are recalled in \Cref{sec:DT}. In particular, the cohomological degrees in Borel--Moore homology are shifted by the virtual dimension of the stacks under consideration \eqref{eq:cohgrading}. Let
\[
 R_{\dd}=\HO^*(\BG G_{\dd},\BoQ)
\]
be the $G_{\dd}$-equivariant cohomology of the point. The Chern classes of universal bundles $\CV_i$ at vertices of $Q$ give a ring homomorphism $R_{\dd}\rightarrow\HO^*(\FM_{\Pi_Q,\dd}^{\theta},\BoQ)$ and so give an action of $R_{\dd}$ on the Borel--Moore homology of $\FM_{\Pi_Q,\dd}^{\theta}$ by cap-product.

Our first main theorem is the following.
\begin{theorem}\label{thm:main}
Let $Q$ be a quiver, let $\dd\in\BoN^{Q_0}\setminus\{0\}$ be a non-zero dimension vector for $Q$, and let $\theta$ be a King stability
parameter (i.e. $\theta\cdot\dd=0$, \Cref{subsection:generic}). Then
\[
 R_{\dd}\cap[\FM_{\Pi_Q,\dd}^{\theta}]^{\rd}
   =\widehat{\Fg}_{\Pi_Q,\dd}^{\theta}
 \subseteq \HO_*^{\BM}(\FM_{\Pi_Q,\dd}^{\theta},\BoQ).
\]
\end{theorem}
\Cref{thm:main} is proven in \Cref{subsection:dimension-argument}.  Here and below, the reduced virtual fundamental class $[\FM_{\Pi_Q,\dd}^{\theta}]^{\rd}$ is defined using the reduced derived structure of $\FM_{\Pi_Q,\dd}^{\theta}$, not the classical stack or the unreduced derived structure, see \Cref{sec:virtual}. The theorem includes all dimension vectors, including those for which the stack $\FM_{\Pi_Q,\dd}^{\theta}$ is not irreducible. For example, when $Q$ is the $1$-Kronecker quiver (type $\mathrm{A}_2$), $\dd=(1,1)$ and $\theta=0$, it is the quotient stack of $\{xy=0\}\subset\BoC^2$ by an action of $\BoC^*\times\BoC^*$. When $\dd$ is a dimension vector such that there exist $\theta$-stable $\Pi_Q$-representations of that dimension, it is known by \cite{CB, BS21} that the stack $\FM_{\Pi_Q,\dd}^{\theta}$ is irreducible and in that case, the classical fundamental class of the stack coincides with the reduced virtual fundamental class. This comes from the fact that generically, in that case, the reduced derived structure agrees with the classical stack, see \Cref{lem:stablefund}.

We write $(- ,-)_Q$ for the symmetrised Euler form of $Q$ (\Cref{subsection:quivers-moment-maps}) and set
\[
 p_Q(\dd)=1-\frac12(\dd,\dd)_Q.
\]
We consider the following grading on the preprojective CoHA:
\[
 \CA_{\Pi_Q,\dd}^k
   =\HO_{-k-(\dd,\dd)_Q}^{\BM}(\FM_{\Pi_Q,\dd},\BoQ).
\]
Equivalently, this is the hypercohomology of the Verdier dual of $\BoQ_{\FM_{\Pi_Q,\dd}}^{\vir} \coloneqq\BoQ_{\FM_{\Pi_Q,\dd}}[-(\dd,\dd)_Q]$; see \Cref{eq:cohgrading}. The shift is the derived virtual dimension of the stack $\FM_{\Pi_Q,\dd}$. The reduced virtual class has degree $-2p_Q(\dd)$ (\Cref{lem:vdim}). When the Kac polynomial $A_{Q,\dd}$ of $Q$ in dimension $\dd$ is nonzero, this is minus twice its degree, see \cite{HLRV}, and also \cite[\S1.15, Proposition]{kac2006root} for the original work of Kac.

\begin{corollary}\label{cor:lowest}
If the Kac polynomial $A_{Q,\dd}(q)$ is nonzero, then
\[
 \Fg_{\Pi_Q,\dd}^{-2p_Q(\dd)}=\BoQ\cdot[\FM_{\Pi_Q,\dd}]^{\rd}.
\]
In particular the reduced virtual class $[\FM_{\Pi_Q,\dd}]^{\rd}$ is nonzero and belongs to the unaffinized (i.e. classical) BPS cohomology $\Fg_{\Pi_Q,\dd}$. If $A_{Q,\dd}=0$, this class is zero. The corresponding assertions hold after restriction to any King stability condition $\theta$ with $\theta\cdot\dd=0$:
\[
  \Fg_{\Pi_Q,\dd}^{\theta,-2p_Q(\dd)}=\BoQ\cdot[\FM_{\Pi_Q,\dd}^{\theta}]^{\rd}.
\]

\end{corollary}

There is also a consequence for smooth quiver varieties, generalising the Kirwan surjectivity theorem of McGerty--Nevins \cite{MN}.

\begin{corollary}\label{cor:smooth}
Assume that every $\theta$-semistable representation of the preprojective algebra $\Pi_Q$ of dimension $\dd$ is stable (for example, $\theta$ is generic for $\dd$, see \Cref{subsection:generic}). Then the cohomology of the stack $\FM_{\Pi_Q,\dd}^{\theta}$ is generated by the Chern classes of the universal bundles. The cohomology of the smooth coarse moduli space $\CM_{\Pi_Q,\dd}^{\theta}$ is generated by tautological classes.
\end{corollary}

We prove \Cref{cor:smooth} in \Cref{subsection:proofs}. By \Cref{sec:deformation}, the version of \Cref{cor:smooth} for deformed preprojective algebras $\Pi_Q^{\lambda}$ holds as well. For the coarse moduli space, tautological classes are obtained by descent from the stack after removing the common scalar torus acting trivially; see \Cref{subsec:gerbe}. The work of McGerty--Nevins \cite{MN} proves tautological generation for smooth \emph{framed} quiver varieties (i.e. there is a vertex with dimension vector $1$). Therefore, \Cref{cor:smooth} is a generalisation to any smooth quiver varieties of their result. Note that the techniques developed in the paper on the tautological generation of the pure cohomology of multiplicative quiver varieties \cite{MNmult} provide a path to obtain \Cref{cor:smooth}, i.e. tautological generation for any \emph{smooth} quiver varieties, not only framed ones, via other means. However, the argument is not written anywhere and so \Cref{cor:smooth} is, to the best of our knowledge, original.

\subsection{Quivers with legs and sign isotypic components}
We recall the strategy of attaching legs to a given quiver to produce smooth quiver varieties, following \cite{CB03Matrices,maffei2005quiver,Letellier,ballandras2023comet}, which rests on the realisation of conjugacy classes of matrices by quiver varieties \cite{kraft1979closures,kraft1981minimal}, see also \cite{kraft1982geometry} for the generalisation to classical groups. For every $i\in Q_0$ with $d_i>0$, attach a leg with dimensions $d_i,\ d_i-1,\ldots,1$, where the first vertex is the original vertex $i$. For convenience, we label the vertices of the leg by $(i,0),\hdots,(i,d_i-1)$. Therefore, with this labeling, $i=(i,0)$. We denote the resulting quiver and dimension vector by $Q^+$ and $\dd^+$ respectively. For example, if $Q$ is the quiver with one vertex and $g$ loops and the dimension vector is $4$, then $Q^+$ is the following quiver:
\[
\begin{tikzcd}[column sep=large]
\bullet_{4} \arrow[loop left, "g\text{ loops}"]
& \bullet_{3} \arrow[l]
& \bullet_{2} \arrow[l]
& \bullet_{1} \arrow[l]
\end{tikzcd}
\]
where we have also indicated the dimension vector $\dd^+$.  We define
\begin{equation}\label{equation:weyl-group}
 N_{\dd}=\sum_i\binom{d_i}{2},\qquad W_{\dd}=\prod_i\mathfrak S_{d_i}.
\end{equation}
Note that $N_{\dd}$ is the dimension of the flag variety of $\GL_{\dd}$ and $W_{\dd}$ is its Weyl group. There exists a generic stability parameter $\Theta$, positive at the added vertices (\Cref{lem:chamber}), which gives a smooth stack $\Y=\FM_{\Pi_{Q^+},\dd^+}^{\Theta}$, since $\Theta$-semistability and stability coincide for dimension vector $\dd^+$.  It follows that $\FM_{\Pi_{Q^+},\dd^+}^{\Theta}$ is a $\mathbf{G}_{\mathrm{m}}$-gerbe over its moduli space $Y\coloneqq \CM_{\Pi_{Q^+},\dd^+}^{\Theta}$ \eqref{equation:smoothstable}. There is a flag correspondence \eqref{eq:flagcorr} defining a map
\[
 \Phi:\CA_{\Pi_Q,\dd}\longrightarrow\HO_*^{\BM}(\Y,\BoQ),
\]
by pullback, pushforward and restriction \eqref{equation:flag-map}, which preserves the normalised CoHA grading, i.e. the shifts by the respective virtual dimensions of $\FM_{\Pi_Q,\dd}$ and $\Y$ (\Cref{lem:Phi-linear}).

On $\Y$, the inward leg maps (pointing to the initial vertex $(i,0)$) are injective, as a consequence of the choice of the stability $\Theta$ (\Cref{lem:injectivity}). The universal representation over $\Y$ gives injective morphisms between universal vector bundles $\CV_{i,r}\rightarrow\CV_{i,r-1}$ over $\Y$ for any $i\in Q_0$, $1\leq r\leq d_i$ (\Cref{subsection:equivariance-tautological-classes}). By convention, we write $\CV_{i,d_i}=0$. We write $\CL_{i,r}=\CV_{i,r-1}/\CV_{i,r}$ for the successive line bundles over $\Y$ and define
$x_{i,r}=c_1(\CL_{i,r})$, the first Chern class of the line bundle $\CL_{i,r}$. We set
\begin{equation}\label{equation:vandermonde}
 P_{\dd}=\BoQ[x_{i,r}:i\in Q_0,\ 1\leq r\leq d_i], \qquad  \Delta_{\dd}=\prod_i\prod_{r<s}(x_{i,r}-x_{i,s})
\end{equation}
for the polynomial ring and the Vandermonde determinant. Then we can identify the invariant ring $P_{\dd}^{W_{\dd}}=R_{\dd}$ where $W_{\dd}$ acts by permuting the $x_{i,r}$, $1\leq r\leq d_i$, independently for any vertex $i\in Q_0$.

\begin{theorem}\label{thm:signintro}
The flag correspondence induces an $R_{\dd}$-linear isomorphism
\[
 \widehat{\Fg}_{\Pi_Q,\dd} \xrightarrow{\sim} \HO^*(\Y,\BoQ)^{\sgn}\cap[\Y]
\]
for the geometric Weyl action. We have
\[
 \Phi([\FM_{\Pi_Q,\dd}]^{\rd})
       =(-1)^{N_{\dd}}\Delta_{\dd}\cap[\Y].
\]
\end{theorem}

The determination of the graded dimension of the sign-isotypic component of the target $\HO^*(\Y,\BoQ)^{\sgn}\cap[\Y]$ is due to Hausel--Letellier--Rodriguez-Villegas \cite{HLRV}, see \Cref{thm:HLRV} and \Cref{cor:signcharacter}. \Cref{thm:signintro} is proven in \Cref{subsection:dimension-argument}. Alternatively, one can avoid the use of the HLRV formula (\Cref{thm:HLRV}) and instead use the realisation of the flag correspondence by leg operators in BPS cohomology, as explained in \Cref{sec:Lie}, to obtain the injectivity of the map $\Phi$ defined in \eqref{equation:flag-map} after restriction to $\widehat{\Fg}_{\Pi_Q,\dd}$.

\subsection{Strategy of the proof}
\label{subsection:strategy-proof}
We briefly outline the strategy to prove our main theorems, \Cref{thm:main} and \Cref{thm:signintro}. The scalar-microlocal construction first gives the inclusion (\Cref{prop:tautprimitive})
\[
 T_{\dd}:=R_{\dd}\cap[\FM_{\Pi_Q,\dd}]^{\rd}
       \subseteq\widehat{\Fg}_{\Pi_Q,\dd}.
\]
The image of $T_{\dd}$ under the flag correspondence is $\kappa(\Delta_{\dd}R_{\dd})\cap[\Y]$, where $\kappa:P_{\dd}\to\HO^*(\Y)$ is the tautological map, sending $x_{i,r}$ to the first Chern class of $\CL_{i,r}$ (\Cref{prop:vfc}). Here $[\Y]$ denotes the classical fundamental class of the smooth irreducible stack $\Y$. Smoothness of $\Y$ implies that it coincides with the reduced virtual fundamental class (\Cref{lem:stablefund}). A dimension-one vertex makes the smooth target a framed quiver stack $\Y$ which is a $\mathbf{G}_{\mathrm{m}}$-gerbe over its good moduli space $Y$, and $Y$ is a quiver variety. Note that if a vertex with dimension $1$ has loops, the quiver variety is a product of an affine space with the quiver variety for the quiver with the loops at the vertex erased and so one may erase the loops for the sake of tautological generation; this comes from the fact that in dimension one, the preprojective relation, defined using commutators, does not involve loops. The stacky version of the McGerty--Nevins theorem on Kirwan surjectivity (\Cref{prop:Kirwan}) therefore implies that the map $\kappa$ is surjective.

The Weyl group action is constructed in an algebraic family (\Cref{subsection:springer-action}), using the formalism of Dobrovolska--Ginzburg--Travkin \cite{DGT}. We develop the part of this formalism that we need here and provide all proofs. We keep a stability condition $\Theta$ fixed, where $\Theta$ satisfies a positivity condition on vertices at legs of $Q^+$ (given by \Cref{lem:chamber}) and vary the central moment parameter of $Q^+$. The full legs identify its parameter space with the space of ordered eigenvalues for the original moment map (\Cref{eq:eigenparameters}). The resulting smooth semiprojective family $\mathcal{X}$ has constant cohomology, by \cite[Theorem~10.2.2]{DGT}, see \Cref{prop:constant}. Its projective map $\pi$ to the ordered spectral base change is small (\Cref{prop:small}), so the generic permutation action extends through intersection cohomology (\Cref{subsection:springer-action}). Specializing gives the Springer action on the cohomology $\HO^*(Y,\BoQ)$ of the moduli space $Y$ of $\mathfrak{Y}$ (\Cref{subsection:springer-action}). All flag quotient line bundles between tautological vector bundles extend over the family $\mathcal{X}$. Consequently this action permutes the $x_{i,r}$, including after passage to the scalar gerbe (\Cref{prop:Wequiv}), since the gerbe is $W$-invariant.

The original vector $\dd\in\BoN^{Q_0}$ may be divisible.  Contrasting with this, the dimension vector $\dd^+$ on the extended quiver $Q^+$ always has at least one entry equal to one, providing a smooth family after choosing a generic stability condition.

Every alternating polynomial in $P_{\dd}$ is a Vandermonde determinant \eqref{equation:vandermonde} times a symmetric polynomial. It follows from surjectivity of $\kappa$ (see \Cref{equation:flag-map} and \Cref{prop:Kirwan}) that
\[
 \Phi(T_{\dd})=\HO^*(\Y)^{\sgn}\cap[\Y].
\]
The HLRV formula (\Cref{thm:HLRV}) identifies the graded dimension of this image with the graded dimension of $\widehat{\Fg}_{\Pi_Q,\dd}$. Since each graded piece is finite-dimensional, the inclusion of $T_{\dd}$ in the affinized BPS cohomology $\widehat{\Fg}_{\Pi_Q,\dd}$ must be an equality. Toda's natural stability comparison for the BPS sheaf/BPS cohomology (\Cref{thm:Toda}) then proves \Cref{thm:main} for arbitrary $\theta$ by restriction.

We also give an independent Lie-theoretic construction of the sign summand of $\HO^*(Y,\BoQ)$ for the $W$-action (\Cref{sec:Lie}). It uses polynomial representations of the Kac--Moody Lie algebras attached to the legs (whose product is isomorphic to the product $\prod_{i\in Q_0}\mathfrak{sl}_{d_i}$) and the identification of BPS Lie algebras with positive halves of generalised Kac--Moody Lie algebras \cite{DHSM}. The generation proof does not require the identification of this Weyl group action with the geometric action on the whole target coming from Springer theory and actually, we do not provide such a comparison. However, comparing these two actions is an interesting question that we plan to study further in the future.

\subsection{Contents of the paper}
The main results (\Cref{thm:main} and \Cref{thm:signintro}) concern the undeformed preprojective algebra $\Pi_Q$. In \Cref{sec:deformation}, we prove the geometric statements for central deformations $\Pi_Q^{\lambda}$, $\lambda\in\BoC^{Q_0}$, using mild generalisations of our methods, see \Cref{conj:deformed}.

We recall the results of cohomological Donaldson--Thomas (DT) theory of preprojective algebras that we use: the cohomological integrality theorem (\Cref{thm:PBW}), the factorisation coproduct (\Cref{subsection:factorisation-coproduct}), the scalar-microlocal Thom isomorphism (\Cref{thm:scalar}), Toda's stability theorem (\Cref{thm:Toda}), Kirwan surjectivity for quiver varieties (\Cref{prop:Kirwan}) following McGerty--Nevins, and the HLRV formula for the sign-part of the cohomology of quiver varieties with legs (\Cref{thm:HLRV} and \Cref{cor:signcharacter}). We also use the simultaneous resolution, the smallness criterion, and the algebraic constancy theorem of Dobrovolska--Ginzburg--Travkin \cite{DGT} for a family of quiver varieties, see \Cref{sec:family}. For the sake of self-containedness, we include the proof of the constancy (\Cref{prop:constant}) after \cite{DGT}. We expect this paper to be crucial for comparisons of BPS cohomologies in various contexts.

\subsection*{Acknowledgements}
The second author would like to thank Tasuki Kinjo for explaining to him virtual fundamental classes of $0$-shifted symplectic stacks at a workshop in Seoul in June 2025. Substantial progress on this project was made by the authors while they were visiting SLMath in Spring 2024, at a conference at IPMU in March 2025, at the CARE conference in Lyon in October 2025, and during the opportunity of a visit of the second author to Edinburgh in April 2026. The second author shared the main  theorem of this work and the ideas underlying the proof in private discussions with Olivier Schiffmann.

\section{Preprojective stacks and tautological classes}\label{sec:notation}
\subsection{Quivers and moment maps}
\label{subsection:quivers-moment-maps}
Let $Q=(Q_0,Q_1)$ be a quiver with set of vertices $Q_0$ and set of arrows $Q_1$. Let $s,t:Q_1\to Q_0$ be the source and target maps. For each arrow $a$, let $a^*$ denote the opposite arrow. The doubled quiver $\overline Q=(Q_0,Q_1\sqcup Q_1^{\mathrm{op}})$ has arrows $a,a^*$ for $a\in Q_1$ where $Q_1^{\mathrm{op}}$ denotes the set of opposite arrows, and we define the preprojective algebra of $Q$ to be the quotient of the path algebra of $\overline{Q}$ by the preprojective relation:
\[
 \Pi_Q=\BoC\overline Q\Big/\Big(\sum_{a\in Q_1}[a,a^*]\Big).
\]
For dimension vectors $\dd,\ee\in\BoN^{Q_0}$, the Euler form and its symmetrized version are given by
\[
 \chi_Q(\dd,\ee)=\sum_{i\in Q_0} d_i e_i-\sum_{a\in Q_1}d_{s(a)}e_{t(a)}, \qquad (\dd,\ee)_Q=\chi_Q(\dd,\ee)+\chi_Q(\ee,\dd).
\]
We let $V_i=\BoC^{d_i}$ be a complex vector space of dimension $d_i$ for $i\in Q_0$ and $G_{\dd}\coloneqq\prod_{i\in Q_0}\mathrm{GL}_{d_i}$. We define the representation space of $\dd$-dimensional representations of $\overline{Q}$ and the Lie algebra of $G_{\dd}$:
\[
 V_{\dd}=\Rep_{\dd}(\overline Q) =\bigoplus_{a\in Q_1}\bigl(\Hom(V_{s(a)},V_{t(a)})\oplus\Hom(V_{t(a)},V_{s(a)})\bigr),\qquad \mathfrak{gl}_{\dd}=\bigoplus_i\End(V_i).
\]
The trace pairing identifies $\mathfrak{gl}_{\dd}$ with its dual $\mathfrak{gl}_{\dd}^*$ and $V_{\dd}$ with the cotangent space $T^*\Rep_{\dd}(Q)$ of the representation space of $\dd$-dimensional representations of the original quiver $Q$. The moment map $\mu\colon V_{\dd}\rightarrow \mathfrak{gl}_{\dd}$, $(x_a,x_{a^*})_{a\in Q_1}\mapsto\sum_{a\in Q_1}[x_a,x_{a^*}]$ has components
\[
 \mu_i(x)=
 \sum_{a:t(a)=i}x_ax_{a^*}
 -
 \sum_{a:s(a)=i}x_{a^*}x_a
\]
at the $i$th vertex for $i\in Q_0$. In particular $\sum_{i\in Q_0}\tr\mu_i(x)=0$ by invariance of the trace under cyclic permutations of its summands. Write
\[
 \mathfrak{gl}_{\dd}^{\,0}
 =\ker\Bigl(\mathfrak{gl}_{\dd}\xrightarrow{\sum_i\tr}\BoC\Bigr),
 \qquad X_{\dd}=\mu^{-1}(0),
 \qquad \FM_{\Pi_Q,\dd}=[X_{\dd}/G_{\dd}]
\]
for the trace-free Lie algebra, the representation space of $\dd$-dimensional representations of the preprojective algebra $\Pi_Q$ and the moduli stack of $\dd$-dimensional representations of $\Pi_Q$. Vertices $i\in Q_0$ of dimension zero, i.e. with $d_i=0$, may be deleted without changing the result of the constructions of this paper. When convenient, we may therefore assume that $d_i>0$ for every $i\in Q_0$. Note that $\FM_{\Pi_Q,\dd}$ has (unreduced) virtual dimension $-(\dd,\dd)_Q=\dim V_{\dd}-\dim G_{\dd}-\dim \mathfrak{gl}_{\dd}$ by a direct calculation.

\subsection{Stability}
\label{subsection:generic}
Let $\theta\in\mathbf{R}^{Q_0}$ be a stability parameter and let
\[
 \Lambda_\theta
 :=\{\dd\in\BoN^{Q_0}:\theta\cdot\dd=0\}
\]
be the corresponding monoid of dimension vectors, where $\theta\cdot\dd\coloneqq\sum_{i\in Q_0}\theta_id_i$. Whenever $\theta$ is fixed, we consider only dimension vectors in $\Lambda_\theta$. We use the convention that a representation $M$ of dimension $\dd\in\Lambda_\theta$ is \emph{$\theta$-semistable} if
\[
 \theta\cdot\dim U\leq0
\]
for every subrepresentation $U\subset M$, where $\dim U=(\dim U_i)_{i\in Q_0}\in\BoN^{Q_0}$ denotes the dimension vector of $U$. It is \emph{stable} if the inequality is strict for every nonzero proper subrepresentation. This is King's convention for stability conditions \cite{King}. A stability parameter $\theta$ is called \emph{generic} for $\dd$ if $\theta\cdot \dd=0$ and
\[
 \theta\cdot\ee\ne0\quad\text{for all }0<\ee<\dd.
\]
Inequalities between dimension vectors are coordinatewise, and $\ee<\dd$ means $\ee\leq\dd$ and $\ee\ne\dd$. It is easily seen that a generic stability parameter exists if and only if $\dd$ is indivisible (i.e. its coordinates are globally coprime). For a fixed $\dd\in\Lambda_\theta$, choose a rational parameter $\theta'\in\BoQ^{Q_0}$ satisfying
\[
 \sgn(\theta'\cdot\ee)=\sgn(\theta\cdot\ee)
 \qquad(0\leq\ee\leq\dd).
\]
Such a parameter exists because the finitely many rational hyperplanes $\{\theta' \in\mathbf{R}^{Q_0}\;\lvert\; \theta'\cdot \ee\}\subset \mathbf{R}^{Q_0}$, for $0< \ee< \dd$ cut out rational polyhedral faces. One of them contains $\theta$ and contains points with rational coordinates.  We may choose $\theta'$ to be one of these points. In particular, $\theta'\cdot\dd=0$, and the stable and semistable loci in dimension $\dd$ are unchanged. Multiplying $\theta'$ by a positive integer gives an integral parameter for GIT. Thus the constructions of moduli spaces below apply for real parameters as well by performing such a replacement one component at a time.

We write $X_{\dd}^{\theta}\subset X_{\dd}$ for the open subset of semistable points and $\FM_{\Pi_Q,\dd}^{\theta}=[X_{\dd}^{\theta}/G_{\dd}]$ for the corresponding moduli stack. Then, $\FM_{\Pi_Q,\dd}^{\theta}$ is the stack of $\theta$-semistable $\dd$-dimensional representations of the preprojective algebra. Its good moduli space is denoted by $\CM_{\Pi_Q,\dd}^{\theta}$. We use the same superscript convention for representations of the tripled quiver: the superscript $\theta$ always denotes the semistable locus for the King stability parameter $\theta$. The zero stability parameter is omitted from the notation, i.e. when $\theta=0$, $\CM_{\Pi_Q,\dd}=\CM_{\Pi_Q,\dd}^{\theta}$.

\subsection{Cohomology and Borel--Moore homology}
\label{subsection:coh-BM}
We give the Borel--Moore homology the natural topological degree: a smooth complex connected variety $X$ of dimension $m$ has a fundamental class $[X]$ in degree $2m$: $[X]\in\HO^{\mathrm{BM}}_{2m}(X,\BoQ)$. The same convention applies to smooth quotient stacks. Thus, since $\BG\Gm$ has dimension $-1$,
\[
 [\BG\Gm]\in\HO_{-2}^{\BM}(\BG\Gm,\BoQ).
\]
For a smooth stack of pure complex dimension $m$, cap product with the fundamental class gives an isomorphism
\[
 \HO^a(\mathfrak X,\BoQ)
 \xrightarrow{\ \sim\ }
 \HO_{2m-a}^{\BM}(\mathfrak X,\BoQ), \quad a\in\BoZ.
\]
This is Poincar\'e duality. For quotient stacks, which include all stacks that this paper is concerned with, these statements can be interpreted by equivariant approximations via the Borel construction \cite{Borel1960TransformationGroups}, see \cite[\S1.1]{anderson2012introduction} and one can use the equivariant derived categories of constructible sheaves developed in \cite{bernstein2006equivariant}.

We put
\begin{equation}\label{eq:cohgrading}
 \CA_{\Pi_Q,\dd}^k
   =\HO_{-k-(\dd,\dd)_Q}^{\BM}(\FM_{\Pi_Q,\dd},\BoQ)=\HO^{k+(\dd,\dd)_Q}(\FM_{\Pi_Q,\dd},\mathbb{D}\BoQ).
\end{equation}
Writing $s=(\dd,\dd)_Q$ and $\mathbb D$ for Verdier duality, we have
\[
 \BoQ_{\FM_{\Pi_Q,\dd}}^{\vir} :=\BoQ_{\FM_{\Pi_Q,\dd}}[-s], \qquad \CA_{\Pi_Q,\dd}^k =\HO^k\!\left(\FM_{\Pi_Q,\dd},\mathbb D\BoQ_{\FM_{\Pi_Q,\dd}}^{\vir}\right).
\]
Here $-s$ is the derived virtual dimension of $\FM_{\Pi_Q,\dd}$, and $\mathbb D\BoQ^{\vir}=(\mathbb D\BoQ)[s]$. The same convention on the semistable open substack gives $\CA_{\Pi_Q,\dd}^{\theta}=\HO^*\!\left(\FM^{\theta}_{\Pi_Q,\dd},\mathbb D\BoQ_{\FM^{\theta}_{\Pi_Q,\dd}}^{\vir}\right)$ for $\dd\in\Lambda_\theta$, and $\vir$ indicates the same shift since the derived virtual dimension of a stack and a non-empty open substack coincide. A class $\alpha\in\HO^{2a}(\FM_{\Pi_Q,\dd}^{\theta},\BoQ)$ acts on $\CA_{\Pi_Q,\dd}^{\theta}$ with degree $2a$; cap product with $\alpha$ gives for any $k\in \BoZ$ a morphism
\[
\CA_{\Pi_Q,\dd}^{\theta,k}\rightarrow\CA_{\Pi_Q,\dd}^{\theta,k+2a}.
\]
We state and prove the tautological generation theorem for rational graded vector spaces (see \Cref{thm:signintro}), that is we work with constructible sheaves with rational coefficients throughout the paper. It will be clear that all constructions are compatible with the natural mixed Hodge structures involved (which are all known to all be pure of Tate type, \cite[Corollary~6.2]{DavInt}).

\subsection{The tautological ring}
Let $\CV_i$ be the universal vector bundle, of rank $d_i$, on the representation stack corresponding to the $i$th vertex of $Q$. Its Chern classes define
\[
 R_{\dd}=\BoQ[c_{i,r}:i\in Q_0,\ 1\leq r\leq d_i]
       =\HO^*(\BG G_{\dd},\BoQ),
 \qquad \deg c_{i,r}=2r.
\]
For the full stack $\FM_{\Pi_Q,\dd}$ of $\dd$-dimensional representations of the preprojective algebra $\Pi_Q$, this is its entire ordinary cohomology as shown by the following lemma, which follows from the equivariant contractibility of $X_{\dd}$.

\begin{lemma}\label{lem:contraction}
The pullback in cohomology by the projection $\FM_{\Pi_Q,\dd}\to\BG G_{\dd}$ induces an isomorphism
\[
 R_{\dd}\xrightarrow{\sim}\HO^*(\FM_{\Pi_Q,\dd},\BoQ).
\]
\end{lemma}
\begin{proof}
Multiplication of every arrow of $\overline{Q}$ by a real number $t\in[0,1]$ defines a $G_{\dd}$-equivariant deformation retraction of $X_{\dd}$ onto the zero representation. The moment map is homogeneous of degree two, so the homotopy stays in $X_{\dd}$: $X_{\dd}$ is $G_{\dd}$-equivariantly contractible. This proves the statement.
\end{proof}

\subsection{The scalar gerbe}\label{subsec:gerbe}
The diagonal subgroup $\Gm\subset G_{\dd}$ acts trivially on $X_{\dd}$. The automorphism group of a $\theta$-stable representation of $\Pi_Q$ coincides with this scalar subgroup. Hence if $\dd$ and $\theta$ are such that $\theta\cdot\dd=0$ and every $\theta$-semistable representation of $\Pi_Q$ of dimension $\dd$ is stable, then the stack $\FM_{\Pi_Q,\dd}^{\theta}$ is a $\Gm$-gerbe over its coarse moduli space $\CM_{\Pi_Q,\dd}^{\theta}$.

For $|\dd|=\sum_i d_i>0$, we set
\begin{equation}\label{eq:scalaru}
 u_{\dd}=\frac1{|\dd|}\sum_i c_1(\CV_i)\in\HO^2(\FM_{\Pi_Q,\dd}^{\theta},\BoQ).
\end{equation}
There is an isomorphism
\begin{equation}\label{eq:gerbesplitting}
 \HO^*(\FM_{\Pi_Q,\dd}^{\theta},\BoQ)\cong \HO^*(\CM_{\Pi_Q,\dd}^{\theta},\BoQ)\otimes\BoQ[u_{\dd}]
\end{equation}
when semistability equals stability for $\theta$.

One may describe the tautological classes on the coarse moduli space $\CM_{\Pi_Q,\dd}^{\theta}$ when stability and semistability coincide, using the quotient $PG_{\dd}=G_{\dd}/\Gm$, where $\Gm$ denotes the diagonal scalar subgroup. Indeed, the representation space of stable representations $X_{\dd}^{\theta\text{-st}}$ is a principal $PG_{\dd}$-bundle over its coarse quotient $\CM_{\Pi_Q,\dd}^{\theta\text{-st}}$, which means that there is a morphism $\CM_{\Pi_Q,\dd}^{\theta\text{-st}}\rightarrow\BG PG_{\dd}$ to the classifying space $\BG PG_{\dd}$ determining this principal bundle. This morphism defines a characteristic-class map from $\HO^*(\BG PG_{\dd},\BoQ)$ to the cohomology of that quotient $\HO^*(\CM_{\Pi_Q,\dd}^{\theta\text{-st}},\BoQ)$.

To describe these classes, recall that
\[
 R_{\dd}=\HO^*(\BG G_{\dd},\BoQ)=\BoQ[x_{i,r}]^{\prod_i\mathfrak S_{d_i}},\qquad u_{\dd} =\frac{1}{|\dd|}\sum_i c_1(\CV_i)=\frac{1}{|\dd|}\sum_{i,r}x_{i,r}, \qquad |\dd|=\sum_i d_i.
\]
Pullback along the morphism $\BG G_{\dd}\to \BG PG_{\dd}$ induced by the natural quotient $G_{\dd}\rightarrow PG_{\dd}$ identifies $H^*(\BG PG_{\dd},\BoQ)$ with the subring of $R_{\dd}$ of polynomials in the $x_{i,r}$ invariant under simultaneous translation of all variables
\[
 x_{i,r}\longmapsto x_{i,r}+t.
\]
Since $u_{\dd}\mapsto u_{\dd}+t$ if all $x_{i,r}$ are translated by $t$, this subring consists precisely of the polynomials in the $x_{i,r}-u_{\dd}$ that are symmetric in the $x_{i,r}$ separately at each vertex.

\section{CoHAs and affinized BPS cohomology}\label{sec:DT}
\subsection{The tripled quiver}
\label{subsection:triple-quiver}
Let $\widetilde Q$ be the quiver obtained from the doubled quiver $\overline Q$ by adding a loop $\omega_i$ at every vertex. Thus $\tilde{Q}$ has set of vertices $Q_0$ and set of arrows $Q_1\sqcup Q_1^{\mathrm{op}}\sqcup\{\omega_i\colon i\in Q_0\}$. We call $\widetilde{Q}$ the \emph{tripled quiver}. The canonical cubic potential is
\[
 W=\sum_{a\in Q_1}\bigl(\omega_{t(a)}aa^*-\omega_{s(a)}a^*a\bigr)=\left(\sum_{i\in Q_0}\omega_i\right)\left(\sum_{a\in Q_1}[a,a^*]\right)\quad \in\BoC\widetilde{Q}/[\BoC\widetilde{Q},\BoC\widetilde{Q}].
\]
Its trace on representations with dimension vector $\dd$ is
\[
 \Tr(W)_{\dd}(x,\omega)=\sum_{i\in Q_0}\tr(\omega_i\mu_i(x)).
\]
The critical equations given by the partial derivatives $\frac{\partial W}{\partial a}$, $\frac{\partial W}{\partial a^*}$ and $\frac{\partial W}{\partial \omega_i}$ say that $\mu(x)=0$ and that the $\omega_i$ commute with all arrow maps $a,a^*$. Thus, the Jacobi algebra of $(\widetilde{Q},W)$ is
\begin{equation}\label{equation:jacobi-algebra}
 \operatorname{Jac}(\widetilde Q,W)\coloneqq\BoC\widetilde{Q}\big/ \left(\frac{\partial W}{\partial a},\frac{\partial W}{\partial a^*},\frac{\partial W}{\partial \omega}\right)\cong\Pi_Q[\omega],
\end{equation}
where $\omega$ is an extra indeterminate. In particular a critical point of $\Tr(W)$ on the stack of representations of $\widetilde{Q}$ is a representation of the preprojective algebra $\Pi_Q$ equipped with an endomorphism (determined by the action of $\omega$).

The dimensional-reduction isomorphism identifies the critical CoHA of $(\widetilde Q,W)$ \cite{DM} with the preprojective CoHA structure on $\CA_{\Pi_Q}=\bigoplus_{\dd}\CA_{\Pi_Q,\dd}$. We use the dimensional-reduction construction of \cite{DavCritical,DavInt}. See also \cite{kinjo2022dimensional} for a more general dimensional-reduction theorem in the framework of shifted cotangent stacks. The CoHA product for the cohomological Hall algebra of $(\widetilde{Q},W)$ and $\CA_{\Pi_Q}$ are not strictly the same: for classes of dimension vectors $\dd$ and $\ee$, the untwisted product is multiplied by $(-1)^{\chi_Q(\dd,\ee)}$ \cite[Appendix]{ren2017cohomological}. Dimensional reduction is taken with the componentwise sign correction of \cite[\S8.2]{DHSM}, so that it respects these products defined at the level of quivers with potentials as in \cite{DM} and at the level of preprojective algebras \cite{yang2018cohomological}. The factorisation coproduct (see \Cref{subsection:factorisation-coproduct}) carries the compatible twist of \cite[\S6.1]{JinThesis}. We refer the reader to these references for the precise conventions, as they do not change any of the subspaces or cap products occurring in the generation theorem.

\subsection{The BPS cohomology subspace and the cohomological integrality isomorphism}
\label{subspace:BPS-subspace}
We briefly recall the embedding of the BPS cohomology inside the Borel--Moore homology of the preprojective stack. Fix a stability parameter $\theta\in\BoQ^{Q_0}$ and a nonzero dimension vector $\dd\in\Lambda_\theta$, and let
\[
 \JH^\theta_{\dd}\colon\FM^\theta_{\widetilde Q,\dd}\longrightarrow \CM^\theta_{\widetilde Q,\dd}
\]
be the good moduli space morphism for the stack of $\dd$-dimensional $\theta$-semistable representations of the tripled quiver $\widetilde{Q}$. On $\BoC$-points, the map $\JH_{\dd}^{\theta}$ sends a $\theta$-semistable representation to the associated $\theta$-polystable representation, obtained from a Jordan--Hölder filtration in the abelian category of $\theta$-semistable representations with dimension vectors in $\Lambda_\theta$.

We write $\Tr(W)^\theta_{\dd}$ for the restriction of the trace of the potential to the stack of semistable points $\FM_{\widetilde{Q},\dd}^{\theta}$, and define
\[
 \mathcal K^\theta_{\dd}:=\JH^\theta_{\dd,*}{}^p\varphi_{\Tr(W)^\theta_{\dd}} \BoQ_{\FM^\theta_{\widetilde Q,\dd}}^{\vir}.
\]
Direct images of constructible complexes are derived by convention. On the smooth stack of representations of the tripled quiver, $\BoQ^{\vir}$ denotes the constant sheaf shifted by the complex dimension of that stack. That is, $\BoQ_{\FM_{\widetilde{Q},\dd}^{\theta}}^{\vir}=\BoQ_{\FM_{\widetilde{Q},\dd}^{\theta}}[\dim\FM_{\widetilde{Q},\dd}^{\theta}]$. This is the virtual normalisation compatible with \eqref{eq:cohgrading} under dimensional reduction. All perverse truncations that we consider in the sequel are taken on the good moduli space $\CM^\theta_{\widetilde Q,\dd}$.

Since $\widetilde Q$ is symmetric, the genericity hypothesis on the stability parameter in the cohomological integrality theorem in \cite{DM} (the symmetry of the Euler form of $\widetilde{Q}$ restricted to $\Lambda_{\theta}$) is satisfied for every $\theta$. In particular, we do not assume that semistability coincides with stability. \cite{DM} gives the vanishing
\[
 {}^p\mathcal H^j(\mathcal K^\theta_{\dd})=0 \qquad\text{for }j<1.
\]
The BPS sheaf is the perverse sheaf
\[
 \calBPS^\theta_{\widetilde Q,W,\dd} := {}^p\mathcal H^1(\mathcal K^\theta_{\dd}) \cong \bigl({}^p\tau^{\leq1}\mathcal K^\theta_{\dd}\bigr)[1].
\]
This is the lowest possibly non-vanishing perverse cohomology sheaf of $\mathcal{K}_{\dd}^{\theta}$. The canonical perverse truncation adjunction morphism therefore gives
\[
 \epsilon^\theta_{\dd}\colon \calBPS^\theta_{\widetilde Q,W,\dd}[-1] \cong {}^p\tau^{\leq1}\mathcal K^\theta_{\dd} \longrightarrow \mathcal K^\theta_{\dd}.
\]
The maps it induces on hypercohomology are injective by \cite[Theorem~C]{DM}, where it is explained that $\epsilon^\theta_{\dd}$ admits a retraction, and therefore so does $\HO^*(\epsilon^\theta_{\dd})$, giving the injectivity. This injectivity is part of the cohomological integrality theorem and is a consequence of the fact that $\mathcal{K}_{\dd}^{\theta}$ splits as the direct sum of its shifted perverse cohomologies. This is true for any symmetric quiver with potential, under the genericity assumption for $\theta$, even though the perverse cohomologies may not be semisimple perverse sheaves/mixed Hodge modules in general. However, note that they are in fact semisimple in our situation of the tripled quiver with its canonical cubic potential, see \cite{DavBPS}. See also \cite[Theorem~8.1 and \S10.1]{DHSM}.

We define the BPS cohomology by applying the derived global sections functor:
\[
 \Fg^\theta_{\Pi_Q,\dd}:=\HO^*\!\left(\CM^\theta_{\widetilde Q,\dd},   \calBPS^\theta_{\widetilde Q,W,\dd}[-1]\right) =\HO^{*-1}\!\left(   \CM^\theta_{\widetilde Q,\dd},\calBPS^\theta_{\widetilde Q,W,\dd}\right).
\]

We briefly explain why dimensional reduction \cite[Appendix~A]{DavCritical} applies with this stability condition: the explanation is needed since the natural morphism $F\colon\FM_{\widetilde{Q},\dd}\rightarrow\FM_{\overline{Q},\dd}$ forgetting the extra loops does not map $\FM_{\widetilde{Q},\dd}^{\theta}$ to $\FM_{\overline{Q},\dd}^{\theta}$ in general. What we explain is that $\mathrm{Crit}(\Tr(W)_{\dd}^{\theta})$ is contained in $F^{-1}(\FM_{\overline{Q},\dd}^{\theta})$. Recall that there is an identification $\operatorname{Jac}(\widetilde Q,W) \cong \Pi_Q[\omega]$ \eqref{equation:jacobi-algebra}. Thus a representation in the critical locus of $\Tr(W)_{\dd}$ is a pair $(M,\omega_M)$ consisting of a $\dd$-dimensional $\Pi_Q$-module $M$ and a $\Pi_Q$-linear endomorphism $\omega_M$. Such a pair is $\theta$-semistable if and only if $M$ is $\theta$-semistable. Indeed, one implication follows because every subrepresentation of the pair is in particular a submodule of $M$. For the converse, if $M$ is unstable, we let $N\subset M$ be a destabilizing subrepresentation. We may take for $N$ the first Harder--Narasimhan submodule (i.e. the maximal submodule of $M$ of maximal slope). Then, the ratio $\theta\cdot\dim N/|\dim N|$ is positive, and $N$ is preserved by every endomorphism of $M$, hence by $\omega_M$, and destabilises the pair.

It follows that the critical locus of $\Tr(W)_{\dd}^{\theta}$ in $\FM^\theta_{\widetilde Q,\dd}$ is contained in the open substack on which the underlying $\overline Q$-representation is $\theta$-semistable, i.e. in $F^{-1}(\FM_{\overline{Q},\dd}^{\theta})$. On this open substack, forgetting the additional loops is a vector bundle projection over $\FM_{\overline{Q},\dd}^{\theta}$, and $\Tr(W)_{\dd}^{\theta}$ is linear along its fibres. Since the vanishing-cycle complex is supported on the critical locus, dimensional reduction on this open substack computes the hypercohomology of $\mathcal K^\theta_{\dd}$; see \cite[\S3.2.2--\S3.2.3]{DavBPS}. With the normalisations of cohomological degrees fixed above, it gives a graded isomorphism
\begin{equation}\label{equation:dimensional-reduction}
 \operatorname{DR}^\theta_{\dd}\colon\HO^*\!\left(   \CM^\theta_{\widetilde Q,\dd},   \mathcal K^\theta_{\dd} \right) \xrightarrow{\sim} \CA^\theta_{\Pi_Q,\dd},
\end{equation}
where $\CA^\theta_{\Pi_Q,\dd}$ denotes the normalised Borel--Moore homology of $\FM^\theta_{\Pi_Q,\dd}$.

We therefore obtain a canonical embedding of the BPS cohomology inside the Borel--Moore homology of the semistable preprojective stack
\begin{equation}
\label{equation:embeddingBPS}
 \iota^\theta_{\dd}\colon \Fg^\theta_{\Pi_Q,\dd} \xhookrightarrow{\ \HO^*(\epsilon^\theta_{\dd})\ } \HO^*\!\left(\CM^\theta_{\widetilde Q,\dd},\mathcal K^\theta_{\dd}\right) \xrightarrow[\sim]{\ \operatorname{DR}^\theta_{\dd}\ } \CA^\theta_{\Pi_Q,\dd}.
\end{equation}

The category of $\theta$-semistable $\Pi_Q$-modules with dimension vectors in $\Lambda_\theta$ is closed under extensions. Its cohomological Hall algebra is
\[
 \CA^\theta_{\Pi_Q}:=\bigoplus_{\dd\in\Lambda_\theta}\CA^\theta_{\Pi_Q,\dd}, \qquad \CA^\theta_{\Pi_Q,0}:=\BoQ.
\]
With the compatible sign conventions for the Hall products, the subspace
\[
 \Fg^\theta_{\Pi_Q}:=\bigoplus_{\dd\in\Lambda_\theta\setminus\{0\}} \Fg^\theta_{\Pi_Q,\dd}\subseteq \CA^\theta_{\Pi_Q}
\]
is closed under the graded commutator Lie bracket, by \cite[Corollary~6.11]{DM}; see also \cite[Theorem~8.1 and \S10.1]{DHSM}. Equipped with this bracket, it is the BPS Lie algebra of $\theta$-semistable $\Pi_Q$-modules.

On a nonzero component $\FM^\theta_{\Pi_Q,\dd}$, put
\[
 \delta_{\dd}:=c_1\!\left(\det\!\left(\bigoplus_{i\in Q_0}\CV_i\right)\right)=\sum_i c_1(\CV_i)=|\dd|u_{\dd}.
\]
The determinant-class action is the action of $\BoQ[u]$ on $\CA_{\Pi_Q}^{\theta}$, with $\deg(u)=2$, given by $u\cdot\alpha=\delta_{\dd}\cap\alpha$. On the dimension-zero component we have $u\cdot1=0$. We write $\iota^\theta\colon \Fg_{\Pi_Q}^{\theta}\rightarrow\CA_{\Pi_Q}^{\theta}$ for the direct sum of the embeddings $\iota^\theta_{\dd}$ for $\dd\in\Lambda_{\theta}$.

We may now recall the cohomological integrality isomorphism for preprojective algebras at the level of graded vector spaces \cite{DavInt}. We do not need the sheaf enhancement of this isomorphism of \emph{loc.cit.} at the level of constructible sheaves over $\CM_{\Pi_Q,\dd}^{\theta}$.
\begin{theorem}[Cohomological integrality]\label{thm:PBW}
For every stability parameter $\theta$, the map
\[
 I^\theta\colon\Fg^\theta_{\Pi_Q}\otimes\BoQ[u]\longrightarrow \CA^\theta_{\Pi_Q},\qquad x\otimes u^k\longmapsto u^k\cdot\iota^\theta(x),
\]
is injective. Symmetrised Hall multiplication induces an isomorphism
\[
 \operatorname{PBW}^\theta\colon\Sym\!\left(\Fg^\theta_{\Pi_Q}\otimes\BoQ[u] \right)\xrightarrow{\sim} \CA^\theta_{\Pi_Q}
\]
of vector spaces graded by the cohomological normalised degree and the dimension vector.
\end{theorem}
Here $\Sym$ denotes the graded symmetric algebra, with the usual Koszul signs, and the map is obtained by applying $I^\theta$ to each factor and averaging the iterated Hall products with these signs. It is not an isomorphism of algebras in general. It is analogous to the PBW isomorphism relating the symmetric algebra over a Lie algebra and its enveloping algebra.
\begin{proof}
Apply \cite[Theorem~C]{DM} to the symmetric quiver with potential $(\widetilde Q,W)$ and its category of $\theta$-semistable representations with dimension vectors in $\Lambda_\theta$. Dimensional reduction identifies the resulting semistable critical CoHA with $\CA^\theta_{\Pi_Q}$, compatibly with the determinant-class action and, with our sign conventions, the Hall multiplication. See also \cite{DavInt}.
\end{proof}

Taking $\theta=0$ recovers the construction for all representations and the full stack $\FM_{\Pi_Q}$, without stability condition.

\Cref{thm:PBW} is the PBW form of the cohomological integrality theorem \cite{DM,DavBPS}. The preprojective formulation and its relation with the BPS Lie algebra are used and recalled in \cite{DHSM,Jin}. The formal variable $u$ in the Hall algebra acts by the additive class $\delta_{\dd}$, whereas $u_{\dd}$ in \eqref{eq:scalaru} is the normalised gerbe class. On a fixed nonzero dimension component, replacing $\delta_{\dd}$ by $u_{\dd}$ only rescales the polynomial generator and leaves the image of the affinized BPS space in $\CA_{\Pi_Q,\dd}^{\theta}$ unchanged. Across dimension vectors, it is $\delta_{\dd}$ that is additive under direct sum and so we use this normalisation in the coproduct calculation below.

\subsection{The factorisation coproduct}
\label{subsection:factorisation-coproduct}

Fix a stability parameter $\theta$, and recall the notation $\Lambda_\theta$ introduced above in \Cref{subsection:generic}. We recall the construction of the factorisation coproduct on $\CA^\theta_{\Pi_Q}$ following \cite{JinThesis,DHKSV}. Throughout, the Hall product and the coproduct carry the compatible sign twists fixed above.

For an analytic open subset $U\subset\BoC$, let $\FM^\theta_{\widetilde Q,\dd}(U)$ be the open substack on which every eigenvalue of every additional loop $\omega_i$ belongs to $U$ and let $\imath_U\colon\FM^\theta_{\widetilde Q,\dd}(U)\rightarrow\FM^\theta_{\widetilde Q,\dd}$ be the (analytic) open inclusion. Denote by $\mathcal H^\theta(U)=\bigoplus_{\dd\in\Lambda_{\theta}}\HO^*(\FM_{\widetilde{Q},\dd}^{\theta}(U),\imath_U^*{^p\varphi}_{\Tr(W)_{\dd}^{\theta}}\BoQ_{\FM_{\widetilde{Q},\dd}^{\theta}}^{\mathrm{vir}})$ its normalised critical cohomology, summed over $\dd\in\Lambda_\theta$. Dimensional reduction identifies $\mathcal H^\theta(\BoC)\cong \CA^\theta_{\Pi_Q}$.

On the critical locus, a representation is a pair $(M,\omega)$, where $M$ is a $\Pi_Q$-module and $\omega$ is a $\Pi_Q$-linear endomorphism. For two disjoint open subsets $U_1,U_2\subset\BoC$, the condition $\operatorname{Spec}(\omega)\subset U_1\sqcup U_2$ gives a canonical spectral decomposition
\[
 (M,\omega)= (M_1,\omega_1)\oplus(M_2,\omega_2), \qquad \operatorname{Spec}(\omega_\nu)\subset U_\nu\quad \nu=1,2,
\]
The summands are $\Pi_Q$-submodules because $\omega$ commutes with the arrow maps.

This decomposition preserves $\theta$-semistability. Indeed, if $M$ is semistable, then
\[
 \theta\cdot\dim M_\nu\leq0, \quad\nu=1,2, \qquad \theta\cdot\dim M_1+\theta\cdot\dim M_2=\theta\cdot\dim M=0.
\]
Thus $\theta\cdot\dim M_\nu=0$ for both summands. Every submodule of $M_\nu$ is a submodule of $M$, so $M_\nu$ is itself semistable. Conversely, a direct sum of $\theta$-semistable modules is $\theta$-semistable, since this category is closed under extensions. In particular, both spectral summands have dimension vectors in $\Lambda_\theta$.

The spectral decomposition is exact. In particular, it also decomposes short exact sequences of critical representations. The factorisation construction therefore applies to the category of $\theta$-semistable objects and to its Hall induction correspondences.

For a nonempty disk $U\subset\BoC$, open restriction induces an isomorphism
\[
 \rho_U\colon \mathcal H^\theta(\BoC) \xrightarrow{\sim}\mathcal H^\theta(U).
\]
To see this, we use the support theorem for the semistable tripled BPS sheaf: it is supported on the locus where all eigenvalues coincide, and is constant in the direction of this common eigenvalue.  This follows either by applying the same argument as \cite[Theorem/Definition~4.1]{DavBPS}, or by combining this theorem with the existence of $\Xi$ as in Theorem \ref{thm:Toda}, or as a very special case of \cite[Theorem 7.2.15]{bu2025cohomology}. The relative cohomological integrality theorem for the $\theta$-semistable category (the version of \Cref{thm:PBW} for $(\widetilde{Q},W)$ at the level of constructible sheaves over $\CM_{\widetilde{Q}}^{\theta}$, see \cite{DM}) consequently reduces the assertion to the isomorphism
\[
 \HO^*(\BoC,\BoQ)\xrightarrow{\sim}\HO^*(U,\BoQ).
\]
Alternatively, one can use the translation invariance of the vanishing cycle sheaf ${^p\varphi_{\Tr(W)_{\dd}^{\theta}}}\BoQ_{\FM_{\widetilde{Q},\dd}^{\theta}}^{\mathrm{vir}}$ for the translation of the $\omega$ coordinate and its $\mathbf{G}_{\mathrm{m}}$ invariance for the action rescaling the $\omega$ coordinate (see \cite{DHKSV}).

The restriction $\rho_{U}$ is compatible with the determinant-class action and the Hall multiplication. We still denote by $\rho_{U}$ the restriction even when $U$ is not a disk (but it is then not necessarily an isomorphism).

For two disjoint disks $U_1,U_2$, spectral factorisation and the Thom--Sebastiani theorem give an isomorphism
\[
 F^\theta_{U_1,U_2}\colon\mathcal H^\theta(U_1\sqcup U_2) \xrightarrow{\sim} \mathcal H^\theta(U_1)\otimes\mathcal H^\theta(U_2).
\]

Using dimensional reduction to identify $\mathcal H^\theta(\BoC)$ with
$\CA^\theta_{\Pi_Q}$, we define \[
 \Delta^\theta_{\mathrm{fac}} :=(\rho_{U_1}^{-1}\otimes\rho_{U_2}^{-1}) \circ F^\theta_{U_1,U_2} \circ\rho_{U_1\sqcup U_2}.
\]
Thus, the factorisation coproduct is the morphism
\[
 \Delta^\theta_{\mathrm{fac}}\colon \CA^\theta_{\Pi_Q} \longrightarrow \CA^\theta_{\Pi_Q}\otimes \CA^\theta_{\Pi_Q}.
\]
Its component on dimension vector $\dd$ takes values in
\[
 \bigoplus_{\substack{\dd_1+\dd_2=\dd\\\dd_1,\dd_2\in\Lambda_\theta}} \CA^\theta_{\Pi_Q,\dd_1}\otimes\CA^\theta_{\Pi_Q,\dd_2}.
\]

The coproduct is independent of the chosen disks, coassociative, graded cocommutative, and compatible with the Hall product. The counit is the projection to the zeroth graded piece with respect to the dimension vector degree. These statements follow from the construction of \cite[\S9]{Jin}, with detailed proofs in \cite[\S6.1]{JinThesis}: the spectral decompositions, restriction maps, and Hall correspondences used there restrict to the semistable categories considered here. We denote by $\Prim(\CA_{\Pi_Q}^{\theta},\Delta_{\mathrm{fac}}^{\theta})\subset\CA_{\Pi_Q}^{\theta}$ the subspace of primitive elements.

\begin{proposition}\label{prop:primitives}
For every stability parameter $\theta\in\mathbf{R}^{Q_0}$, there is an identification of the subspace of primitive elements of $\CA_{\Pi_Q}^{\theta}$ with the affinized BPS cohomology $\widehat{g}_{\Pi_Q}^{\theta}$, defined as the image of $I^{\theta}$ (\Cref{thm:PBW}):
\[
 \Prim\!\left(\CA^\theta_{\Pi_Q},\Delta^\theta_{\mathrm{fac}}\right) = \widehat{\Fg}^{\theta}_{\Pi_Q} := \im\!\left(\Fg^\theta_{\Pi_Q}\otimes\BoQ[u]\xrightarrow{\,I^\theta\,}\CA^\theta_{\Pi_Q}\right).
\]
\end{proposition}

\begin{proof}
Write $\Delta=\Delta^\theta_{\mathrm{fac}}$. The support lemma \cite[Lemma~4.1]{DavInt} places the spectral support of $\calBPS^\theta_{\widetilde Q,W,\dd}$ on the small diagonal, so $\mathrm{supp}(\calBPS_{\widetilde Q,W,\dd})\subset \CM_{\overline{Q},\dd}^{\theta}\times\mathbf{A}^1\subset \CM_{\widetilde{Q},\dd}^{\theta}$ where the copy of $\mathbf{A}^1$ prescribes a single eigenvalue for all extra loops $\omega_i$, $i\in Q_0$. Its restriction to a locus with two nonzero spectral summands therefore vanishes. Consequently, for every $x\in\Fg^\theta_{\Pi_Q}$,
\[
 \Delta(\iota^\theta(x)) = \iota^\theta(x)\otimes1 +1\otimes\iota^\theta(x).
\]

The determinant-class action is a coderivation:
\[
 \Delta(u\cdot a)
 =
 (u\otimes1+1\otimes u)\cdot\Delta(a).
\]
Indeed, under a direct sum the total tautological bundle splits as the sum of the tautological bundles on the two factors, so its first Chern class is the sum of their first Chern classes. Since $u\cdot1=0$, every $u^k\cdot\iota^\theta(x)$ is primitive. This proves
\[
 V:=\widehat{\Fg}^{\theta}_{\Pi_Q} \subseteq \Prim(\CA^\theta_{\Pi_Q},\Delta).
\]

Equip $\Sym(V)$ with the standard graded cocommutative coalgebra structure for which every element of $V$ is primitive. Since $\Delta$ is compatible with multiplication and the inclusion of $V$ has primitive image, symmetrised Hall multiplication defines a coalgebra morphism
\[
 \Sym(V)\longrightarrow\CA^\theta_{\Pi_Q}.
\]
By \Cref{thm:PBW}, this morphism is an isomorphism of graded vector spaces, hence an isomorphism of coalgebras. Over $\BoQ$, the primitive subspace of $\Sym(V)$ is exactly $V$. The asserted equality follows. Note that we did not need to use the Milnor--Moore theorem for cocommutative connected bialgebras to obtain this identification.
\end{proof}

The primitive subspace $\mathrm{Prim}(\CA_{\Pi_Q}^{\theta},\Delta_{\mathrm{fac}}^{\theta})$ is closed under the graded commutator bracket: this is true for any cocommutative bialgebra, by a direct calculation. It therefore gives the affinized BPS subspace $\widehat{\Fg}^{\theta}_{\Pi_Q}$ the structure of a Lie subalgebra of $\CA^\theta_{\Pi_Q}$. This characterisation is intrinsic to the factorisation coproduct and requires no choice of a PBW splitting (i.e. a splitting of the PBW map $\mathrm{PBW}^{\theta}$ in \Cref{thm:PBW}).

\subsection{Graded dimensions}
Let $A_{Q,\dd}(q)$ be the Kac polynomial counting absolutely indecomposable representations over finite fields \cite{kac2006root,MR0677715,MR0607162,MR0557581}. The BPS/Kac identity is
\begin{equation}\label{eq:BPSKac}
 \sum_k\dim\Fg_{\Pi_Q,\dd}^k\,t^{k/2}=A_{Q,\dd}(t^{-1}).
\end{equation}
See \cite{DavInt,DavBPS} and \cite[\S9]{DHSM}. In particular the BPS cohomology is finite-dimensional and concentrated in (non-positive) even degrees. We obtain
\begin{equation}\label{eq:affKac}
 \sum_k\dim\widehat{\Fg}_{\Pi_Q,\dd}^{\,k}t^{k/2}=\frac{A_{Q,\dd}(t^{-1})}{1-t}
\end{equation}
since $\widehat{\Fg}_{\Pi_Q,\dd}\cong \Fg_{\Pi_Q,\dd}\otimes\BoQ[u]$ and $u$ sits in cohomological degree $2$. If $A_{Q,\dd}\ne0$, it is monic of degree $p_Q(\dd)$ \cite[\S1.15, Proposition]{kac2006root}. Therefore the smallest cohomological degree in $\Fg_{\Pi_Q,\dd}$, and in its affinization $\widehat{\Fg}_{\Pi_Q,\dd}$, is $-2p_Q(\dd)$, with multiplicity one.

\subsection{Stability conditions and BPS cohomology}\label{subsec:Toda}
The tripled quiver is symmetric. Toda's theorem \cite{Toda}, which we recall here for a general symmetric quiver with potential, therefore applies to its canonical cubic potential.

\begin{theorem}[Toda]\label{thm:Toda}
For a symmetric quiver with potential $(Q,W)$, a stability parameter $\theta\in\mathbf{R}^{Q_0}$, and a dimension vector $\dd\in\BoN^{Q_0}$ such that $\theta\cdot\dd=0$, the natural adjunction morphism induces an isomorphism between its BPS sheaf on $\CM_{Q,\dd}$ and the direct image of the BPS sheaf on $\CM_{Q,\dd}^{\theta}$ under the GIT affinization map. Namely, there is a morphism $\pi\colon\CM_{Q,\dd}^{\theta}\rightarrow\CM_{Q,\dd}$ and $\Xi\colon \pi_*\calBPS_{Q,W,\dd}^{\theta}\cong \calBPS_{Q,W,\dd}$. For the tripled quiver with the canonical cubic potential, this isomorphism together with the dimensional reduction gives an identification of the BPS spaces
\[
 r_\theta:\Fg_{\Pi_Q,\dd}\xrightarrow{\sim}\Fg_{\Pi_Q,\dd}^{\theta} \quad\text{if }\theta\cdot\dd=0.
\]
\end{theorem}

The assertion valid for any symmetric quiver with potential is \cite[Lemma~4.10]{Toda}; see also \cite[Proposition~10.1]{DHSM} and \cite[Proposition~6.3]{Jin}. The morphism is induced by the adjunction map associated with the open embedding $\imath\colon\FM_{\tilde{Q},\dd}^{\theta}\rightarrow\FM_{\tilde{Q},\dd}$. Indeed, we have $\imath^*{}^p\varphi_{\Tr(W)_{\dd}} \BoQ_{\FM_{\widetilde Q,\dd}}^{\vir}\cong {}^p\varphi_{\Tr(W)^\theta_{\dd}} \BoQ_{\FM^\theta_{\widetilde Q,\dd}}^{\vir}$ by open restriction of vanishing cycle sheaf complexes.
Toda's isomorphism $r_{\theta}$ is obtained by applying ${}^p\tau^{\leq 1}\JH_*$ to the adjunction ${}^p\varphi_{\Tr(W)_{\dd}} \BoQ_{\FM_{\widetilde Q,\dd}}^{\vir}\rightarrow \imath_*{}^p\varphi_{\Tr(W)^\theta_{\dd}} \BoQ_{\FM^\theta_{\widetilde Q,\dd}}^{\vir}$ and taking derived global sections, using the isomorphism $\pi_*{}^p\tau^{\leq 1}\JH^{\theta}_*\varphi_{\Tr(W)^\theta_{\dd}} \BoQ_{\FM^\theta_{\widetilde Q,\dd}}^{\vir}\cong {}^p\tau^{\leq 1}\JH_*\varphi_{\Tr(W)_{\dd}} \BoQ_{\FM_{\widetilde Q,\dd}}^{\vir}$; see \cite{Toda} for details.

\begin{corollary}\label{cor:Todaaff}
Let $\theta\in\mathbf{R}^{Q_0}$ be a stability parameter and $\dd\in\BoN^{Q_0}$ such that $\theta\cdot\dd=0$. Open restriction in Borel--Moore homology $\HO^{\mathrm{BM}}(\FM_{\Pi_Q,\dd})\rightarrow\HO^{\mathrm{BM}}(\FM_{\Pi_Q,\dd}^{\theta})$ induces an isomorphism
\[
 r_\theta:\widehat{\Fg}_{\Pi_Q,\dd}\xrightarrow{\sim}\widehat{\Fg}_{\Pi_Q,\dd}^{\theta}.
\]
It commutes with the action of every tautological class for which that action stays in the indicated affinized BPS cohomology subspace.
\end{corollary}
The precaution in the last statement is justified by the fact we do not know a priori that the affinized BPS cohomology is preserved by the action of tautological classes. This fact will be a consequence of \Cref{thm:main}, which is proven in \Cref{sec:proof}.
\begin{proof}
The determinant line bundles restrict to the same line bundles on the open $\FM_{\Pi_Q,\dd}^{\theta}\subset\FM_{\Pi_Q,\dd}$. Affinizing \Cref{thm:Toda} with this action gives the isomorphism. Alternatively, spectral factorisation commutes with restriction to the $\theta$-semistable locus: a direct summand of a $\theta$-semistable object is $\theta$-semistable and has dimension vector in $\Lambda_\theta$. The cap-product compatibility is the usual compatibility with open restriction.
\end{proof}

\section{Reduced virtual fundamental classes}\label{sec:virtual}
\subsection{The reduced derived zero locus}
\label{subsection:reduced-derived-zero-locus}
Fix $\dd\neq0$. The equation $\sum_i\tr\mu_i=0$ holds identically on $V_{\dd}$. Thus the moment map takes values in the $G_{\dd}$-representation
\[
 \mathfrak{gl}_{\dd}^{\,0} := \ker\left(\bigoplus_{i\in Q_0}\mathfrak{gl}_{d_i}   \xrightarrow{\ \sum_i\tr\ }\BoC\right),
\]
of dimension $\dim G_{\dd}-1$. Replacing the target $\mathfrak{gl}_{\dd}$ of the moment map by $\mathfrak{gl}_{\dd}^{0}$ removes the identically zero total-trace equation, without changing the classical zero locus.

\begin{definition}
\label{definition:reduced-stack}
The reduced derived preprojective stack is the derived fiber product
\[
 \FM_{\Pi_Q,\dd}^{\rd}:=\bigl[   (V_{\dd}\times^{\mathbf R}_{\mathfrak{gl}_{\dd}^{\,0}}\{0\})/G_{\dd}\bigr].
\]
Its classical truncation is $\FM_{\Pi_Q,\dd}$, and its virtual dimension is
\[
 \nu_{\dd} := \dim V_{\dd} -\dim\mathfrak{gl}_{\dd}^{\,0} -\dim G_{\dd} =\dim V_{\dd}-2\dim G_{\dd}+1.
\]
Its equivariant virtual fundamental class is denoted by
\[
 v_{\dd} :=[\FM_{\Pi_Q,\dd}]^{\rd} \in H^{\mathrm{BM}}_{2\nu_{\dd}}       (\FM_{\Pi_Q,\dd},\BoQ),
\]
where we use unshifted Borel--Moore homological degrees here. We may call $[\FM_{\Pi_Q,\dd}]^{\rd}$ the \emph{reduced virtual fundamental class} of $\FM_{\Pi_Q,\dd}$. For $\theta\cdot\dd=0$, we denote by $v_{\dd}^{\theta}$ its restriction to the $\theta$-semistable open for the natural restriction morphism $H^{\mathrm{BM}}_{2\nu_{\dd}}(\FM_{\Pi_Q,\dd},\BoQ)\rightarrow H^{\mathrm{BM}}_{2\nu_{\dd}}(\FM_{\Pi_Q,\dd}^{\theta},\BoQ)$. Equivalently, $v_{\dd}^{\theta}$ is the reduced virtual fundamental class of $\FM_{\Pi_Q,\dd}^{\theta}$.
\end{definition}

We explain more explicitly the definition of the reduced virtual fundamental class, fix $\dd\ne0$ and put
\[
 \mathcal X_{\dd}:=[V_{\dd}/G_{\dd}],\qquad\mathcal E_{\dd}:=[V_{\dd}\times\mathfrak{gl}_{\dd}^{\,0}/G_{\dd}]\longrightarrow\mathcal X_{\dd}.
\]
The stack $\mathcal X_{\dd}$ is smooth, and $\mathcal E_{\dd}$ is a vector bundle over $\mathcal{X}_{\dd}$ of rank
\[
 r_{\dd} :=\dim\mathfrak{gl}_{\dd}^{\,0}=\sum_i d_i^2-1.
\]
The moment map defines a section $s_{\mu}$ of this vector bundle. Writing $\iota_{\dd}\colon\FM_{\Pi_Q,\dd}\hookrightarrow\mathcal X_{\dd}$ for the classical zero locus of $s_{\mu}$, we have the Cartesian square
\[
\begin{array}{ccc}
 \FM_{\Pi_Q,\dd}
   &\xrightarrow{\ \iota_{\dd}\ }&\mathcal X_{\dd}\\
 {\scriptstyle\iota_{\dd}}\big\downarrow
   &&\big\downarrow{\scriptstyle s_{\mu}}\\
 \mathcal X_{\dd}
   &\xrightarrow{\ 0_{\mathcal E_{\dd}}\ }&\mathcal E_{\dd}.
\end{array}
\]
Replacing this fibre product by the derived fibre product gives $\FM_{\Pi_Q,\dd}^{\rd}$ defined in \Cref{definition:reduced-stack}.

The zero section $0_{\mathcal{E}_{\dd}}\colon \mathcal{X}_{\dd}\rightarrow\mathcal{E}_{\dd}$ is a regular embedding of codimension $r_{\dd}$. Its refined Gysin homomorphism, formed with respect to the above square, is
\[
 0_{\mathcal E_{\dd},s_{\mu}}^{!}\colon A_k(\mathcal X_{\dd})_{\BoQ} \longrightarrow A_{k-r_{\dd}}(\FM_{\Pi_Q,\dd})_{\BoQ}.
\]
Here the Chow groups of the quotient stacks are understood equivariantly \cite{EG98}. We define
\[
 v_{\dd} := 0_{\mathcal E_{\dd},s_{\mu}}^{!}[\mathcal X_{\dd}] \in A_{\nu_{\dd}}(\FM_{\Pi_Q,\dd})_{\BoQ}, \qquad \nu_{\dd} =\dim V_{\dd}-\dim G_{\dd}-r_{\dd}.
\]
We use the same notation for its cycle class in Borel--Moore homology, with the normalisations fixed above.

Concretely, by \cite{fulton2013intersection}, the section $s_{\mu}$ induces a closed embedding of the normal cone
\[
 j_{\dd}\colon
 C_{\FM_{\Pi_Q,\dd}/\mathcal X_{\dd}}
 \hookrightarrow
 \mathcal E_{\dd}|_{\FM_{\Pi_Q,\dd}}.
\]
Deformation to the normal cone specialises $[\mathcal X_{\dd}]$ to the fundamental cycle of this cone. The refined pullback is then the ordinary zero-section Gysin pullback in the restricted vector bundle:
\[
 v_{\dd}= 0_{\mathcal E_{\dd}|_{\FM_{\Pi_Q,\dd}}}^{!}
 \,j_{\dd,*} [C_{\FM_{\Pi_Q,\dd}/\mathcal X_{\dd}}].
\]
Thus $v_{\dd}$ is a class on the zero locus $\FM_{\Pi_Q,\dd}$ itself. Its pushforward to the ambient stack $\mathcal{X}_{\dd}$ satisfies
\[
 \iota_{\dd,*}v_{\dd}=c_{r_{\dd}}(\mathcal E_{\dd})
 \cap[\mathcal X_{\dd}].
\]
If the reduced moment-map equations form a regular sequence, the normal cone is the restricted vector bundle and $v_{\dd}$ agrees with the ordinary fundamental cycle $[\FM_{\Pi_Q,\dd}]$ of $\FM_{\Pi_Q,\dd}$.

The adjective ``reduced'' refers to the fact that we remove the identically vanishing total-trace equation from the derived zero-locus construction, i.e. we replace $\mathfrak{gl}_{\dd}$ by $\mathfrak{gl}_{\dd}^0$ for the target of the moment map. This removes one trivial equation and increases the virtual dimension by one.  The word ``reduced'' here has no relation with the operation of taking the quotient by nilpotent elements in rings of functions on schemes.
\begin{lemma}\label{lem:vdim}
One has
\[
 \vdim\FM_{\Pi_Q,\dd}^{\rd}=1-(\dd,\dd)_Q, \qquad v_{\dd}\in \HO_{2-2(\dd,\dd)_Q}^{\BM}(\FM_{\Pi_Q,\dd},\BoQ).
\]
Its degree in $\CA_{\Pi_Q,\dd}$ is $(\dd,\dd)_Q-2=-2p_Q(\dd)$.
\end{lemma}
\begin{proof}
The dimension of $V_{\dd}$ is $2\sum_{a\in Q_1}d_{s(a)}d_{t(a)}$. The rank of the equation bundle is $\dim\mathfrak{gl}_{\dd}^0=\sum_i d_i^2-1$, and quotienting by $G_{\dd}$ subtracts another $\sum_i d_i^2$. This gives the virtual dimension. The other assertions follow from the topological Borel--Moore homological degree (\Cref{subsection:coh-BM}) and the shift \eqref{eq:cohgrading}.
\end{proof}

\begin{remark}
Using all of $\mathfrak{gl}_{\dd}$ instead as the target of the moment map would add a zero section of a trivial line bundle. Its Euler class is zero. Thus the ordinary unreduced virtual fundamental class vanishes.

\end{remark}

\begin{lemma}\label{lem:stablefund}
On the open substack of $\FM_{\Pi_Q,\dd}$ of stable representations, the reduced virtual fundamental class restricts to the ordinary fundamental class of $\FM_{\Pi_Q,\dd}^{\theta\text{-}\mathrm{st}}$.

If the stack $\FM_{\Pi_Q,\dd}$ contains stable representations of some dimension vector $\dd\in\Lambda^{\theta}$ for some stability condition $\theta$, then the ordinary fundamental class of $\FM_{\Pi_Q,\dd}^{\theta}$ coincides with the reduced virtual fundamental class.
\end{lemma}
\begin{proof}
The differential of the Hamiltonian moment map $\mu\colon V_{\dd}\rightarrow\mathfrak{gl}_{\dd}$ is dual, through the symplectic form, to the infinitesimal action of $\mathfrak{gl}_{\dd}$ on $V_{\dd}$. A stable representation  has infinitesimal stabiliser equal to the diagonal scalar line $\BoC\id\subset\mathfrak{gl}_{\dd}$. Consequently the differential is surjective onto $\mathfrak{gl}_{\dd}^{\,0}$, the annihilator of that line. The reduced equations are regular there, and over the stable open substack, the reduced derived structure is the ordinary smooth structure.

The second statement follows from the first and the fact that under the existence of stable representations, the classical stack $\FM_{\Pi_Q,\dd}^{\theta}$ is irreducible of dimension $1-(\dd,\dd)_Q$, see \cite{CB}, which is also the reduced virtual dimension.
\end{proof}

\subsection{Scalar microlocal homology}

We explain how the reduced virtual class arises from the critical description of $\Pi_Q[\omega]$ (\Cref{subsection:triple-quiver}). The argument is the global moment-map version of \cite[Proposition~2.13]{DHKSV}. We recall here the zero-section comparison used in its proof.

We fix a dimension vector $\dd\neq0$ and a stability parameter $\theta\in\BoQ^{Q_0}$ satisfying $\theta\cdot\dd=0$ and put
\[
 \CF_{\dd}^{\theta} :={}^p\varphi_{\Tr(W)_{\dd}^{\theta}} \BoQ_{\FM_{\widetilde Q,\dd}^{\theta}}^{\vir}.
\]

For an analytic open $U\subset\BoC$, consider the spectral open $\FM_{\widetilde Q,\dd}^{\theta}(U)$ introduced in \Cref{subsection:factorisation-coproduct}.
The scalar locus
\[
 \Lambda_{\mathrm{scal},\dd}^{\theta}(U):= \{(M,z\,\id_M):     M\in\FM_{\Pi_Q,\dd}^{\theta},\ z\in U\}\cong\FM_{\Pi_Q,\dd}^{\theta}\times U
\]
is closed in the critical locus of $\Tr(W)_{\dd}^{\theta}$ over this open. Let
\[
 \iota_U: \Lambda_{\mathrm{scal},\dd}^{\theta}(U) \hookrightarrow\FM_{\widetilde Q,\dd}^{\theta}(U)
\]
be its inclusion. We define
\[
 C_{\mathrm{scal},\dd}^{\theta}(U) := \HO^*\!\left(   \Lambda_{\mathrm{scal},\dd}^{\theta}(U),\iota_U^!\left(\CF_{\dd}^{\theta}|_{\FM_{\widetilde Q,\dd}^{\theta}(U)}\right)\right).
\]
Since $\iota_U$ is closed, the adjunction morphism $\iota_{U,*}\iota_U^!\to\id$ gives a map
\begin{equation}\label{equation:forgetting-supports}
C_{\mathrm{scal},\dd}^{\theta}(U)\rightarrow\mathcal{H}^{\theta}(U)
\end{equation}
 to the critical cohomology over $U$.

Write $\mathcal H_{\dd}^{\theta}(U)$ for the dimension-$\dd$ summand of $\mathcal H^{\theta}(U)$. For a nonempty disk $U$, forgetting supports \eqref{equation:forgetting-supports}, the disk-restriction isomorphism $\rho_U$ (\Cref{subsection:factorisation-coproduct}), and dimensional reduction \eqref{equation:dimensional-reduction} give the composite
\[
 \epsilon_U: C_{\mathrm{scal},\dd}^{\theta}(U)\longrightarrow \mathcal H_{\dd}^{\theta}(U) \xrightarrow{\ \rho_U^{-1}\ } \mathcal H_{\dd}^{\theta}(\BoC) \xrightarrow{\ \operatorname{DR}_{\dd}^{\theta}\ } \CA_{\Pi_Q,\dd}^{\theta}.
\]

\begin{theorem}[Scalar-microlocal Thom isomorphism]
\label{thm:scalar}
For every nonempty disk $U$, there is a natural isomorphism
\begin{equation}\label{equation:thomiso}
 \operatorname{Th}_U: \HO^a(\FM_{\Pi_Q,\dd}^{\theta},\BoQ) \xrightarrow{\sim} C_{\mathrm{scal},\dd}^{\theta}(U)^{a-2p_Q(\dd)}.
\end{equation}
Under this identification,
\[
 \epsilon_U\bigl(\operatorname{Th}_U(\alpha)\bigr) =\alpha\cap v_{\dd}^{\theta}\qquad\text{for any $\alpha\in\HO^*(\FM_{\Pi_Q,\dd}^{\theta},\BoQ)$}.
\]
In particular, the unit $1\in\HO^*(\FM_{\Pi_Q,\dd}^{\theta},\BoQ)$ maps to $v_{\dd}^{\theta}$. These identifications commute with restriction to smaller disks $U'\subset U$ and to semistable opens.
\end{theorem}
\begin{proof}
Let $B$ be an algebraic variety, $\CE\to B$ be a vector bundle of rank $r$, and let $s\colon B\rightarrow\CE$ be a section with classical zero locus $Z\coloneqq s^{-1}(0)$. Consider the fibrewise linear function given by the contraction with $s$
\[
 h:\operatorname{Tot}(\CE^\vee)\longrightarrow\BoC, \qquad h(b,\xi)=\langle s(b),\xi\rangle.
\]
Put $\nu=\dim B-r$ (the virtual dimension of $s^{-1}(0)$), and let
\[
 \CF_h := {}^p\varphi_h \BoQ_{\operatorname{Tot}(\CE^\vee)}^{\vir}, \qquad\BoQ_{\operatorname{Tot}(\CE^\vee)}^{\vir}\coloneqq\BoQ_{\operatorname{Tot}(\CE^\vee)}[\dim\CE].
\]
We regard $\CF_h$ as a constructible complex on the critical locus $\operatorname{Crit}(h)$ of $h$. Denote by
\[
 j:Z= \operatorname{Crit}(h)\cap B=s^{-1}(0)\hookrightarrow\operatorname{Crit}(h),\qquad q:\operatorname{Crit}(h)\longrightarrow Z
\]
the zero section and the projection. Then, by dimensional reduction,
\[
 q_!\CF_h\cong\BoQ_Z[\nu], \qquad q_*\CF_h\cong\mathbb D\BoQ_Z[-\nu].
\]
Under these identifications, the functorial morphism $q_!\rightarrow q_*$ forgetting supports gives the morphism
\begin{equation}\label{equation:map-forgetting-supports}
 \BoQ_Z[\nu]\longrightarrow\mathbb D\BoQ_Z[-\nu]
\end{equation}
corresponding to the virtual fundamental class of the derived zero locus of $s$. In particular, the induced map
\[
 \HO^a(Z,\BoQ) \longrightarrow\HO_{2\nu-a}^{\BM}(Z,\BoQ)
\]
is the cap product with the virtual fundamental class $0_{\CE}^![B]$ of $s^{-1}(0)$, defined via the Cartesian square
\[\begin{tikzcd}
	{s^{-1}(0)} & B \\
	B & {\mathcal{E}}
	\arrow[from=1-1, to=1-2]
	\arrow[from=1-1, to=2-1]
	\arrow["\lrcorner"{anchor=center, pos=0.125}, draw=none, from=1-1, to=2-2]
	\arrow["s"', from=1-2, to=2-2]
	\arrow["{0_{\mathcal{E}}}", from=2-1, to=2-2]
\end{tikzcd}\]

We refer to \cite[\S2.2.3 and Corollary~2.8]{DHKSV} for more details. These operations are compatible with smooth pullback, so the same construction applies equivariantly and on smooth quotient stacks.

We apply this comparison to the reduced moment-map section. The trace pairing gives a $G_{\dd}$-equivariant decomposition
\[
 \mathfrak{gl}_{\dd} = \BoC\cdot\id_{\dd}\oplus\mathfrak{gl}_{\dd}^{\,0}.
\]
Accordingly, we write
\[
 \omega=z\,\id_{\dd}+\eta,\qquad z=\frac{1}{|\dd|}\sum_i\tr(\omega_i),\qquad\eta\in\mathfrak{gl}_{\dd}^{\,0}.
\]
The identity $\sum_i\tr\mu_i=0$ gives
\[
 \Tr(W)_{\dd}(x,\omega) =\sum_i\tr\bigl(\eta_i\mu_i(x)\bigr).
\]
Thus the trace of the potential $\Tr(W)_{\dd}$ is independent of the scalar coordinate $z$.

Let $V_{\dd}^{\theta}\subset V_{\dd}$ denote the open subspace of the space of $\dd$-dimensional representations of $\overline{Q}$ containing the semistable representations. The induced fibre-linear function is associated with the section $\mu$ of the vector bundle
\[
 \CE_{\dd}^{0,\theta} := [V_{\dd}^{\theta}\times\mathfrak{gl}_{\dd}^{\,0}/G_{\dd}]\longrightarrow[V_{\dd}^{\theta}/G_{\dd}].
\]
Its derived zero locus is the $\theta$-semistable open of $\FM_{\Pi_Q,\dd}^{\theta,\rd}$.

By \Cref{lem:vdim}, the reduced virtual dimension is
\[
 \nu_{\dd} := \vdim\FM_{\Pi_Q,\dd}^{\rd} = 1-(\dd,\dd)_Q = 2p_Q(\dd)-1.
\]
The scalar locus is the product $\FM_{\Pi_Q,\dd}^{\theta}\times\mathbf{A}^1$ of the zero section for the reduced fibre-linear function with the $z$-line. Smooth pullback for normalised vanishing cycles and the zero-section comparison therefore give
\[
 \iota_U^!\left(\CF_{\dd}^{\theta}|_{\FM_{\widetilde Q,\dd}^{\theta}(U)} \right) \cong \BoQ_{\FM_{\Pi_Q,\dd}^{\theta}}[\nu_{\dd}]      \boxtimes\BoQ_U[1] \cong \BoQ_{\Lambda_{\mathrm{scal},\dd}^{\theta}(U)}       [2p_Q(\dd)].
\]
Since $\HO^*(U,\BoQ)=\BoQ$ for a disk, taking hypercohomology proves the isomorphism in \Cref{equation:thomiso}.

To identify the corresponding map forgetting supports \eqref{equation:map-forgetting-supports}, first take the whole scalar line $U=\BoC$. The factor $\BoQ_{\BoC}[1]$, which is Verdier self-dual, contributes the same shift $[1]$ to its source and target. After dimensional reduction, the zero-section comparison identifies the map with the virtual-class morphism
\[
 \BoQ_{\FM_{\Pi_Q,\dd}^{\theta}}[\nu_{\dd}+1] \longrightarrow \mathbb D\BoQ_{\FM_{\Pi_Q,\dd}^{\theta}}[1-\nu_{\dd}].
\]
It therefore sends the unit to the refined zero-section class $v_{\dd}^{\theta}$ defined in \S \ref{subsection:reduced-derived-zero-locus}. The construction is $\HO^*(\FM_{\Pi_Q,\dd}^{\theta},\BoQ)$-linear, so it sends every $\alpha$ to $\alpha\cap v_{\dd}^{\theta}$.

Restriction to a disk $U\subset\BoC$ commutes with the zero-section comparison and with forgetting supports. Using $\rho_U$ gives the same formula for $\epsilon_U$. All constructions also commute with restriction along open embeddings.
\end{proof}

\begin{proposition}\label{prop:tautprimitive}
For every nonzero dimension vector,
\[
 \HO^*(\FM^{\theta}_{\Pi_Q,\dd},\BoQ)\cap v_{\dd}\subseteq\widehat{\Fg}^{\theta}_{\Pi_Q,\dd}.
\]
That is, the tautological subspace $\HO^*(\FM^{\theta}_{\Pi_Q,\dd},\BoQ)\cap v_{\dd}$ of $\CA^{\theta}_{\Pi_Q}$ is contained in the primitive subspace or equivalently the affinized BPS cohomology.
\end{proposition}
\begin{proof}
Let $a\in\HO^*(\FM^{\theta}_{\Pi_Q,\dd},\BoQ)$ and put $x=a\cap v_{\dd}$. Choose disjoint disks $U_1,U_2$ inside a larger disk $U$. By \Cref{thm:scalar}, the critical-cohomology class corresponding to $x$ under dimensional reduction lies in the image of $C^{\theta}_{\mathrm{scal},\dd}(U)$, namely cohomology with supports on $\Lambda^{\theta}_{\mathrm{scal}}(U)$.

Consider a decomposition $\dd=\dd_1+\dd_2$ with
$\dd_1,\dd_2\ne0$. The component
$\Delta^{\theta}_{\mathrm{fac},\dd_1,\dd_2}$ is constructed by restriction
to the spectral locus of critical representations
\[
 (M,\omega)
 =
 (M_1,\omega_1)\oplus(M_2,\omega_2),
 \qquad
 \dim M_\nu=\dd_\nu,\quad
 \operatorname{Spec}(\omega_\nu)\subset U_\nu,\quad \nu=1,2
\]
followed by spectral factorisation and Thom--Sebastiani.

This locus is disjoint from $\Lambda^{\theta}_{\mathrm{scal}}(U)$: a scalar endomorphism $z\,\id_M$ has only the eigenvalue $z$, whereas both nonzero summands above have nonempty spectra contained in disjoint disks. The map from cohomology with supports on $\Lambda^{\theta}_{\mathrm{scal}}(U)$ therefore becomes zero after restriction to this spectral locus. Compatibility with dimensional reduction gives
\[
 \Delta^{\theta}_{\mathrm{fac},\dd_1,\dd_2}(x)=0 \qquad\text{if }\dd_1,\dd_2\ne0.
\]

Only the components indexed by $(\dd,0)$ and $(0,\dd)$ remain. Under the disk-restriction identifications used to define the coproduct, these are respectively $x\otimes1$ and $1\otimes x$. Hence
\[
 \Delta^{\theta}_{\mathrm{fac}}(x)=x\otimes1+1\otimes x.
\]
Thus $x$ is primitive, and \Cref{prop:primitives} implies $x\in\widehat{\Fg}^{\theta}_{\Pi_Q,\dd}$.
\end{proof}

The same support proof is used in \cite[\S7.5]{DHKSV} for sheaves on symplectic surfaces. It applies to every cohomology multiple, without requiring that the tautological ring preserve the primitive subspace a priori.

\begin{lemma}\label{lem:restrictionvirtual}
Let $\theta\in\mathbf{R}^{Q_0}$ be a stability parameter. For the map in Borel--Moore homology induced by the restriction to the semistable locus $r_{\theta}\colon \HO^{\mathrm{BM}}(\FM_{\Pi_Q,\dd})\rightarrow\HO^{\mathrm{BM}}(\FM_{\Pi_Q,\dd}^{\theta})$, we have
\[
 r_\theta(a\cap v_{\dd})=a|_{\FM_{\Pi_Q,\dd}^{\theta}}\cap v_{\dd}^{\theta}.
\]
\end{lemma}
\begin{proof}
Refined zero-section pullback defining the virtual class, \Cref{subsection:reduced-derived-zero-locus}, is easily seen to commute with open restriction. The cap product has the same compatibility. By combining both statements, we obtain the lemma.
\end{proof}

\section{Quivers with legs}\label{sec:legs}
\subsection{The extended dimension vector}
Fix a quiver $Q$ with dimension vector $\dd\in\mathbf{N}^{Q_0}$.  For $n=d_i>0$, introduce vertices $(i,r)$, $1\leq r<n$, and arrows
\[
 (i,n-1)\longrightarrow\cdots\longrightarrow(i,1)                 \longrightarrow(i,0)=i.
\]
Give $(i,r)$ dimension $n-r$. There are $n-1$ added vertices. The chain has $n$ vertices when the original vertex $i=(i,0)$ is included.  We perform this operation at every $i\in Q_0$.  We denote the resulting data consisting of the extended quiver with the extended dimension vector by $(Q^+,\dd^+)$. Although it is not apparent from the notation, the extended quiver $Q^+$ depends on $\dd$.  

If $\dd\neq 0$, then $\dd^+$ has a coordinate equal to $1$, and is therefore indivisible. Note that if all coordinates of $\dd$ are $0$ or $1$, then $Q^+=Q$ and so the vertex with coordinate $1$ may carry loops. This does not cause any problems as noted in \Cref{subsection:strategy-proof}: we may discard the loops at $1$-dimensional vertices without changing the cohomology. Put $N=N_{\dd}$ and $W=W_{\dd}$ defined as in \eqref{equation:weyl-group}.

\begin{lemma}\label{lem:legdimension}
One has
\[
 (\dd^+,\dd^+)_{Q^+}=(\dd,\dd)_Q-2N,\qquad p_{Q^+}(\dd^+)=p_Q(\dd)+N.
\]
\end{lemma}
\begin{proof}
This is a direct calculation. The contribution of the additional vertices at a leg attached to a vertex $i$ such that $d_i=n$ to the symmetrised Euler form is
\[
 2\sum_{r=1}^{n-1}(n-r)^2-2\sum_{r=1}^{n-1}(n-r+1)(n-r)=-2\sum_{r=1}^{n-1}(n-r)=-2\binom n2.
\]
We obtain the formulas by summing over the set of original vertices $Q_0$.
\end{proof}

\subsection{Injectivity of legs}
\label{subsection:injectivity-legs}
Fix $i\in Q_0$. We concentrate on the leg at the $i$th vertex and we drop the label $i$ to simplify the notation, so for example we denote $V_r=V_{i,r}$. Write
\[
 a_r\colon V_r\longrightarrow V_{r-1},\qquad
 b_r\colon V_{r-1}\longrightarrow V_r, \qquad 1\leq r\leq d_i-1.
\]
Then, the $a_r$ are the maps pointing towards the initial vertex $i=(i,0)$ while the $b_r$ point in the opposite direction. With the moment-map convention of \Cref{sec:notation}, the relations are
\begin{equation}\label{eq:legrelations}
 b_ra_r=a_{r+1}b_{r+1}\quad(1\leq r<n-1),\qquad b_{n-1}a_{n-1}=0,\qquad \mu_i(x)+a_1b_1=0.
\end{equation}
We define
\begin{equation}
\label{AFlag_def}
A_r=a_1\cdots a_r,\qquad A_0=\id,\qquad \eta=a_1b_1,\qquad F_r=\mathrm{im}A_r,
\end{equation}
so $A_r$ is the composition of the maps corresponding to the path $(i,r)\rightarrow (i,r-1)\rightarrow\hdots\rightarrow(i,0)=i$.

\begin{lemma}\label{lem:injectivity}
Suppose $\Theta_{i,r}>0$ at all added vertices $(i,r)$, $i\in Q_0$, $1\leq r\leq d_i-1$, and $\Theta\cdot\dd^+=0$. Every $\dd^+$-dimensional $\Theta$-semistable representation of $\Pi_{Q^+}$ has all maps $a_r$, $1\leq r\leq d_i-1$ injective.
\end{lemma}
\begin{proof}
Define $K_r=\ker A_r$.  Note that $a_rK_r\subset K_{r-1}$.  The relations give
\[
 A_rb_r=(a_1b_1)A_{r-1}.
\]
It follows that $b_rK_{r-1}\subset K_r$. Thus the collection of vector space $K_r$, extended by zero at the original vertices and other legs, form a subrepresentation $\rho$.  Note that if one of the $a_r$ is not injective, then some $K_r\ne0$. Then $\rho$ is nonzero and proper, since $K_0=0$, and its $\Theta$-weight is strictly positive. This contradicts semistability.
\end{proof}

\begin{lemma}\label{lem:chamber}
There is a rational stability parameter $\Theta\in\BoQ^{Q^+_0}$, positive at every added vertex, which is generic for $\dd^+$ and satisfies $\Theta\cdot\dd^+=0$.
\end{lemma}
\begin{proof}
We assign arbitrary fixed positive real numbers at the added vertices. By adjusting the coordinates at the original vertices, we can arrange the equality $\Theta\cdot\dd^+$.  We deduce that the cone of stability parameters positive at the added leg vertices intersects $(\dd^+)^\perp$ nontrivially. For a proper nonzero integral vector $\ee\leq\dd^+$, the equation $\Theta\cdot\ee=0$ cannot hold identically on the hyperplane $(\dd^+)^{\perp}$ of stability conditions: otherwise $\ee$ would be proportional to $\dd^+$, and since one of the coordinates of $\dd^+$ is equal to one this would imply $\ee=\dd^+$. We may therefore choose rational $\Theta$ that avoids the finitely many resulting proper hyperplanes $\ee^{\perp}$ for proper non-zero integral vectors $\ee\leq\dd^+$.
\end{proof}

We fix a stability parameter $\Theta\in\mathbf{R}^{Q_0^+}$ as in \Cref{lem:chamber}. Then
\begin{equation}
\label{equation:smoothstable}
 \Y=\FM_{\Pi_{Q^+},\dd^+}^{\Theta}
\end{equation}
is a smooth stack, since the stability parameter $\Theta$ is generic for $\dd^+$ so that stability equals semistability. It may be empty. It is a scalar gerbe over a smooth variety $Y$, where $Y=\CM_{\Pi_{Q^+},\dd^+}^{\Theta}$ is the good moduli space. The inclusion provided  by \Cref{lem:injectivity} is all that will be used in the sequel: the stack $\Y$ is contained in the open substack of $\FM_{\Pi_{Q^+},\dd^+}$ defined by the condition that the maps $a_r$ are injective. They then specify a flag in the original vertex vector space $V_{i,0}$.

\subsection{The flag model}
\label{subsection:flag-model}
Let
\[
\CB=\prod_i\Flag(V_i)=G_{\dd}/B
\]
be the variety of decreasing full flags of the vector spaces $V_i$, $i\in Q_0$, where $B\subset G_{\dd}$ is the corresponding Borel subgroup. Its dimension is $N$. A cotangent vector at a flag $V_i=F_{i,0}\supset F_{i,1}\supset\cdots\supset F_{i,d_i}=0$ is identified, by the trace pairing, with an endomorphism $\eta$ of $V_i$ satisfying
\begin{equation}
\label{equation:lowering-leg}
 \eta F_{i,r}\subset F_{i,r+1}.
\end{equation}
We denote by $\Zleg$ the open substack of the preprojective stack $\FM_{\Pi_{Q^+},\dd^+}$ of the extended quiver on which the inward arrows $a_r$, $1\leq r\leq n-1$, pointing towards the original vertex $(i,0)$, are injective.

\begin{proposition}\label{prop:flagmodel}
There is an equivalence of stacks
\[
 \Zleg\cong\bigl[\{(x,F,\eta)\in V_{\dd}\times T^*\CB:                         \mu(x)+\eta=0\}/G_{\dd}\bigr].
\]
\end{proposition}
\begin{proof}
For a representation of $\Pi_{Q^+}$ with injective legs, with inwards maps $a_r$ at the vertex $i$, we define $A_r,F_r,\eta$ as in \eqref{AFlag_def}.  Repeated use of \eqref{eq:legrelations} shows $\eta_i A_r=A_{r+1}b_{r+1}$. Therefore $\eta F_r\subset F_{r+1}$. The preprojective relation at the original-vertex is $\mu_i(x)+\eta=0$.

Conversely, from $(F,\eta)$, identify $V_r$ with $F_r$. Take $a_r$ to be the inclusion and
\[
 b_r=\eta|_{F_{r-1}}\colon F_{r-1}\longrightarrow F_r.
\]
The lowering condition (\Cref{equation:lowering-leg}) gives the preprojective relation on every leg, and the original moment map equation for the quiver $Q$ gives the remaining relation. These operations commute with base change. Choices of identifications $V_r\cong F_r$ are exactly the changes of basis at the added vertices. Thus the constructions are inverse on quotient stacks and provide the equivalence of stacks of the proposition.
\end{proof}

\subsection{The derived flag model}\label{subsec:derivedflag}
We will need the compatibility of \Cref{prop:flagmodel} with reduced derived structures, which does not follow just by identifying classical points as we did in \Cref{subsection:flag-model}. We show that the moment map equations on the legs of the quiver are regular, which shows that the leg equations contribute no derived structure, so
that the reduced derived structure on $\mathcal{Z}$ is the derived zero locus of $\mu+\eta$ on $[V_{\dd} \times T^*\CB/G_{\dd}]$.

On a fixed flag $F_0\supset F_1\supset\hdots\supset F_{n-1}\supset F_n=0$, let
\[
 \mathcal U=\bigoplus_{r=1}^{n-1}\Hom(F_{r-1},F_r), \qquad \mathcal T=\bigoplus_{r=1}^{n-1}\End(F_r).
\]
Define a linear map $D:\mathcal U\to\mathcal T$ by
\[
 D(b)_r=b_r|_{F_r}-a_{r+1}b_{r+1},
\]
where the second term is zero at the last vertex and $a_r\colon F_r\rightarrow F_{r-1}$ is the inclusion.

\begin{lemma}\label{lem:linear-elimination}
The map $D$ is surjective, and its kernel is canonically the cotangent space at the flag. In families this is an exact sequence of vector bundles
\[
 0\longrightarrow T^*\CB \longrightarrow\mathcal U \xrightarrow{D}\mathcal T \longrightarrow0
\]
over the product of flag varieties $\CB$.
\end{lemma}
\begin{proof}
Given $(c_r)_{1\leq r\leq n-1}\in\mathcal{T}=\bigoplus_{r=1}^{n-1}\End(F_r)$, choose $b_{n-1}$ extending $c_{n-1}$ from $F_{n-1}$ to $F_{n-2}$. Working successively from the end vertex of the $i$th leg towards the initial vertex $i=(i,0)$, choose $b_r$ whose restriction to $F_r$ is $c_r+a_{r+1}b_{r+1}$. Such extensions always exist. This proves fibrewise surjectivity of $D$ in the case in which $Q$ has one vertex, and since the ranks are constant, it proves the surjectivity of $D$ as a morphism of vector bundles in this case too.  Iterating the argument for quivers with multiple vertices proves the general case.

If $D(b)=0$, set $\eta=a_1b_1$. Then, we have $b_1|_{F_1}=a_2b_2,\hdots,b_{n-2}|_{F_{n-2}}=a_{n-1}b_{n-1}, b_{n-1}|_{F_{n-1}}=0$. Since all $a_r$ are injective, $a_1b_1$ determines $b$ and the equation $D(b)=0$ gives $a_1b_1|_{F_{r}}=a_1\hdots a_{r+1}b_{r+1}$. Since $b_{r+1}$ sends $F_r$ to $F_{r+1}$, it follows that $a_1b_1$ is a lowering endomorphism of $F_0$ with respect to the flag $F_{\bullet}$. Therefore, the kernel of $D$ can be identified with the lowering endomorphisms of $F_{\bullet}$. The rank of $\ker(D)$ is the difference
\[
 \sum_{r=1}^{n-1}(n-r)(n-r+1)-\sum_{r=1}^{n-1}(n-r)^2=\binom n2=N,
\]
as required.
\end{proof}

Consequently the derived vanishing of the leg moment equations does not produce any non-trivial derived structure: these equations eliminate independent linear directions. After this elimination, the remaining reduced equation is $\mu+\eta=0$ with values in $\mathfrak{gl}_{\dd}^{\,0}$. Indeed the sum of the traces of commutators is zero, and $\eta$ is traceless since it is a lowering endomorphism. We obtain an identification of the reduced derived injective open $\Zleg\subset\FM_{\Pi_{Q^+},\dd^+}$ with the derived zero locus of $\mu+\eta$ in the flag model. This proves the derived compatibility needed below since the equivalence of \Cref{prop:flagmodel} naturally upgrades to an equivalence of derived stacks.

\section{The flag correspondence and the Vandermonde determinant}\label{sec:correspondence}
\subsection{The correspondence}
\label{subsection:the-correspondence}
Let $\CE$ be the rank-$N$ vector bundle on $\Zleg$ (the substack of $\FM_{\Pi_{Q^+},\dd^+}$ of modules for which the inward arrows $a_{i,r}$ are injective) pulled back
from $T^*\CB$ via the second projection $\Zleg\rightarrow [\CB/G_{\dd}]$: it fits into the Cartesian square
\[\begin{tikzcd}
	{\mathcal{E}} & {[(T^*\mathcal{B})/G_{\dd}]} \\
	{\mathfrak{Z}} & {[\mathcal{B}/G_{\dd}]}
	\arrow[from=1-1, to=1-2]
	\arrow[from=1-1, to=2-1]
	\arrow["\lrcorner"{anchor=center, pos=0.125}, draw=none, from=1-1, to=2-2]
	\arrow[from=1-2, to=2-2]
	\arrow[from=2-1, to=2-2]
\end{tikzcd}\]
The cotangent coordinate $\eta$ in $\Zleg$ (see \Cref{prop:flagmodel}) is a section $s$ of $\CE$. Its classical zero locus is the quotient stack
\[
 \Fr=[X_{\dd}\times\CB/G_{\dd}].
\]
There is a correspondence
\begin{equation}\label{eq:flagcorr}
 \FM_{\Pi_Q,\dd}\xleftarrow{\ p\ }\Fr\xrightarrow{\ j\ }\Zleg.
\end{equation}
The map $p$ is induced by the $G_{\dd}$-equivariant projection $X_{\dd}\times\CB\rightarrow X_{\dd}$ and so is a smooth proper flag bundle, with fiber $\CB$, and is of relative dimension $N$, and $j$ is the closed immersion of the zero-locus of the section $\eta$ of the vector bundle $\CE$ over $\Zleg$. If $\iota:\Y\hookrightarrow\Zleg$ is the stable open subset for the stability parameter $\Theta$ chosen by \Cref{lem:chamber} (see \Cref{lem:injectivity} for the inclusion), we define
\begin{equation}\label{equation:flag-map}
 \Phi=\iota^*j_*p^!: \HO_m^{\BM}(\FM_{\Pi_Q,\dd})\longrightarrow \HO_{m+2N}^{\BM}(\Y).
\end{equation}
We may refer to $\Phi$ as the \emph{flag map}.
\begin{lemma}\label{lem:Phi-linear}
The map $\Phi$ \eqref{equation:flag-map} is $R_{\dd}$-linear and preserves the normalised CoHA grading, i.e. it induces a cohomologically graded map $\HO^{\mathrm{BM}}_*(\FM_{\Pi_Q,\dd},\mathbb{D}\BoQ^{\mathrm{vir}})\rightarrow\HO^{\mathrm{BM}}_*(\Y,\mathbb{D}\BoQ^{\vir})$.
\end{lemma}
\begin{proof}
We first prove $R_{\dd}$-linearity. Recall that $R_{\dd}=\HO^*(\BG G_{\dd},\BoQ)$ is generated by the Chern classes of the universal bundles at the original vertices $i\in Q_0$, without involving the extra vertices at legs of $Q^+$.

For every $i\in Q_0$, let $\CV_i$ denote the corresponding universal bundle, with a superscript indicating the stack over which it is considered when necessary. By construction of the flag correspondence \eqref{eq:flagcorr}, there are canonical isomorphisms
\[
p^*\CV_i^{\FM_{\Pi_Q,\dd}}\simeq \CV_i^{\Fr}\simeq j^*\CV_i^{\Zleg},\qquad\iota^*\CV_i^{\Zleg}\simeq \CV_i^{\Y}.
\]
Indeed, the morphism $p$ forgets the flags, while $j$ identifies the correspondence with the zero locus of the cotangent coordinate. Neither operation changes the universal bundles at the original vertices.

Consequently, for every tautological class $f\in R_{\dd}$, its pullbacks to the correspondence satisfy $p^*f=j^*f$.  Let
\[
\alpha\in\HO_m^{\BM}(\FM_{\Pi_Q,\dd},\BoQ).
\]
Since $p$ is smooth of relative dimension $N$, smooth Gysin pullback commutes with cap products:
\[
p^!(f\cap\alpha)=p^*f\cap p^!\alpha.
\]
Moreover, since $j$ is a closed immersion, the projection formula for proper pushforward gives
\[
j_*(j^*f\cap\beta)=f\cap j_*\beta
\]
for every $\beta\in\HO_*^{\BM}(\Fr,\BoQ)$. Finally, restriction to an open subset commutes with cap products.

Applying these identities successively, we obtain
\[
\begin{aligned}
\Phi(f\cap\alpha)
&=\iota^*j_*p^!(f\cap\alpha)\\
&=\iota^*j_*(p^*f\cap p^!\alpha)\\
&=\iota^*j_*(j^*f\cap p^!\alpha)\\
&=\iota^*(f\cap j_*p^!\alpha)\\
&=\iota^*f\cap\iota^*j_*p^!\alpha\\
&=f\cap\Phi(\alpha).
\end{aligned}
\]
This proves that $\Phi$ is $R_{\dd}$-linear.

We now check the degree statement. If $m=-k-(\dd,\dd)_Q$, then \Cref{lem:legdimension} gives
\[
 m+2N=-k-(\dd^+,\dd^+)_{Q^+}.
\]
Since $N$ is the relative dimension of $p$, the target has the same normalised virtual degree $k$ \eqref{eq:cohgrading} as the source.
\end{proof}
\subsection{Flag Chern roots}
On $[\CB/G_{\dd}]$, and hence on $\Zleg$ by pulling back along the second projection $\Zleg\rightarrow[\CB/G_{\dd}]$, put
\[
 \CL_{i,r}=F_{i,r-1}/F_{i,r},\qquad x_{i,r}=c_1(\CL_{i,r}),\qquad 1\leq r\leq d_i,
\]
where $(F_{i,r})_{i\in Q_0,1\leq r\leq d_i-1}$ are the tautological flags over $\CB$. The total Chern classes are
\[
 c(\CV_{i,r})=\prod_{s=r+1}^{d_i}(1+x_{i,s}).
\]
In particular the tautological classes of the universal bundles $\CV_i$ at the original vertices $i\in Q_0$ are the elementary symmetric polynomials in the $x_{i,r}$.

The associated graded of the cotangent bundle over $\CB$ is
\[
 \operatorname{gr}T^*\CB=\bigoplus_i\bigoplus_{r<s}           \Hom(\CL_{i,r},\CL_{i,s}).
\]
One may check this by writing a strictly flag-lowering matrix in a basis
adapted to the flag. Hence, since taking associated graded does not change Chern classes,
\begin{equation}\label{eq:euler}
 c_N(\CE)=\prod_i\prod_{r<s}(x_{i,s}-x_{i,r})           =(-1)^N\Delta_{\dd}\in\HO^{2N}(\Zleg).
\end{equation}
Here we use the Vandermonde determinant $\Delta_{\dd}$ \eqref{equation:vandermonde}.
\subsection{Image of the fundamental class under the flag correspondence}
\begin{proposition}\label{prop:vfc}
The flag correspondence $\Phi$ \eqref{equation:flag-map} satisfies
\[
 \Phi(v_{\dd})=(-1)^N\Delta_{\dd}\cap[\Y].
\]
More generally,
\[
 \Phi(f\cap v_{\dd})=(-1)^N\kappa(f\Delta_{\dd})\cap[\Y],\qquad f\in R_{\dd},
\]
where $\kappa\colon R_d\rightarrow\HO^*(\Y)$ is the universal map obtained by pullback in cohomology from $\Y\rightarrow \BG G_{\dd}$.
\end{proposition}
\begin{proof}
By \Cref{subsec:derivedflag} and \Cref{prop:flagmodel}, the reduced derived structure on $\Zleg$ is the zero locus of $\mu+\eta\colon [V_{\dd}\times T^*\CB]\rightarrow [\mathfrak{gl}_{\dd}^0/G_{\dd}]$. The derived zero locus of its induced section $s=\eta\colon \Zleg\rightarrow\CE$ is
\begin{equation}\label{eq:derivedF}
 \Fr^{\rd}=\FM_{\Pi_Q,\dd}^{\rd}\times_{\BG G_{\dd}}[\CB/G_{\dd}].
\end{equation}
Indeed, the equations $\mu+\eta=0$ and $\eta=0$ together are the equation $\mu=0$; the cotangent coordinate and its equation eliminate one another. This is therefore an identification of the derived fibre products.

Smooth pullback in \eqref{eq:derivedF} gives
\[
 p^!v_{\dd}=[\Fr^{\rd}]^{\vir}.
\]

We now apply the virtual zero-section formula. Recall that $\CE$ is a vector bundle of rank $N$ on $\Zleg^{\rd}$, and that the cotangent coordinate $\eta$ defines a section $s$ of $\CE$. By \eqref{eq:derivedF}, its derived zero locus is $\Fr^{\rd}$.

We recall the general formula for the pushforward. Let $\mathcal X$ be a quasi-smooth derived stack, let $E$ be a vector bundle of rank $r$ on $\mathcal X$, and let $s$ be a section with derived zero locus
\[
j:Z^{\rd}(s)\longrightarrow\mathcal X.
\]
Then
\begin{equation}\label{eq:virtualzero} j_*[Z^{\rd}(s)]^{\vir}=c_r(E)\cap[\mathcal X]^{\vir}.
\end{equation}

Indeed, let
\[
0,\Gamma_s:\mathcal X\longrightarrow\mathrm{Tot}(E)
\]
be the zero section and the section $s$. The derived zero locus is their derived intersection. We have
\[
[Z^{\rd}(s)]^{\vir}=0_{\Gamma_s}^![\mathcal X]^{\vir}.
\]
Here $0_{\Gamma_s}^!$ denotes the refined Gysin map associated with the zero section and the morphism $\Gamma_s$.

By compatibility of refined Gysin pullback with proper pushforward,
\[
j_*0_{\Gamma_s}^![\mathcal X]^{\vir}=0^!\Gamma_{s,*}[\mathcal X]^{\vir}.
\]
Homotopy invariance for the vector bundle $E$ gives
\[
\Gamma_{s,*}[\mathcal X]^{\vir}=0_*[\mathcal X]^{\vir}.
\]
The self-intersection formula for the zero section therefore implies
\[
\begin{aligned}
j_*[Z^{\rd}(s)]^{\vir}
&=0^!0_*[\mathcal X]^{\vir}\\
&=c_r(E)\cap[\mathcal X]^{\vir}.
\end{aligned}
\]
This proves \eqref{eq:virtualzero}.  Applying \eqref{eq:virtualzero} to the section $s=\eta$ of $\CE$ on $\Zleg^{\rd}$ gives
\[
j_*[\Fr^{\rd}]^{\vir}=c_N(\CE)\cap[\Zleg^{\rd}]^{\vir}.
\]
Consequently,
\[
j_*p^!v_{\dd}=c_N(\CE)\cap[\Zleg^{\rd}]^{\vir}.
\]

Virtual fundamental classes commute with restriction to open substacks. Moreover, by \Cref{lem:stablefund}, the reduced virtual fundamental class on $\Y$ is its ordinary fundamental class. Thus
\[
\begin{aligned}
\Phi(v_{\dd})
&=\iota^*j_*p^!v_{\dd}\\
&=\iota^*\bigl(
c_N(\CE)\cap[\Zleg^{\rd}]^{\vir}
\bigr)\\
&=c_N(\CE|_{\Y})\cap[\Y]\\
&=(-1)^N\Delta_{\dd}\cap[\Y],
\end{aligned}
\]
where the last equality follows from \eqref{eq:euler}. The general assertion follows from \Cref{lem:Phi-linear}.
\end{proof}

\section{Families of quiver varieties}\label{sec:family}
We now keep the stability parameter $\Theta\in\mathbf{Q}^{Q^+_0}$ given by \Cref{lem:chamber} fixed and vary the complex moment parameter for the extended quiver. This gives the comparison between the target $\Y$ of the flag correspondence and the generic spectral varieties $Y_{\sigma}=f^{-1}(\sigma)$, see \eqref{eq:diagram} used in the sign character formula in \cite{HLRV}. The construction and its geometric properties are due to Dobrovolska--Ginzburg--Travkin \cite[\S9.2]{DGT}. We give the details needed for the tautological classes and their specialization and to identify how they transform under the Weyl group action.

\subsection{Ordered spectra and the doubled-quiver quotient}
\label{family_setup_sec}
Let
\[
 G=PG_{\dd}=G_{\dd}/\Gm,\qquad \Fg=\mathrm{Lie}(G),\qquad\mathfrak h=\left\{(\sigma_{i,r})_{i\in Q_0, 1\leq r\leq d_i}:\sigma_{i,r}\in\BoC, \quad\sum_{i,r}\sigma_{i,r}=0\right\}.
\]
The trace pairing identifies the dual $\mathfrak g^*$ with the matrices in $\mathfrak{gl}_{\dd}$ of total trace zero. Then $\mathfrak h$ is a Cartan subspace and $\mathfrak g^*/\!\!/G=\mathfrak h/\!\!/W$, where $W=\prod_{i\in Q_0}\mathfrak{S}_{d_i}$, by the Chevalley isomorphism. Let
\[
 V=V_{\dd}=\Rep_{\dd}(\overline Q),\qquad M=V/\!\!/G,\qquad \chi:M\longrightarrow\mathfrak h/\!\!/W
\]
where $\chi$ is the map induced by the moment map and the adjoint quotient. The affine variety $M$ parametrizes semisimple representations of the double quiver $\overline{Q}$ and $\chi$ remembers the unordered set of eigenvalues of the moment map $\mu(x)$ for $x\in V$.

Denote by $\mathfrak h^{\rm gen}\subset\mathfrak{h}$ the open subset of spectra with distinct eigenvalues at each vertex and no nonempty proper selection of eigenvalues of total sum zero. Let
\begin{equation}
\label{equation:definitionZ}
 Z=(M\times_{\mathfrak h/\!\!/W}\mathfrak h)_{\rm red}
\end{equation}
be the reduced subscheme of the fiber product $M\times_{\mathfrak{h}/\!\!/W}\mathfrak{h}$. There is a finite map $a:Z\to M$ (the pullback of $\mathfrak{h}\rightarrow\mathfrak{h}/\!\!/W$), and $W$ acts on $Z$ by permuting the ordered eigenvalues (action on the $\mathfrak{h}$-factor). There is a morphism $p\colon Z\rightarrow\mathfrak{h}$ given by projection on the second factor. We write $Z^{\rm gen}\coloneqq p^{-1}(\mathfrak{h}^{\rm{gen}})$ for its restriction to $\mathfrak h^{\rm gen}$.

\subsection{Central parameters and eigenflags}
We recall the classical connection between deformations of moment maps of quivers with legs and Hamiltonian reductions at regular semisimple orbits, (following \cite[Lemma~9.1]{CB03Normality}, \cite[\S~1 and \S~3]{CB03Matrices}, \cite[Proposition~5.2.2]{Letellier}). Let
$\xi=(\xi_{i,r})_{i\in Q_0, 0\leq r\leq d_i-1}\in\BoC^{Q_0^+}$ be a complex moment parameter for $Q^+$, satisfying $\xi\cdot\dd^+=0$. The deformed preprojective equations on a given leg are, using the notation dropping the label of the leg introduced in \Cref{subsection:injectivity-legs},
\begin{equation}\label{eq:complexleg}
 b_ra_r-a_{r+1}b_{r+1}=-\xi_{i,r}\id_{V_r}, \qquad\mu_i(x)+a_1b_1=\xi_{i,0}\id_{V_i}.
\end{equation}
The minus sign before $\xi_{i,r}$ comes from our convention regarding the preprojective relation and the definition of the quiver $Q^+$, \Cref{subsection:quivers-moment-maps} and \Cref{sec:legs}. At the last vertex the term $a_{r+1}b_{r+1}$ is set to be zero.
Define
\begin{equation}\label{eq:eigenparameters}
 \sigma_{i,r}=\sum_{k=0}^r\xi_{i,k}, \qquad 0\leq r\leq d_i-1.
\end{equation}
Then
\[
 \sum_{i,r}\sigma_{i,r}    =\sum_i\sum_{k=0}^{d_i-1}(d_i-k)\xi_{i,k}=\xi\cdot\dd^+=0.
\]
The transformation is invertible:  $\xi_{i,0}=\sigma_{i,0}$ and $\xi_{i,r}=\sigma_{i,r}-\sigma_{i,r-1}$ for $r\geq 1$.  There is an identification between the space of central deformation parameters $\xi$ for $Q^+$ and the space of ordered eigenvalues in $\mathfrak{h}$.
\begin{lemma}\label{lem:eigenflag}
Suppose that $\sigma_{i,0},\ldots,\sigma_{i,d_i-1}$ are pairwise distinct for each $i$, and that
\[
 \sum_i\sum_{r=0}^{d_i-1}\sigma_{i,r}=0.
\]
Then every representation satisfying the $\xi$-deformed preprojective relations \eqref{eq:complexleg} has injective inward leg maps $a_r\colon V_r\rightarrow V_{r-1}$. Forgetting the legs induces an equivalence of stacks
\[
 [\mu_{\dd^+}^{-1}(\xi)/G_{\dd^+}]\xrightarrow{\sim} [\mu_{\dd}^{-1}(\mathcal O_\sigma)/G_{\dd}],
\]
where $\mathcal O_\sigma\subset\mathfrak{gl}_{\dd}^0$ is the $G_{\dd}$-orbit of the tuple
\[
 \bigl(\operatorname{diag}(\sigma_{i,0},\ldots,\sigma_{i,d_i-1})\bigr)_{i\in Q_0}.
\]
Injectivity induces a full flag $F_{\bullet}$ of $V_i$, such that for each $i$, the flag quotient $F_{i,r-1}/F_{i,r}$ is canonically identified with the $\sigma_{i,r-1}$-eigenline bundle, for $1\leq r\leq d_i$.
\end{lemma}

\begin{proof}
We fix $i\in Q_0$, put $n=d_i$, drop the index $i$ and use the notations of \Cref{subsection:injectivity-legs}. The case $n=1$ is immediate as there are no legs. For $n\geq2$, write
\[
 a_r:V_r\longrightarrow V_{r-1},\qquad b_r:V_{r-1}\longrightarrow V_r, \qquad \dim V_r=n-r,\qquad 0\leq r\leq n-1.
\]
Set $V_n=0$ and $a_n=b_n=0$. The relations \eqref{eq:complexleg} read
\[
 \mu_i+a_1b_1=\xi_0\id, \qquad b_ra_r=a_{r+1}b_{r+1}-\xi_r\id \quad (1\leq r\leq n-1).
\]
Here, by \eqref{eq:eigenparameters},
\[
 \xi_0=\sigma_0, \qquad \xi_r=\sigma_r-\sigma_{r-1}.
\]

We prove by downward induction on $1\leq r\leq n-1$ that $b_ra_r$ is invertible, with eigenvalues
\[
 -\xi_r,\;-(\xi_r+\xi_{r+1}),\;\ldots,\; -(\xi_r+\cdots+\xi_{n-1}),
\]
or equivalently
\[
 \sigma_{r-1}-\sigma_s, \qquad r\leq s\leq n-1.
\]
For $r=n-1$, this follows from $b_{n-1}a_{n-1}=-\xi_{n-1}\id$. For the induction step, invertibility of $b_{r+1}a_{r+1}$ gives
\[
 V_r=\im(a_{r+1})\oplus\ker(b_{r+1}).
\]
On the first summand, $a_{r+1}b_{r+1}$ is conjugate to $b_{r+1}a_{r+1}$; on the second, which has dimension one, it vanishes. Subtracting $\xi_r\id$ gives the required eigenvalues. They are all nonzero by distinctness of the $\sigma_s$. Thus every $a_r$ is a split injection.

Put
\[
 F_0=V_0,\qquad F_r=\im(a_1\cdots a_r)\quad(1\leq r\leq n-1), \qquad F_n=0,
\]
and let $\eta=a_1b_1$. Then, $F_r$ defines a full flag in $V_i$. Iterating the preprojective relations on the leg \eqref{eq:complexleg} shows that $\eta$ preserves this flag and acts on $F_{r-1}/F_r$ by $-\sum_{k=1}^{r-1}\xi_k$. Consequently, $\mu_i=\xi_0\id-\eta$ acts on this quotient by
\[
 \xi_0+\sum_{k=1}^{r-1}\xi_k=\sigma_{r-1}.
\]
The $\sigma_r$ are therefore all eigenvalues of $\mu_i$, which implies the identity
\[
 \prod_{s=0}^{n-1}(\mu_i-\sigma_s\id)=0.
\]
Since the $\sigma_s$ are distinct, the spectral projectors give a decomposition into eigenlines, and
\begin{equation}\label{equation:direct-sum-decomposition}
 F_r=\bigoplus_{s=r}^{n-1}\ker(\mu_i-\sigma_s\id).
\end{equation}

Conversely, suppose that $\mu_i$ has the prescribed spectrum given by $\sigma_{i,r}$, $i\in Q_0$, $0\leq r\leq d_{i-1}$. Define $F_r$ by the preceding formula \eqref{equation:direct-sum-decomposition}, take $V_r=F_r$, and let $a_r\colon V_r\rightarrow V_{r-1}$ be the inclusion. Define
\[
 b_r= \left.(\sigma_{r-1}\id-\mu_i)\right|_{F_{r-1}} : F_{r-1}\longrightarrow F_r.
\]
This has image in $F_r$ because it vanishes on the $\sigma_{r-1}$-eigenline and preserves all the other eigenlines. We have
\[
 \mu_i+a_1b_1=\sigma_0\id=\xi_0\id
\]
and
\[
 a_{r+1}b_{r+1}-b_ra_r=(\sigma_r-\sigma_{r-1})\id_{F_r} =\xi_r\id_{F_r},
\]
including the terminal relation with $a_n=b_n=0$. Thus \eqref{eq:complexleg} holds.

All constructions commute with base change: the spectral projectors are polynomials in $\mu_i$ whose denominators are nonzero constants. Moreover, the injections $a_1\cdots a_r$ identify every leg vector bundle with $F_r$, and the deformed preprojective relations uniquely determine the reverse maps under these identifications. Hence forgetting and reconstructing the legs are mutually inverse on objects and morphisms in families. This proves the claimed equivalence of stacks.
\end{proof}

We say that $\sigma$ is \emph{generic} if it is regular semisimple, has total trace zero, and no nonempty proper selection of its eigenvalues has sum zero.

\begin{lemma}\label{lem:generic-eigen}
For the parameter $\xi=\Theta$ chosen in \Cref{lem:chamber}, the ordered spectrum $\sigma$ is generic.
\end{lemma}
\begin{proof}
On each leg the differences $\xi_{i,r}=\sigma_{i,r+1}-\sigma_{i,r}$ in \eqref{eq:eigenparameters} are positive. Hence the eigenvalues $\sigma_{i,r}$ are distinct at each vertex.
For subsets $S_i\subset\{1,\ldots,d_i\}$, $i\in Q_0$, set
\[
 e_{i,0}=|S_i|,\qquad e_{i,r}=|\{s\in S_i:s>r\}|.
\]
Then $0\leq\ee\leq\dd^+$ and
\[
 \sum_i\sum_{s\in S_i}\sigma_{i,s}=\xi\cdot\ee.
\]
A nonempty proper selection gives $0<\ee<\dd^+$. The genericity of $\Theta$ excludes a zero sum.
\end{proof}

Therefore, generic spectra $\sigma$ exist.
\subsection{The smooth family of quiver varieties}
\label{subsection:smooth-GIT-family}
We recall the construction of the smooth GIT family of deformed quiver varieties following \cite{DGT}. Put $\alpha=\dd^+$. The parameter transformation
\eqref{eq:eigenparameters} identifies
$\{\xi\in \BoC^{Q_0^+}:\xi\cdot\alpha=0\}$ with $\mathfrak h$.
A smooth quasi-projective variety $X$ is \emph{semiprojective} for a
$\Gm$-action if its fixed locus is proper and every point $x\in X$ has a
limit $\lim_{s\mapsto 0}s\cdot x$. The \emph{core} of $X$ is the proper union of the
repelling strata, equivalently the locus of points whose orbit also
has a limit as $s\to\infty$. We refer to \cite{HRV} for the study of this class of varieties.

Consider the family
 \begin{equation}
 \label{Xcal_def}
 \mathcal X=
 \{(x^+,\xi)\in \mathrm{Rep}_{\dd^+}(\overline{Q^+})\times\BoC^{Q_0^+}_0:\mu_{Q^+}(x^+)=\xi\}^{\Theta\text{-st}}
       /\bigl(\prod_v\GL_{\alpha_v}/\Gm\bigr)
 \xrightarrow{\ f\ }\mathfrak h.
\end{equation}
Here, we use the stability parameter $\Theta$ given by \Cref{lem:chamber} and define $\BoC^{Q_0^+}_0\coloneqq\{\xi\in \BoC^{Q_0^+}\mid \xi\cdot\dd^+=0\}$. The map $f$ only remembers the parameter $\xi$, seen as an element of $\mathfrak{h}$ via the parameter transformation \eqref{eq:eigenparameters}.
Forgetting the leg maps and remembering the spectrum gives morphisms $\mathcal{X}\rightarrow M$ and $\mathcal{X}\rightarrow\mathfrak{h}$ and so, by smoothness of $\mathcal{X}$, a map $\pi$ to the reduced fiber product $Z$:
\begin{equation}\label{eq:diagram}
 \begin{tikzcd}
 \mathcal X \arrow[r,"\pi"] \arrow[dr,"f"'] & Z \arrow[d,"p"] \arrow[r,"a"] & M\arrow[d,"\chi"]\\
 &\mathfrak h\arrow[r]&\mathfrak h/\!\!/W.
 \end{tikzcd}
\end{equation}
The reduced zero fibre $Y=f^{-1}(0)_{\mathrm{red}}$ is the good moduli space of the stack $\Y$ introduced in \eqref{equation:smoothstable}.

\begin{proposition}\label{prop:geometry}
The morphism $\pi$ is projective, $f$ is smooth, and $\pi$
is an isomorphism over $Z^{\rm gen}$. The total space $\mathcal X$
is semiprojective for the action scaling the doubled arrows $a,a^*\in \overline{Q_1^+}$.
\end{proposition}

\begin{proof}
The stability parameter $\Theta$ is generic for the dimension vector $\dd^+$. Consequently stable and semistable points coincide throughout the family. Therefore, $\mathcal{X}$ is a smooth algebraic variety. The stabilizer of a stable representation is the scalar subgroup $\mathbf{G}_{\mathrm{m}}$, and the differential $d\mu_{Q^+}$ of the moment map $\mu_{Q^+}\colon \mathrm{Rep}_{\dd^+}(\overline{Q^+})\rightarrow\mathrm{Lie}\left(\prod_{v}\GL_{\alpha_v}/\mathbf{G}_{\mathrm{m}}\right)^*$, which is dual to the infinitesimal action of $\mathrm{Lie}\left(\prod_{v}\GL_{\alpha_v}/\mathbf{G}_{\mathrm{m}}\right)$, is therefore surjective. Therefore, the differential of the function $\mu_{Q^+}-\xi\colon (\mathrm{Rep}_{\dd^+}(Q^+))^{\Theta\text{-st}}\times\BoC^{Q_0^+}_0\rightarrow \mathrm{Lie}\left(\prod_{v}\GL_{\alpha_v}/\mathbf{G}_{\mathrm{m}}\right)^*$ is also surjective. This proves the smoothness of the map
\[
 \tilde{\mathcal{X}}\coloneqq\{(x^+,\xi)\in \mathrm{Rep}_{\dd^+}(\overline{Q^+})\times\BoC^{Q_0^+}_0:\mu_{Q^+}(x^+)=\xi\}^{\Theta\text{-st}}\xrightarrow{\ \tilde{f}\ }\mathfrak h.
\]
Now, there is a commutative diagram where the vertical map is a principal $\prod_v \GL_{\alpha_v}/\mathbf{G}_{\mathrm{m}}$-bundle:
\[\begin{tikzcd}
	{\tilde{\mathcal{X}}} & {\mathfrak{h}} \\
	{\mathcal{X}}
	\arrow["{\tilde{f}}", from=1-1, to=1-2]
	\arrow[from=1-1, to=2-1]
	\arrow["f"', from=2-1, to=1-2]
\end{tikzcd}\]
This proves smoothness of $f$. Note that the genericity of $\Theta$ for the dimension vector $\dd^+$ is the only property of $\Theta$ that we used to prove the smoothness of $f$.

The positive weights of the stability parameter $\Theta$ force the inward maps $a_r$ to be injective for all central parameters $\xi$, by an argument generalising \Cref{lem:injectivity}. Indeed, on one leg at the $i$th vertex, dropping the label $i$ and using notations as in \Cref{subsection:injectivity-legs}, put $A_r=a_1\cdots a_r$ and $\eta=a_1b_1$. The deformed equations \eqref{eq:complexleg} give
\[
 A_rb_r=\left(\eta+\sum_{1\leq k<r}\xi_{i,k}\id\right)A_{r-1}.
\]
Hence the kernels $K_r=\ker A_r$, with zero vector spaces at the original vertices, again form a subrepresentation, exactly as in \Cref{lem:injectivity}. Quotienting by changes of basis on the legs leaves a complete flag preserved by the original moment map, whose successive eigenvalues are \eqref{eq:eigenparameters}.

We now prove projectivity of $\pi\colon\mathcal{X}\rightarrow Z=(M\times_{\mathfrak{h}/\!\!/W}\mathfrak{h})_{\mathrm{red}}$. Define
\[
 \mathcal R=\{(x^+,\xi)\in \Rep_{\dd^+}(\overline{Q^+})\times\mathbb{C}^{Q^+_0}_0 \mid\mu_{Q^+}(x^+)=\xi\},\qquad H=\prod_{v\in Q_0^+}\GL_{\alpha_v}/\Gm,\]
and put $B=\mathbb C[\mathcal R]^H$. After multiplying $\Theta$ by a positive integer so that the result has integer coordinates, let $\chi_\Theta$ be the corresponding character of $H$ giving the $H$-linearization of the trivial line bundle over $\mathcal{R}$. Since stability and semistability coincide, the relative GIT construction gives
\[
 \mathcal X=\operatorname{Proj}_{\operatorname{Spec}B}\left( \bigoplus_{n\geq0}\mathbb C[\mathcal R]^{H,\chi_\Theta^n}\right).
\]
By classical geometric invariant theory \cite{mumford1994geometric}, the graded algebra $\bigoplus_{n\geq0}\mathbb C[\mathcal R]^{H,\chi_\Theta^n}$ is finitely generated over $B$, so the canonical morphism
\[
 q:\mathcal X\longrightarrow\operatorname{Spec}B
\]
is projective.

It remains to relate this affine quotient to $Z$. By the first fundamental theorem for quiver invariants \cite[Theorem~1]{le1990semisimple} (see also \cite{donkin1994polynomial} for its generalisation to any characteristic, which is not needed here), $\mathbb C[\mathcal R]^H$ is generated by the central parameters $\xi$ and traces of oriented cycles in the doubled quiver $\overline{Q^+}$ (note that $H$ acts trivially on central parameters).

The deformed leg relations allow us to eliminate every excursion into a leg. For example, with $\eta=a_1b_1$, the preceding identity gives
\[
 a_1\cdots a_rb_r\cdots b_1=\prod_{k=0}^{r-1} \left(\eta+\sum_{j=1}^{k}\xi_{i,j}\id\right).
\]
Moreover,
\[
 \eta=\xi_{i,0}\id-\mu_i(x).
\]
Using these relations and cyclicity of the trace, every invariant involving leg arrows reduces to a polynomial in the central parameters and traces of cycles in the original doubled quiver. Thus the natural homomorphism
\[
 \mathbb C[M]\otimes\mathbb C[\mathfrak h]\longrightarrow B
\]
is surjective.

We claim that this homomorphism factors through the coordinate ring of $M\times_{\mathfrak h/\!\!/W}\mathfrak h$, i.e. through
\[
 \mathbb C[M]\otimes_{\mathbb C[\mathfrak h]^W}\mathbb C[\mathfrak h].
\]
Indeed, the deformed moment map equations imply
\begin{equation}
\label{equation:determinant-identity}
 \det(t-\mu_i(x))=\prod_{r=1}^{d_i}(t-\sigma_{i,r}).
\end{equation}
This identity holds on the entire affine moment family $\{(x^+,\xi)\in \mathrm{Rep}_{\dd^+}(\overline{Q^+})\times\BoC^{Q_0^+}:\mu_{Q^+}(x^+)=\xi\}$, without any injectivity assumption of the inward maps in the legs. To see this, use $\det(t\id+ab)=t^{\dim V-\dim W}\det(t\id+ba)$ for maps $a:W\to V$ and $b:V\to W$, together with
\[
 b_ra_r=a_{r+1}b_{r+1}-\xi_{i,r}\id.
\]
Since $\dim V_{i,r-1}=\dim V_{i,r}+1$, successive applications and $\xi_{i,r}=\sigma_{i,r+1}-\sigma_{i,r}$ give
\[
 \begin{aligned}
 &\det_{V_{i,r-1}}
 \bigl((t-\sigma_{i,r})\id+a_rb_r\bigr)\\
 &\qquad =
 (t-\sigma_{i,r})
 \det_{V_{i,r}}
 \bigl((t-\sigma_{i,r+1})\id+a_{r+1}b_{r+1}\bigr).
 \end{aligned}
\]
Iterating, with the last leg space $V_{i,d_i}$ equal to zero, we obtain the claimed identity \eqref{equation:determinant-identity}.

The algebra $\mathbb C[\mathfrak h]^W$ is generated by the coefficients of the characteristic polynomials at the original vertices, and the identity \eqref{equation:determinant-identity} shows that their images in $B$ agree under the two natural homomorphisms.

We therefore obtain a surjection
\[
 \mathbb C[M]\otimes_{\mathbb C[\mathfrak h]^W} \mathbb C[\mathfrak h]\twoheadrightarrow B.
\]
It follows that $(\operatorname{Spec}B)_{\rm red}$ is a closed subscheme of $Z$. Since $\mathcal X$ is smooth, hence reduced, the projective morphism $q$ factors through $(\operatorname{Spec}B)_{\rm red}$, and the induced morphism remains projective. Composing with the closed immersion into $Z$ gives $\pi$, which is therefore again projective.

Over $\mathfrak h^{\rm gen}$, the flag given by the injective inward leg maps $a_r$ is the unique ordered eigenflag for the moment map. Moreover, a subrepresentation of the original doubled representation would select eigenvalues of total sum zero, by taking traces of its moment map. Genericity excludes such a subrepresentation: the representation of the original doubled quiver $\overline{Q}$ must be simple. By \Cref{lem:eigenflag}, the deformed leg equations give a unique extension, up to changes of bases on the legs. This extension is $\Theta$-stable: a nonzero proper subrepresentation must have all the original spaces, since the original doubled representation is simple, and must omit a nonzero part of the added spaces. Its weight is therefore strictly negative since the weight of $\dd^+$ is $\Theta\cdot\dd^+=0$. The construction is compatible with families, so $\pi$ is an isomorphism over $Z^{\rm gen}$.

Finally, scaling every doubled arrow by $s$ scales $\xi$, and hence $\sigma$, by $s^2$. The affine quotient $\operatorname{Spec}B$ has a
nonnegative grading, with degree-zero part $\mathbb C$, and therefore contracts to its vertex $\operatorname{Spec}B_0$. Since $q$ is projective, the valuative criterion for properness shows that every $\Gm$-orbit in $\mathcal X$ has a limit as $s\to0$.

Every fixed point of $\mathcal{X}$ maps to the vertex of $\operatorname{Spec}B$. Thus the fixed locus of $\mathcal X$ is a closed subset of the projective fibre of $q$ over this vertex, and is therefore proper. These are the semiprojectivity conditions, and the proposition is proved.
\end{proof}

\subsection{Constancy of the cohomology}
For clarity, we spell out a weight argument for the restriction maps in our situation, to show the constancy of the cohomology of the family $\mathcal{X}\rightarrow\mathfrak{h}$ following \cite{DGT}. We use the standard convention that $\BoQ(-1)$ has weight two, so that a Tate twist by $(-1)$ increases weights by two.

\begin{proposition}\label{prop:constant}
For every $\sigma\in\mathfrak h$, restriction by the inclusion of the fiber $i_{\sigma}\colon \mathcal{X}_{\sigma}\coloneqq f^{-1}(\sigma)\rightarrow\mathcal{X}$ is an isomorphism
\[
 i_\sigma^*: \HO^*(\mathcal X,\BoQ)\xrightarrow{\sim}\HO^*(\mathcal X_\sigma,\BoQ).
\]
These are ring isomorphisms, compatible with mixed Hodge structures, and the action of Chern classes of respective tautological bundles. The sheaves $R^k f_*\BoQ$ and $R^k f_!\BoQ$ are constant, that is the complexes of constructible sheaves $f_*\BoQ_{\mathcal{X}}$ and $f_!\BoQ_{\mathcal{X}}$ are constant over $\mathfrak{h}$.
\end{proposition}

\begin{proof}
We give the argument of \cite[Theorem~10.2.2]{DGT}, spelling out the localization sequence, the weight argument, and the passage to derived direct images. All cohomology groups have rational coefficients.

\smallskip
\noindent\emph{The common core.}
Let $C$ be the proper core of $\mathcal X$, i.e. the repelling locus of the action scaling all arrows of $\overline{Q^+}$ with weight one and the parameter $\xi$ with weight $2$, and put $Y=f^{-1}(0)$. The action of $\mathbf{G}_{\mathrm{m}}$ rescaling all arrows of $\overline{Q^+}$ induces an action on $\mathfrak h$ of weight $2$ so that $f$ is equivariant. Consequently, if a point $x\in\mathcal X$ has a limit under $t\to\infty$, then $f(x)=0$: otherwise the orbit $t^2f(x)$ would have no limit in $\mathfrak h$. It follows that $C\subseteq Y=f^{-1}(0)$. Conversely, the condition defining the core of $Y$ is the same as that defining the core of $\mathcal X$ for points of $Y$, since $Y$ is a closed $\mathbf{G}_{\mathrm{m}}$-invariant subvariety. Thus $Y$ and $\mathcal X$ have the same core $C$.

The variety $Y$ is smooth because $f$ is smooth, and it is semiprojective: its fixed locus is a closed subvariety of the proper fixed locus of $\mathcal X$, and the limit as $t\to0$ of a point of $Y$ remains in $Y$. The cohomology-to-core theorem gives restriction isomorphisms
\begin{equation}\label{equation:common-core-XY}
 \HO^*(\mathcal X,\BoQ)\xrightarrow{\sim}\HO^*(C,\BoQ),\qquad \HO^*(Y,\BoQ)\xrightarrow{\sim}\HO^*(C,\BoQ).
\end{equation}
Their compatibility with the inclusion $Y\hookrightarrow\mathcal X$ shows that restriction from $\mathcal X$ to $Y$ is an isomorphism. Moreover, the cohomology of both smooth semiprojective varieties is pure of weight $n$ in cohomological degree $n$. These are \cite[Theorem~1.3.1 and Corollary~1.3.2]{HRV} in our specific context.

\smallskip
\noindent\emph{Reduction to a one-parameter family.}
Fix $\sigma\ne0$ and form the Cartesian square
\[
\begin{tikzcd}
 X'\arrow[r,"q"]\arrow[d,"g"'] & \mathcal X\arrow[d,"f"]\\
 \mathbf A^1\arrow[r,"z\mapsto z^2\sigma"'] & \mathfrak h.
\end{tikzcd}
\]
Thus
\[
 X'=\{(z,x)\in\mathbf{A}^1\times\mathcal{X}\mid f(x)=z^2\sigma\}, \qquad g(z,x)=z.
\]
Since $g$ is a base change of $f$, it is smooth; in particular, $X'$ is smooth. The action
\[
 t\cdot(z,x)=(tz,t\cdot x)
\]
is well defined because $f(t\cdot x)=t^2f(x)$. Every point $x\in X'$ has a limit $\lim_{t\mapsto 0}t\cdot x$, inherited from the corresponding limit in $\mathcal X$. Its fixed locus is $(X')^{\mathbf{G}_{\mathrm{m}}}=\{0\}\times\mathcal X^{\Gm}$ and is proper. Hence $X'$ is semiprojective again. A point of its core must have $z=0$, and its second coordinate must belong to $C$. Therefore the core of $X'$ is $\{0\}\times C$, identified with $C$.

We identify $g^{-1}(0)$ with $f^{-1}(0)=Y$, and put
\[
 F=\mathcal X_\sigma,\qquad U=X'\setminus Y.
\]
The $\mathbf{G}_{\mathrm{m}}$-action gives an algebraic isomorphism
\begin{equation}\label{equation:product-U}
 \Gm\times F\xrightarrow{\sim}U,\qquad (z,y)\longmapsto(z,z\cdot y),
\end{equation}
whose inverse sends $(z,x)$ to $(z,z^{-1}\cdot x)$. The common-core comparison also gives
\begin{equation}\label{equation:common-core-comparison}
 i^*:\HO^*(X',\BoQ)\xrightarrow{\sim}\HO^*(Y,\BoQ), \qquad i:Y\hookrightarrow X',
\end{equation}
and $\HO^k(X',\BoQ)$ and $\HO^k(Y,\BoQ)$ are pure of weight $k$, again by \cite{HRV}.

\smallskip
\noindent\emph{The localization sequence and its weights.}
Write $j:U\hookrightarrow X'$ for the open immersion. Since $g$ is smooth, $Y=g^{-1}(0)$ is a smooth effective Cartier divisor in $X'$. Purity for this codimension-one immersion gives
\[
 i^!\BoQ_{X'}\cong\BoQ_Y[-2](-1).
\]
Applying derived global sections to the localization triangle
\[
 i_*i^!\BoQ_{X'}\longrightarrow\BoQ_{X'} \longrightarrow Rj_*\BoQ_U\xrightarrow{+1}
\]
gives the long exact sequence of mixed Hodge structures
\begin{equation}\label{eq:gysin}
 \HO^{k-2}(Y)(-1)\xrightarrow{i_*}\HO^k(X')
 \xrightarrow{j^*}\HO^k(U)
 \xrightarrow{\partial}\HO^{k-1}(Y)(-1)
 \xrightarrow{i_*}\HO^{k+1}(X').
\end{equation}

Exactness expresses $\HO^k(U)$ as an extension of the cokernel of the first map in \eqref{eq:gysin} by the kernel of the last map. The cokernel is a quotient of the pure weight-$k$ structure $\HO^k(X')$. The kernel is a subobject of $\HO^{k-1}(Y)(-1)$, which is pure of weight $k+1$. It follows that
\[
 \operatorname{Gr}^W_\ell\HO^k(U)=0\qquad\text{unless }\ell=k\text{ or }\ell=k+1.
\]

On the other hand, the product description of $U$ \eqref{equation:product-U} and the identities $\HO^0(\Gm)=\BoQ$, $\HO^1(\Gm)=\BoQ(-1)$ give
\[
 \HO^k(U)\cong\HO^k(F)\oplus\HO^{k-1}(F)(-1)
\]
as mixed Hodge structures. In degree $k=j+1$, the second summand is $\HO^j(F)(-1)$. The preceding weight bound shows that this summand has weights only $j+1$ and $j+2$, so that $\HO^j(F)$ has weights only $j-1$ and $j$. Since $F$ is smooth, $\HO^j(F)$ has weights at least $j$ by standard Hodge theory \cite{DeligneHodgeII}. Therefore it is pure of weight $j$.

\smallskip
\noindent\emph{Vanishing of the Gysin maps and identification of restriction.}
Every class $a\in\HO^m(Y)$ extends to a class $\widetilde a\in\HO^m(X')$, by the common-core comparison \eqref{equation:common-core-comparison}. The projection formula gives
\[
 i_*a=i_*(i^*\widetilde a\smile1) =\widetilde a\smile i_*1.
\]
Here $i_*1$ is the divisor class of $Y$, namely $c_1(\mathcal O_{X'}(Y))$, with the Tate twists understood. The divisor $Y$ is principal: it is the zero divisor of the global regular function $z=g$ on $X'$. Hence $\mathcal O_{X'}(Y)$ is trivial and $i_*1=0$. All Gysin maps given by $i_*$ in \eqref{eq:gysin} consequently vanish, and that sequence yields
\[
 0\longrightarrow\HO^k(X')\xrightarrow{j^*}\HO^k(U) \xrightarrow{\partial}\HO^{k-1}(Y)(-1)\longrightarrow0.
\]
The left term is pure of weight $k$ and the right term is pure of weight $k+1$. Strictness of morphisms of mixed Hodge structures therefore identifies $\operatorname{im}(j^*)=W_k\HO^k(U)$.  Using the purity of $\HO^*(F)$ and the K\"unneth decomposition, this weight subspace is precisely
\[
 W_k\HO^k(U)=\HO^k(F)\otimes\HO^0(\Gm).
\]
Restriction to $\{1\}\times F\subset U$ is the identity on this summand and vanishes on the summand involving $\HO^1(\Gm)$. Thus the composite
\[
 \HO^k(X')\longrightarrow\HO^k(U)\longrightarrow\HO^k(F)
\]
is an isomorphism. This composite is restriction to $g^{-1}(1)\subset X'$.

The projection $q:X'\to\mathcal X$ is the identity on the identified cores $C$. Its pullback on cohomology is therefore an isomorphism, by compatibility of both cohomology-to-core isomorphisms. The composite
\[
 \HO^*(\mathcal X)\xrightarrow{q^*}\HO^*(X')      \longrightarrow\HO^*(F)
\]
is exactly $i_\sigma^*$. This proves the required restriction isomorphism for $\sigma\ne0$; the case $\sigma=0$ was proved by the common-core comparison \eqref{equation:common-core-XY}. Since these maps are induced by algebraic inclusions, they are ring homomorphisms and morphisms of mixed Hodge structures. Being bijective, they are isomorphisms with both structures.

\smallskip
\noindent\emph{Constancy of the derived direct images.}
Put $B=\mathfrak h$, $r=\dim B$, and $V=R\Gamma(\mathcal X,\BoQ)$. Write $V_B$ for the constant complex on $B$ obtained by pulling back $V$ from a point, that is $V_B=V\otimes \BoQ_{B}$. Let
\[
 p:B\times\mathcal X\longrightarrow B,\qquad \gamma:\mathcal X\longrightarrow B\times\mathcal X, \quad x\longmapsto(f(x),x)
\]
be the projection and graph morphisms, respectively.  The graph $\gamma$ is a closed immersion and $p\circ\gamma=f$. Restriction to the graph gives a morphism of complexes
\[
 \BoQ_{B\times\mathcal X}\longrightarrow \gamma_*\BoQ_{\mathcal X}.
\]
Applying $p_*$, and using the product projection to identify $p_*\BoQ_{B\times\mathcal X}$ with $V_B$, gives a natural map
\[
 \eta:V_B\longrightarrow f_*\BoQ_{\mathcal X}.
\]

We analyse this map on costalks. For $\sigma\in B$, let $u_\sigma:\{\sigma\}\hookrightarrow B$ be the inclusion. Extraordinary base change gives
\[
 u_\sigma^!f_*\BoQ_{\mathcal X}\cong R\Gamma(\mathcal X_\sigma,i_\sigma^!\BoQ_{\mathcal X}).
\]
Since $f$ is smooth, $\mathcal X_\sigma$ is a smooth subvariety of codimension $r$ and purity gives
\[
 i_\sigma^!\BoQ_{\mathcal X} \cong\BoQ_{\mathcal X_\sigma}[-2r](-r).
\]
Likewise $u_\sigma^!V_B\cong V[-2r](-r)$.  Under these identifications, $u_\sigma^!\eta$ is the restriction map
\[
 R\Gamma(\mathcal X,\BoQ)[-2r](-r)\longrightarrow R\Gamma(\mathcal X_\sigma,\BoQ)[-2r](-r).
\]

We have already seen that $\HO^*(\mathcal{X})\rightarrow \HO^*(\mathcal{X}_{\sigma})$ is an isomorphism, and it follows that $u_\sigma^!\eta$ is an isomorphism for every $\sigma$. Costalks detect isomorphisms of constructible complexes, so $\eta$ is an isomorphism.

All the constructions are algebraic and descend from their counterparts in the derived category of mixed Hodge modules. The same map is an isomorphism there as well, since the underlying rational realisation from complexes of mixed Hodge modules to constructible complexes is conservative.

Finally, let $n$ be the relative complex dimension of $f$. Verdier duality and smoothness of $\mathcal X$ and $B$ give
\[
 f_!\BoQ_{\mathcal X}\cong (V^\vee)_B[-2n](-n),
\]
where $V^\vee$ denotes the derived dual of $V$ over $\BoQ$. Indeed, $\dim\mathcal X=n+r$, $\mathbb D_{\mathcal X}\BoQ_{\mathcal X}=\BoQ_{\mathcal X}[2(n+r)](n+r)$, and $\mathbb D_B(V_B)=(V^\vee)_B[2r](r)$. This proves constancy of both derived direct images, and hence of all their cohomology sheaves. Compatibility with tautological Chern classes follows from their functoriality under restriction.
\end{proof}

\begin{remark}\label{rem:emptyfamily}
If $Y$ is empty, then $\mathcal X$ is empty: an orbit limit of any point of $\mathcal X$ would belong to $Y$. The generic spectral variety $f^{-1}(\sigma)$, $\sigma \in \mathfrak{h}^{\mathrm{gen}}$ is then empty as well, by \Cref{prop:geometry}. The statements about small resolutions in \Cref{subsection:smallness} will be made under the assumption $Y\ne\varnothing$.
\end{remark}

\section{Springer theory and tautological classes}\label{sec:Weyl}
\subsection{Smallness of the simultaneous resolution}
\label{subsection:smallness}
Define $Y\coloneqq f^{-1}(0)\subset \mathcal{X}\xrightarrow{f} \mathfrak{h}$ as in the previous section, with $\mathfrak{h}$ as defined in \S \ref{family_setup_sec} and $\mathcal{X}$ as in \eqref{Xcal_def}.  Throughout \S \ref{sign_per_sec}--\ref{sign_per_sec} we assume $Y\ne\varnothing$. Put $r=\dim\mathfrak h=|\dd|-1$ and $D=\dim\mathcal X=2p_{Q^+}(\dd^+)+r$. We normalize intersection complexes by demanding $\IC_T|_{T^{\rm sm}}=\BoQ_{T^{\rm sm}}[\dim T]$, that is we use the perverse normalization.

\begin{proposition}\label{prop:small}
The map $\pi\colon \mathcal{X}\rightarrow Z$ in \eqref{eq:diagram} is a small resolution. In particular,
\[
 \pi_*\BoQ_{\mathcal X}[D]\cong\IC_Z.
\]
\end{proposition}
\begin{proof}
First, $Z^{\rm gen}$ is dense in $Z$. Indeed, smoothness makes the image of $f$ open. This image contains zero, since $Y=f^{-1}(0)\neq\emptyset$ by assumption, so equivariance of $f$ for the positive-weight $\mathbf{G}_{\mathrm{m}}$-action on $\mathfrak{h}$ implies that the image of $f$ is all of $\mathfrak h$. Thus the generic spectral locus in the irreducible variety $M$ is nonempty: it contains $a\pi (x)$ for $x\in f^{-1}(\mathfrak{h}^{\mathrm{gen}})$. The finite flat morphism $\mathfrak h\to\mathfrak h/\!\!/W$ shows that its inverse image by $a$ is dense in $Z$. Indeed, put $q:\mathfrak h\to\mathfrak h/\!\!/W$ and $M^{\rm gen}=\chi^{-1}(q(\mathfrak h^{\rm gen}))$. Since $\mathfrak h^{\rm gen}$ is $W$-invariant, its image under the finite quotient $q$ is open and has inverse image $\mathfrak h^{\rm gen}$. Thus $Z^{\rm gen}=a^{-1}(M^{\rm gen})$. Surjectivity of $f$ implies that $M^{\rm gen}$ is nonempty. Since $M$ is irreducible, the open subset $M^{\rm gen}\subset M$ is dense. Consider the unreduced fibre product $\widetilde Z=M\times_{\mathfrak h/\!\!/W}\mathfrak h$. The projection $\widetilde Z\to M$ is finite flat, being a base change of $q$, and hence is open. Reduction does not change the underlying topological space, so $a:Z\to M$ is also open. For every nonempty open subset $O\subset Z$, the nonempty open subset $a(O)\subset M$ meets $M^{\rm gen}$. Consequently $O$ meets $a^{-1}(M^{\rm gen})=Z^{\rm gen}$, proving that $Z^{\rm gen}$ is dense in $Z$.

By \Cref{prop:geometry}, $\pi$ is an isomorphism over $Z^{\rm gen}$. Its inverse image $f^{-1}(\mathfrak h^{\rm gen})$ is dense in $\mathcal X$, because $f$ is smooth and therefore open. Hence $\pi$ is an isomorphism between the dense open subsets $f^{-1}(\mathfrak{h}^{\rm gen})\rightarrow Z^{\rm gen}$. Moreover, its image is closed by properness (\Cref{prop:geometry}) and contains $Z^{\rm gen}$, so $\pi$ is surjective. Consequently $\pi$ is birational.

Now apply \cite[Proposition~10.1.5]{DGT} to \Cref{prop:geometry}. The mechanism is as follows. The constant complex on $\mathcal X$ is locally acyclic over $\mathfrak h$, because $f$ is smooth. That is, for any smooth curve $C\subset\mathfrak{h}$ and the base-change
\[\begin{tikzcd}
	{\mathcal{X}_C} & {Z_C} & C \\
	{\mathcal{X}} & Z & {\mathfrak{h}}
	\arrow["{\pi_C}", from=1-1, to=1-2]
	\arrow["{f_C}", bend left=30, from=1-1, to=1-3]
	\arrow[from=1-1, to=2-1]
	\arrow["\lrcorner"{anchor=center, pos=0.125}, draw=none, from=1-1, to=2-2]
	\arrow["{p_C}", from=1-2, to=1-3]
	\arrow[from=1-2, to=2-2]
	\arrow["\lrcorner"{anchor=center, pos=0.125}, draw=none, from=1-2, to=2-3]
	\arrow["\imath", from=1-3, to=2-3]
	\arrow["\pi", from=2-1, to=2-2]
	\arrow["f"', bend right=30, from=2-1, to=2-3]
	\arrow["p", from=2-2, to=2-3]
\end{tikzcd}\]
the vanishing cycle sheaf $\varphi_{f_C}\BoQ_{\mathcal{X}}=0$.
Proper pushforward preserves this property since the vanishing cycle sheaf functor commutes with proper maps. By the decomposition theorem, $\pi_*\BoQ[D]$ is a direct sum of shifted intersection complexes. On $Z^{\rm gen}$ it is the shifted constant sheaf since $\pi$ is an isomorphism over this open subset. Any additional summand would consequently have support over $\mathfrak h\setminus\mathfrak h^{\rm gen}$. Pulling back to a curve transverse to a general point of its support would produce a nonzero complex supported on the special fibre; its vanishing cycles could not vanish. This contradicts local acyclicity. Thus only the middle extension from $Z^{\rm gen}$ remains, which proves the assertion: $\pi$ is small and $\pi_*\BoQ_{\mathcal X}[D]\cong\IC_Z$.
\end{proof}

\subsection{The specialised Springer action}
\label{subsection:springer-action}
We explain how to obtain the Springer action on the cohomology of $Y$. The $W$-action on $Z$ coming from the $W$-action on the $\mathfrak{h}$-factor (see \eqref{equation:definitionZ}) gives a $W$-equivariant structure on $\IC_Z$, determined by its action on the constant sheaf on $Z^{\rm gen}$ and then using intermediate extension. Then,
\[
 K=a_*\IC_Z=(a\pi)_*\BoQ_{\mathcal X}[D]
\]
has a Weyl action extending the deck action on the generic ordered-spectrum covering.

Let
\[
 M_{\rm nil}=(\chi^{-1}(0))_{\mathrm{red}},\qquad\rho:Y\longrightarrow M_{\rm nil}
\]
be the central restriction of $a\pi$. It is important that $M_{\rm nil}$ allows \emph{nilpotent} moment maps. It is larger than the zero-moment moduli space $M_0=\mu^{-1}(0)/\!\!/G$. The reduced fibre of $Z$ over zero is $M_{\rm nil}$.  Indeed, since the reduced fibre of $q\colon \mathfrak{h}\rightarrow\mathfrak{h}/\!\!/W$ over $0\in \mathfrak{h}/\!\!/W$ is $0\in\mathfrak{h}$ we find $(Z\times_{\mathfrak{h}} \{0\})_{\mathrm{red}}\cong (M\times_{\mathfrak{h}/\!\!/W}\{0\})_{\mathrm{red}}=M_{\mathrm{nil}}$ . Consider the Cartesian square
\[
\begin{tikzcd}[column sep=large, row sep=large]
Y \arrow[r, "j"] \arrow[d, "\rho"']
    & \mathcal X \arrow[d, "a\pi"] \\
M_{\rm nil} \arrow[r, "i"']
    & M.
\arrow[ul, phantom, "\lrcorner", very near end]
\end{tikzcd}
\]
By proper base change, we obtain
\[
i^*K=i^*(a\pi)_*\BoQ_{\mathcal X}[D]\cong \rho_*j^*\BoQ_{\mathcal X}[D]\cong \rho_*\BoQ_Y[D].
\]
Consequently,
\[
 K|_{M_{\rm nil}}[-r]\cong\rho_*\BoQ_Y[D-r].
\]
Consequently the latter complex, and hence $\HO^*(Y,\BoQ)$, carries a Weyl group action. The central map $\rho:Y\to M_{\rm nil}$ is semismall, again by the relative smallness argument. Indeed, since $\pi$ is small, the complex
\[
P=\pi_*\BoQ_{\mathcal X}[D]\cong\IC_Z
\]
is perverse. Moreover, the smoothness of $f:\mathcal X\to\mathfrak h$ implies that $\BoQ_{\mathcal X}[D]$ is universally locally acyclic (ULA) over $\mathfrak h$. Since $\pi$ is proper, the same holds for $P$ as proper morphisms preserve ULA complexes.

Writing $i_0:Z_0\coloneqq p^{-1}(0)\hookrightarrow Z$ for the inclusion of the central fibre, local acyclicity, the perverse exactness of nearby cycles and the triangle connecting nearby cycles, vanishing cycles and restriction $i_0^*P$ imply that $i_0^*P[-r]$ is perverse. By proper base change in the diagram
\[\begin{tikzcd}
	Y & {Z_0} \\
	{\mathcal{X}} & Z
	\arrow[from=1-1, to=1-2]
	\arrow[from=1-1, to=2-1]
	\arrow["\lrcorner"{anchor=center, pos=0.125}, draw=none, from=1-1, to=2-2]
	\arrow["{i_0}", from=1-2, to=2-2]
	\arrow["\pi"', from=2-1, to=2-2]
\end{tikzcd}\]
we obtain
\[
i_0^*P[-r]\cong \pi_{0,*}\BoQ_Y[D-r],
\]
where $\pi_0:Y\to Z_0$ is the central restriction of $\pi$. Recall that $r=\dim \mathfrak{h}$.

Since $Y$ is smooth of dimension $D-r$, the perversity of this direct image implies that $\pi_0$ is semismall. As $Z_0\to M_{\rm nil}$ is finite, the same holds for $\rho$. This is also a consequence of \cite[Proposition~11.1.5(ii)]{DGT}.

We thus obtain the specialized Springer action on $\rho_*\BoQ_Y[D-r]$.

\begin{proposition}\label{prop:compare}
Under the restriction isomorphisms of \Cref{prop:constant}, the specialized action agrees with the generic action obtained by relabeling eigenvalues and then using the constant cohomology sheaf. In particular it agrees with the action in \cite[\S2.2]{HLRV}.
\end{proposition}
\begin{proof}
The equivariance of $\IC_Z$ gives morphisms between the pushforwards $p_*\IC_Z$ to $\mathfrak h$ covering the $W$-action $\sigma\mapsto w\sigma$. These pushforwards identify with those of the constant complex on $\mathcal X$, by \Cref{prop:small}, and their cohomology sheaves are constant by \Cref{prop:constant}. A morphism between constant sheaves is determined on any nonempty connected open. On the generic open the construction is exactly the deck action. Its value at zero is therefore the specialization of that same action.

The generic ordered family $f^{-1}(\mathfrak{h}^{\mathrm{gen}})$ is the one appearing in \cite{HLRV}: both parametrize a doubled-quiver representation with an ordered regular-semisimple eigenflag for its moment map. HLRV's canonical parallel transport and the identifications here are the identifications of the same constant cohomology local system.
\end{proof}

This proves constancy on the ordered base of eigenvalues $\mathfrak{h}$ and identifies the resulting representation on the unordered base. The kernel of $\pi_1(\mathfrak h^{\rm gen}/W)\to W$ acts trivially.

\subsection{Equivariance of tautological classes}
\label{subsection:equivariance-tautological-classes}
Let $P=P_{\dd}$, $R=P^W=R_{\dd}$, and let $\kappa:P\to\HO^*(\Y,\BoQ)$ be the map defined by the Chern classes of the tautological vector bundles $\mathcal{V}_{i,r}$, $i\in Q_0$ over $\Y$. Write $\XX\to\mathcal X$ for the scalar gerbe of stable representations in the total family, that is $\mathfrak{X}=\widetilde{\mathcal{X}}/G_{\dd^+}$ following notation of the proof of \Cref{prop:geometry}. It carries the vertex tautological vector bundles $\mathcal V_{i,r}$. The injective leg maps give line bundles
\[
 \mathcal L_{i,r}=\mathcal V_{i,r-1}/\mathcal V_{i,r}, \quad 1\le r\le d_i,\qquad \mathcal V_{i,d_i}=0.
\]
Their first Chern classes are denoted $x_{i,r}=c_1(\mathcal{L}_{i,r})$. Set $|\dd|=\sum_i d_i$ and
\[
 u_{\dd}=\frac1{|\dd|}\sum_i c_1(\mathcal V_{i,0}) =\frac1{|\dd|}\sum_{i,r}x_{i,r}.
\]
The class $u_{\dd}$ has weight one on the scalar gerbe. The classes $\bar x_{i,r}=x_{i,r}-u_{\dd}$ descend rationally to the algebraic variety $\mathcal X$. We have $\HO^*(\XX)=\HO^*(\mathcal X)\otimes\BoQ[u_{\dd}]$, and the same statement holds on every fibre of $f\colon\mathcal{X}\rightarrow\mathfrak{h}$ by base-change in this smooth family.

\begin{proposition}\label{prop:Wequiv}
The specialized Weyl group action of $W$ on $\HO^*(\Y)\cong \HO^*(Y)\otimes\BoQ[u_d]$ constructed in \Cref{subsection:springer-action} acts by graded ring automorphisms and satisfies
\[
 w(x_{i,r})=x_{i,w(r)},\qquad w(u_{\dd})=u_{\dd},\qquad w\kappa(f)=\kappa(wf)\quad(f\in P).
\]
\end{proposition}
\begin{proof}
At a generic spectrum $\sigma\in\mathfrak{h}^{\mathrm{gen}}$, $\mathcal L_{i,r}$ is the eigenline of the $r$-th ordered eigenvalue. Relabeling spectra permutes these lines. Their Chern classes in $\HO^*(\mathcal{X}_{\sigma})$ are restrictions of the classes on the whole family, so they are constant sections by \Cref{prop:constant}. We then apply \Cref{prop:compare} to specialize the identities to the fiber over zero which is $Y$.

At generic parameters the operators are induced by isomorphisms of varieties permuting the eigenvalues, followed by constant-family identification. They preserve cup products and the unit. These identities between maps of constant cohomology sheaves remain true at zero. We extend the action to the gerbe $\mathfrak{X}$ by fixing $u_{\dd}$. Poincar\'e duality in the smooth family gives the corresponding action on Borel--Moore homology and preserves the fundamental classes.
\end{proof}

\subsection{Kirwan surjectivity on the quiver variety of the extended quiver}
The following is a direct consequence of \cite{MN}. For convenience, we use the statement at the level of the stack $\Y$ instead of its moduli space $Y$.
\begin{proposition}\label{prop:Kirwan}
The tautological map
\[
 \kappa:P\longrightarrow\HO^*(\Y,\BoQ)
\]
is surjective.
\end{proposition}
\begin{proof}
Choose a vertex $\infty$ of $Q^+$ of dimension one. If it is the endpoint of a leg, it has no loops. Removing it and regarding its incident arrows as framing arrows is the Crawley-Boevey trick in reverse. We refer to \cite{CB} and \cite[\S3.2]{MN} for this description. The scalar quotient of the gauge group is
\[
 G_{\dd^+}/\Gm\cong\prod_{v\in Q^+_0, v\ne\infty}\GL_{d_v^+}, \qquad [(g_v)]\longmapsto(g_vg_\infty^{-1})_{v\ne\infty}.
\]
The equation at $\infty$ is implied by the other moment equations and their total trace. The quotient $Y$ is therefore a smooth framed Nakajima quiver variety.

If all original dimensions are at most one, there are no legs. We choose an original vertex with dimension $1$ as $\infty$. Any loops at this vertex give a factor of affine space in the preprojective moduli space: their maps are scalars and their commutators vanish, so that the loops at that vertex do not contribute to the moment map. The loops at the vertex of dimension $1$ therefore do not affect ordinary cohomology.

By \cite[Theorem~1.2]{MN}, the cohomology of $Y$ is generated by the Chern classes of the framed tautological bundles. On the stack these are the vector bundles $\CV_v\otimes\CV_\infty^{-1}$. The class $c_1(\CV_\infty)$ supplies the scalar gerbe variable. Thus $\HO^*(\Y)$ is generated by the Chern classes of all vertex bundles. Their flag filtrations express every such class as a polynomial in the $x_{i,r}$.
\end{proof}

\subsection{The sign convention}
We use the following convention regarding representations of the symmetric groups arising in Springer theory. For the classical Grothendieck resolution, deck action on the regular semisimple locus followed by restriction to the nilpotent cone gives the Springer action with trivial representation at a regular nilpotent and sign on the top cohomology over zero.

\subsection{The sign perverse sheaf}
\label{sign_per_sec}
The construction also produces a perverse sheaf
\[
 \mathcal S_{\dd}=\left(\rho_*\BoQ_Y[2p_{Q^+}(\dd^+)]\right)^{\sgn}
\]
on $M_{\rm nil}$. Its support is contained in the zero-moment subspace $M_0$, by \cite[Theorem~1.3.2(ii), \S9.4]{DGT} and its derived global sections recover the Kac polynomials of $Q$. This support statement does not by itself identify $\mathcal S_{\dd}$ with the BPS sheaf, although it is expected that $\mathcal{S}_{\dd}$ indeed coincides with the BPS sheaf of $Y$. The proof below compares cohomology through the actual flag correspondence and the sign character formula. The current paper does not rely on the perverse sheaf $\mathcal S_{\dd}$.

\subsection{The sign-isotypic dimension formula}
\label{sign_isotypic_sec}
Set $p^+=p_{Q^+}(\dd^+)=p_Q(\dd)+N$ (\Cref{lem:legdimension}). The following theorem is the cohomological interpretation of the coefficients of Kac polynomials that allowed Hausel--Letellier--Rodriguez-Villegas to prove the Kac positivity conjecture \cite{HLRV}.

\begin{theorem}[Hausel--Letellier--Rodriguez-Villegas]\label{thm:HLRV}
For the specialized Springer action,
\[
 A_{Q,\dd}(q)=\sum_i\dim\HO_c^{2i}(Y,\BoQ)^{\sgn}q^{i-p^+}.
\]
The odd cohomology of $Y$ vanishes.
\end{theorem}
The generic-fibre character formula is \cite[Corollary~1.6 and Theorem~2.1]{HLRV}. Apply \Cref{prop:constant,prop:compare} to transport it to $Y$. Odd vanishing also follows directly from \Cref{prop:Kirwan}, since all Chern classes have even degree. If $Y$ is empty, \Cref{rem:emptyfamily} reduces the formula to the empty generic fibre and in that case, both sides are $0$.

\begin{corollary}\label{cor:signcharacter}
In the normalised CoHA grading \eqref{eq:cohgrading}, we can identify the Hilbert series
\[
 \Hilb\bigl(\HO^*(\Y)^{\sgn}\cap[\Y]\bigr)=\frac{A_{Q,\dd}(t^{-1})}{1-t}.
\]
\end{corollary}
\begin{proof}
The complex dimension of $Y$, when nonempty, is by a direct calculation $2p^+$. Poincar\'e duality pairs $\HO_c^{2i}(Y)$ with $\HO^{4p^+-2i}(Y)$ and the sign character is self-dual. Setting $j=2p^+-i$ in \Cref{thm:HLRV} gives
\[
 \sum_j\dim\HO^{2j}(Y)^{\sgn}t^j=t^{p^+}A_{Q,\dd}(t^{-1}).
\]
The stack $Y$, being a $\mathbf{G}_{\mathrm{m}}$-gerbe over $Y$, has complex dimension $2p^+-1$, and its fundamental class has normalised degree $-2p^+$. Capping with that class therefore removes the factor $t^{p^+}$. The invariant scalar variable contributes $1/(1-t)$. If $Y$ is empty, the HLRV formula (\Cref{thm:HLRV}) gives zero on both sides.
\end{proof}

\section{Proof of the tautological generation of affinized BPS cohomology}\label{sec:proof}
\subsection{Alternating polynomial lifts}
We first give an elementary lemma.
\begin{lemma}\label{lem:alternating}
Let $P=\BoQ[x_{i,r}:i\in Q_0, 1\leq r\leq d_i]$, with $W=\prod_i\mathfrak S_{d_i}$, and let $\kappa:P\twoheadrightarrow H$ be a $W$-equivariant surjective map of graded algebras. Then
\[
 H^{\sgn}=\kappa(\Delta_{\dd}P^W).
\]
\end{lemma}
\begin{proof}
Let $h\in H^{\sgn}$, and choose a polynomial $f$ lifting it. The alternating projector defined by $\operatorname{Alt}(f)
=\frac1{|W|}\sum_{w\in W}\sgn(w)wf$ gives another lift. An alternating polynomial vanishes when two variables in the same block coincide. Thus it is divisible by each $x_{i,r}-x_{i,s}$, $i\in Q_0$. These linear polynomials are pairwise nonassociate irreducibles, so unique factorisation in $P$ gives $\operatorname{Alt}(f)=\Delta_{\dd}g$. Since both $\operatorname{Alt}(f)$ and $\Delta_{\dd}$ transform by sign, $g\in P^W$. Conversely, every $\Delta_{\dd}g$ with $g\in P^W$ is alternating.
\end{proof}

\subsection{The dimension argument}
\label{subsection:dimension-argument}
\begin{proof}[Proof of \Cref{thm:main} at zero stability]
Let $T_{\dd}=R_{\dd}\cap v_{\dd}$.
By \Cref{prop:tautprimitive},
\[
 T_{\dd}\subseteq\widehat{\Fg}_{\Pi_Q,\dd}.
\]
By \Cref{prop:vfc},
\[
 \Phi(T_{\dd})=\kappa(\Delta_{\dd}R_{\dd})\cap[\Y].
\]
By applying Kirwan surjectivity, \Cref{prop:Wequiv}, and \Cref{lem:alternating}, we can identify this image with $\HO^*(\Y)^{\sgn}\cap[\Y]$. Equations \eqref{eq:affKac} and \Cref{cor:signcharacter} then identify its graded dimension with that of $\widehat{\Fg}_{\Pi_Q,\dd}$.

Since $\Phi$ preserves the virtually shifted grading (\Cref{equation:flag-map}), in every degree $k\in\BoZ$ we have
\[
 \dim T_{\dd}^k\leq\dim\widehat{\Fg}_{\Pi_Q,\dd}^{\,k} =\dim\Phi(T_{\dd})^k\leq\dim T_{\dd}^k.
\]
All dimensions are finite. Thus every inequality is an equality, and $T_{\dd}=\widehat{\Fg}_{\Pi_Q,\dd}$.
\end{proof}

In \Cref{sec:Lie}, we give a different way to prove \Cref{thm:main} that does not rely in the HLRV formula, but instead on the identification of the flag map with an iterated Lie bracket in the cohomological Hall algebra, \Cref{prop:Heckeflag}. This allows one to prove that the flag map is injective on the BPS cohomology $\widehat{\Fg}_{\Pi_Q,\dd}$, using \Cref{prop:Liesign} and restriction to the semistable part, which does not change the BPS cohomology. Then, using that the flag map sends the BPS cohomology inside the sign-invariant part of the target, we obtain the same chain of inequality as in the proof of \Cref{thm:main} above, concluding the proof.

\begin{proof}[Proof of \Cref{thm:signintro}]
The preceding proof shows that $\Phi$ maps $\widehat{\Fg}_{\Pi_Q,\dd}=T_{\dd}$ onto the sign isotypic component in $\HO^*(\Y,\BoQ)$. The equality of graded dimensions also shows that its kernel on this space is zero. Linearity and the image of the generator were proved in \Cref{lem:Phi-linear,prop:vfc}.
\end{proof}

\subsection{Stability, the lowest degree, and the smooth case}
\label{subsection:proofs}
\begin{proof}[Proof of \Cref{thm:main} for arbitrary stability]
Restrict the equality of \Cref{thm:main} at zero stability using \Cref{cor:Todaaff}. By \Cref{lem:restrictionvirtual}, restriction sends its generator $v_{\dd}$ to $v_{\dd}^{\theta}$ and commutes with all tautological cap products. This gives the required equality of subspaces, by the commutative diagram:
\[\begin{tikzcd}
	{R_{\dd}\cdot v_{\dd}^{\theta}} & {\widehat{\mathfrak{\Fg}}_{\Pi_Q,\dd}^{\theta}} & {\HO^{\mathrm{BM}}_*(\FM_{\Pi_Q,\dd}^{\theta})} \\
	{R_{\dd}\cdot v_{\dd}} & {\widehat{\mathfrak{\Fg}}_{\Pi_Q,\dd}} & {\HO^{\mathrm{BM}}_*(\FM_{\Pi_Q,\dd})}
	\arrow[from=1-1, to=1-2]
	\arrow[from=1-2, to=1-3]
	\arrow["\simeq", from=2-1, to=1-1]
	\arrow["\simeq"', from=2-1, to=2-2]
	\arrow["\simeq"', from=2-2, to=1-2]
	\arrow[from=2-2, to=2-3]
	\arrow[from=2-3, to=1-3]
\end{tikzcd}\]
since the map $R_{\dd}\cdot v_{\dd}^{\theta}\rightarrow \widehat{\mathfrak{\Fg}}_{\Pi_Q,\dd}^{\theta}$ must then be an isomorphism.
\end{proof}

\begin{proof}[Proof of \Cref{cor:lowest}]
If $A_{Q,\dd}=0$, the primitive space is zero by \eqref{eq:affKac}. Thus $v_{\dd}=0$. Otherwise, \Cref{thm:main} implies $v_{\dd}\ne0$. By \Cref{lem:vdim}, its degree is $-2p_Q(\dd)$. Monicity of the Kac polynomial says that this is the one-dimensional lowest-degree piece of affinized BPS. Since the determinant class $u$ has positive degree $2$, that piece already lies in unaffinized BPS $\Fg_{\Pi_Q,\dd}\subset\widehat{\Fg}_{\Pi_Q,\dd}\cong \Fg_{\Pi_Q,\dd}\otimes\BoQ[u]$. Restriction and Toda's theorem (\Cref{thm:Toda}) prove the stable versions for any stability parameter $\theta$.
\end{proof}

\begin{proof}[Proof of \Cref{cor:smooth}]
If the semistable stack is empty, there is nothing to prove. Otherwise it is smooth by \Cref{lem:stablefund}. The semistable PBW theorem has no nonzero decomposable summand in dimension $\dd$: if semistable objects of nonzero dimensions $\dd_1,\ldots,\dd_m\in\Lambda_\theta$, $m>1$, existed with sum $\dd$, their direct sum would be strictly semistable in dimension $\dd$. This contradicts the hypothesis. Hence
\[
 \HO^{\mathrm{BM}}_*(\FM_{\Pi_Q,\dd}^{\theta},\BoQ^{\vir})=\CA_{\Pi_Q,\dd}^{\theta}=\widehat{\Fg}_{\Pi_Q,\dd}^{\theta}.
\]
The reduced virtual class is the smooth fundamental class by \Cref{lem:stablefund}. The generation theorem (\Cref{thm:main}) and Poincar\'e duality now give a surjection $R_{\dd}\twoheadrightarrow\HO^*(\FM_{\Pi_Q,\dd}^{\theta},\BoQ)$.

Finally use the splitting \eqref{eq:gerbesplitting}. Modulo $u_{\dd}$ the surjection becomes a surjection onto the cohomology of the coarse moduli space. Removing the scalar character from the Chern classes gives precisely the rational tautological classes of \Cref{subsec:gerbe} over the moduli space. This then proves that $\HO^*(\mathcal{M}_{\Pi_Q,\dd}^{\theta},\BoQ)$ is generated by tautological classes.
\end{proof}

\subsection{A presentation by tautological relations}
Let $J_{\Theta}=\ker(P_{\dd}\xrightarrow{\kappa}\HO^*(\Y))$. It is a homogeneous $W$-stable ideal.

\begin{corollary}\label{cor:presentation}
Let
\[
 I_{\dd} =\{f\in R_{\dd}:\Delta_{\dd}f\in J_{\Theta}\}.
\]
Then
\[
 R_{\dd}/I_{\dd}\xrightarrow{\ \sim\ }\widehat{\Fg}_{\Pi_Q,\dd}, \qquad \overline f\longmapsto f\cap v_{\dd},
\]
where the unit on the left is placed in degree $-2p_Q(\dd)$. In particular $I_{\dd}$ is independent of the chosen generic stability parameter $\Theta$ lying in the positive-leg chamber (i.e. whose coordinates at leg-vertices are positive) given by \Cref{lem:chamber}.
\end{corollary}
\begin{proof}
Surjectivity is \Cref{thm:main}. By \Cref{thm:signintro}, an element $f\cap v_{\dd}$ vanishes if and only if its image $\Phi(f\cap v_{\dd})$ by the flag correspondence vanishes. By \Cref{prop:vfc}, this is equivalent to $\kappa(\Delta_{\dd}f)=0$. Thus $I_{\dd}$ is the annihilator of the intrinsic class $v_{\dd}$ in $R_{\dd}$, which also proves chamber independence.
\end{proof}

For $A_{Q,\dd}\ne0$, this quotient is a finite free $\BoQ[u_{\dd}]$-module of rank $A_{Q,\dd}(1)$. With the unshifted grading of the polynomial ring, its Hilbert series is
\[
 \Hilb(R_{\dd}/I_{\dd})=\frac{t^{p_Q(\dd)}A_{Q,\dd}(t^{-1})}{1-t}.
\]
After choosing the determinant normalisation,
\[
 \Fg_{\Pi_Q,\dd}\cong R_{\dd}/(I_{\dd}+(u_{\dd}))
\]
as graded vector spaces, again placing the unit in degree $-2p_Q(\dd)$. This description equips affinized BPS with a cyclic tautological-module presentation. The identification of its Lie bracket in that quotient ring is an interesting question. Moreover, this isomorphism provides a ``ringification'' of the BPS cohomology $\Fg_{\Pi_Q,\dd}$ that would be interesting to study, and try to connect with Hausel's quest of a ringification of intersection cohomology \cite{Hausel2024} in some representation-theoretic contexts. In particular, studying the indecomposable summands of this algebra could provide new insights on the BPS cohomology.

\section{Sign isotypic component from the BPS Lie algebra}\label{sec:Lie}
\subsection{The real-root action}

We explain an other Lie-theoretic description of a Weyl group action on quiver varieties, deduced from \cite{DHSM} and the description of BPS Lie algebras as generalised Kac--Moody Lie algebras.

Every vertex added to $Q$ to form the quiver with legs $Q^+$ is loop-free. The real simple roots supported on a leg of length $d_i-1$ form a root system of type $A_{d_i-1}$. The generalised Kac--Moody Lie algebra description of the BPS Lie algebra \cite{DHSM} therefore provides an integrable action of
\[
 \mathfrak s=\bigoplus_{i:d_i>1}\mathfrak{sl}_{d_i}.
\]
on the full BPS Lie algebra $\Fg_{\Pi_{Q^+}}$, which is identified with a positive part of a Borcherds--Kac--Moody Lie algebra. We use the adjoint action in the full generalised Kac--Moody algebra whose positive part is $\Fg_{\Pi_{Q^+}}$.

Fix the original dimension vector $\dd$, and set
\[
 \mathcal V_{\dd}=\bigoplus_{\ee\in\BoN^{Q_0^+\setminus Q_0}}                      \Fg_{\Pi_{Q^+},(\dd,\ee)}.
\]
The space $\mathcal{V}_{\dd}$ is stable under the positive and negative real-root operators generating $\mathfrak{s}$. Indeed, their action changes only added coordinates, i.e. the dimension vector $\ee$, while an original coordinate stays positive, and this action preserves roots by definition. A root obtained in this way through the action of $\mathfrak{s}$ cannot become a negative root, since the original vertex remains positive. Each real-root $\mathfrak{sl}_2$ copy acts locally nilpotently, by integrability, as a consequence of the generalised Serre relations involved in the definition of the positive part $\Fg_{\Pi_{Q^+}}$, see \cite{DHSM}.

\subsection{Polynomial representations of type $A$ Lie algebras}

We recall the identification between dimension vectors along the legs of $Q^+$ and weights of representations of type $A$ Lie algebras.

Consider the leg attached to a vertex $i$ of $Q$, and let $n=d_i$. We label the vertices of this leg by $1,\ldots,n-1$, starting from the vertex adjacent to $i$. Given a dimension vector $(e_1,\ldots,e_{n-1})$ along this leg, we set $e_0=n$ and $e_n=0$.

Let $\varepsilon_1,\ldots,\varepsilon_n$ denote the standard basis of the weight lattice of the diagonal torus of $\mathrm{GL}_n$. We associate to the leg dimension vector the weight
\begin{equation}\label{eq:GLweight}
\mathrm{wt}(\mathbf e)=\sum_{r=1}^n (e_{r-1}-e_r)\varepsilon_r.
\end{equation}
In other words, writing $\mathrm{wt}(\mathbf e)=(a_1,\ldots,a_n)$, we have $a_r=e_{r-1}-e_r$. Conversely, the leg dimensions are recovered from the weight by
\[
e_r=a_{r+1}+\cdots+a_n,\qquad 1\leq r\leq n-1.
\]
In particular, $\sum_{r=1}^n a_r=n$.

Recall that the weights of a homogeneous polynomial representation of $\mathrm{GL}_n$ of degree $n$ belong to
\[
\left\{(a_1,\ldots,a_n)\in\mathbb Z_{\geq 0}^n\;\middle|\; \sum_{r=1}^n a_r=n\right\}.
\]
Under the correspondence above, the nonnegativity of these weights is equivalent to
\[
n=e_0\geq e_1\geq\cdots\geq e_{n-1}\geq e_n=0.
\]

We now describe how the simple-root operators act under this identification. For $1\leq r\leq n-1$, increasing $e_r$ by one changes the corresponding weight by
\[
\mathrm{wt}(\mathbf e+\mathbf 1_r)=\mathrm{wt}(\mathbf e)-\varepsilon_r+\varepsilon_{r+1},
\]
where $\mathbf 1_r$ denotes the dimension vector supported at the $r$-th vertex of the leg. Writing $\alpha_r=\varepsilon_r-\varepsilon_{r+1}$ for the $r$-th simple root of $\mathfrak{sl}_n$, this is precisely the weight shift associated with the negative simple-root operator
\[
f_r=E_{r+1,r}\in\mathfrak{sl}_n.
\]
Thus, with our conventions, the simple-root generator associated with the $r$-th vertex of the leg is identified with $f_r$. In particular, the positive-root Lie algebra associated with the leg is identified with the lower-triangular nilpotent subalgebra of $\mathfrak{sl}_n$.

Recall that a finite-dimensional representation of $\GL_n$ is \emph{polynomial} if its matrix coefficients are polynomial functions of the entries of $\GL_n$, or equivalently, if the action extends to the monoid $\mathrm{Mat}_n$ of all $n\times n$ matrices. It is \emph{homogeneous of degree $m$} if scalar matrices $t\mathrm{Id}_n$ act by multiplication by $t^m$. The irreducible polynomial representations are the Schur modules $\mathbb S_\lambda(\mathbb C^n)$, indexed by partitions $\lambda=(\lambda_1\geq\cdots\geq\lambda_n\geq0)$. Their homogeneous degree is $|\lambda|=\sum_r\lambda_r$.

\begin{lemma}\label{lem:polynomial}
In each cohomological degree, $\mathcal V_{\dd}$ extends to a polynomial representation of $\prod_i\GL_{d_i}$, homogeneous of degree $d_i$ in the $i$-th factor.
\end{lemma}

\begin{proof}
Consider a single factor $\mathfrak{sl}_n$, where $n=d_i$. Recall that its weights are represented by the tuples $(a_1,\ldots,a_n)$ defined in \eqref{eq:GLweight}, which satisfy $\sum_r a_r=n$.

Since the representation is integrable, its set of weights is invariant under the Weyl group $S_n$, acting by permutations of the coordinates $a_r$. For every occurring weight, the last coordinate satisfies $a_n=e_{n-1}\geq0$, since it is a dimension at the last vertex of the leg. By Weyl group invariance, any coordinate can be moved to the last position, and hence $a_r\geq0$ for every $r$.

Thus all occurring weights belong to the finite set
\[
\left\{(a_1,\ldots,a_n)\in\mathbb Z_{\geq0}^n\;\middle|\; \sum_{r=1}^n a_r=n\right\}.
\]
Since each weight space in a fixed cohomological degree is finite-dimensional, the representation in that degree is finite-dimensional, and hence semisimple.

The highest weights of its irreducible summands are therefore weakly decreasing nonnegative $n$-tuples of sum $n$, namely partitions of $n$. By the classification of finite-dimensional representations of $\GL_n$, each such summand extends to the polynomial representation $\mathbb S_\lambda(\mathbb C^n)$, homogeneous of degree $n$. The condition on the sum of the weights determines the extension from $\mathfrak{sl}_n$ to $\GL_n$ uniquely.

Applying the same argument to every leg, and using that the actions of the different factors commute, gives the required polynomial representation $\mathcal{V}_{\dd}$ of $\prod_i\GL_{d_i}$.
\end{proof}

The weight with all added dimensions zero at the extra leg vertices is $(n,0,\ldots,0)$. The full-legged weight is $(1,\ldots,1)$, which has zero weight for $\mathfrak{sl}_n$.

\subsection{The determinant-one Weyl action}
Let $T\subset\SL_n$ be the diagonal maximal torus and $N(T)$ its normaliser. Since the representation is finite-dimensional, the $\mathfrak{sl}_n$-action integrates to an $\SL_n$-action. The normaliser $N_{\SL_n}(T)$ preserves the zero-weight space, and $T$ acts trivially on it. Consequently, the action of $N_{\SL_n}(T)$ factors through the Weyl group $W=N_{\SL_n}(T)/T\simeq S_n$.

More explicitly, if $e_r,f_r$ are the standard Chevalley generators, the simple reflection $s_r\in W$ is represented by the Tits operator
\[
\dot{s}_r=\exp(f_r)\exp(-e_r)\exp(f_r),
\]
as in \cite[\S 3.8, Lemma 3.8]{Kac1990}. These operators satisfy the braid relations, while $\dot{s}_r^2=\alpha_r^\vee(-1)\in T$. Since $T$ acts trivially on the zero-weight space, the restrictions of the Tits operators satisfy the Coxeter relations and therefore define a genuine action of $W$. This agrees with the action obtained from the normaliser.

\begin{lemma}\label{lem:signlemma}
Let $V$ be a homogeneous polynomial $\GL_n$-representation of degree $n$. Let $W=\mathfrak S_n$ act on $V_{(1,\ldots,1)}$ using determinant-one representatives. Put
\[
 D_n=E_{n-1}^{(1)}E_{n-2}^{(2)}\cdots E_1^{(n-1)},\qquad E_r^{(m)}=\frac{E_r^m}{m!},\qquad E_r=E_{r+1,r}\in\mathfrak{sl}_r.
\]
Then $D_n$ induces an isomorphism
\[
 V_{(n,0,\ldots,0)}\xrightarrow{\sim}V_{(1,\ldots,1)}^{\sgn}.
\]
For $n=1$ this is the identity map.
\end{lemma}
\begin{proof}
Decompose
\[
 V=\bigoplus_{\lambda\vdash n}M_\lambda\otimes S_\lambda(\BoC^n).
\]
The weight $(n,0,\ldots,0)$ occurs only for $\lambda=(n)$, and then with multiplicity one. Indeed, in the semistandard-tableau description all entries would have to be $1$, so no column can have length greater than one.

The multilinear weight $(1,\ldots,1)$ in $S_\lambda(\BoC^n)$, under permutation matrices, is the Specht representation of shape $\lambda$. This follows directly from the Young-symmetriser construction: the multilinear part of $(\BoC^n)^{\otimes n}$ has basis indexed by permutations, and applying the Young symmetriser gives the corresponding Specht module.

A permutation matrix has determinant $\sgn(w)$. To obtain a representative in $\SL_n$, multiply by a diagonal matrix with determinant $\sgn(w)$. On the multilinear weight, that diagonal correction acts by the product of its diagonal entries, hence by $\sgn(w)$. The Weyl representation on this weight is therefore
\[
 \sgn\otimes\operatorname{Sp}_\lambda.
\]
It contains sign exactly when $\lambda=(n)$, with multiplicity one.

It remains to check the map on this summand. In $S_{(n)}(\BoC^n)=\Sym^n(\BoC^n)$, write $E_r=z_{r+1}\partial_{z_r}$. Successive divided powers give
\[
 D_n(z_1^n)=n!\,z_1\cdots z_n\ne0.
\]
The map is thus an isomorphism on every multiplicity space $M_{(n)}$.
\end{proof}

\begin{proposition}\label{prop:Liesign}
The product of the leg operators induces an isomorphism
\[
 \Fg_{\Pi_Q,\dd}\xrightarrow{\sim}\Fg_{\Pi_{Q^+},\dd^+}^{\sgn_{\mathrm{Lie}}}.
\]
\end{proposition}
\begin{proof}
The extension by zero at the new vertices identifies $\Fg_{\Pi_Q,\dd}$ with the weight having all added coordinates zero in $\mathcal{V}_{\dd}$. We apply \Cref{lem:polynomial,lem:signlemma} to every factor. The full-legged dimension vector $\dd^+$ gives weight zero for every leg $\mathfrak{sl}_{d_i}$. We conclude using \Cref{lem:signlemma}.
\end{proof}

For affinization, use the class
\[
 u_{\dd}=\frac1{|\dd|}\sum_{i\in Q_0}c_1(\CV_i)
\]
on every summand with the original dimensions fixed. It has scalar weight one, so it gives a polynomial generator in the integrality theorem. Its multiples of BPS classes are primitive by the same spectral-support argument as in \Cref{prop:primitives}. The original tautological vector bundles at vertices of $Q$ are unchanged by the Hecke correspondences on the legs of the quiver, so this class commutes with those operators. Consequently \Cref{prop:Liesign} also holds for affinized BPS cohomology. Toda's comparison theorem (\Cref{thm:Toda}) transports the resulting sign summand to the smooth stable stack by restriction.

\subsection{Hecke operators and the flag correspondence}
We use the cohomological Hall algebra convention in which the CoHA product $a*b$ is defined using the stack of extensions with subobject of class $b$. Let $e_r$ be the degree-zero simple-root class at a newly added vertex of a leg. The operator $E_r$ on the old root spaces is given by the adjoint action $\ad(e_r)$.

For an isolated loop-free vertex, the dimensional-reduction shuffle product gives
\begin{equation}\label{eq:dividedpower}
 e^{*m}=\varepsilon_m\,m!\,[\BG\GL_m],\qquad \varepsilon_m\in\{1,-1\}.
\end{equation}
Indeed, the corresponding one-loop critical shuffle formula has no root-factor denominator or numerator. Multiplication is the sum over shuffles, and multiplying $m$ unit classes gives $m!$ copies of the unit class in rank $m$. The Euler-form and dimensional-reduction sign conventions can only change this by the sign $\varepsilon_m$. Thus action by cohomological Hall Lie bracket of divided powers of \Cref{lem:signlemma} add the whole space at a new vertex at once, up to a nonzero scalar. This is the only normalisation used below.

\begin{proposition}\label{prop:Heckeflag}
After restriction to the stack of semistable representations for a stability parameter $\Theta$ in a positive-leg chamber (given for example by \Cref{lem:chamber}), the product of the divided adjoint operators in \Cref{prop:Liesign} agrees with the flag correspondence map $\Phi$ \eqref{equation:flag-map}, up to a nonzero scalar depending on Hall normalisations.
\end{proposition}
The normalizations are irrelevant to our study, since the only relevant information is the image of the flag map $\Phi$, which only requires its identification up to a scalar.
\begin{proof}
For a single operator,
\[
 \frac{\ad(e)^m}{m!}(a)=\sum_{q=0}^m\frac{(-1)^q}{(m-q)!\,q!}                  e^{*(m-q)}*a*e^{*q}.
\]
Any term with $q>0$ parametrises an object having a nonzero subobject supported on added vertices. Such an object is unstable for a stability parameter $\Theta$ positive on these vertices, as chosen in \Cref{lem:chamber}. The same argument applies after successive adjoint operations: a term having a nonzero total right-hand leg factor still has a leg-supported subobject. Only the completely left-ordered product survives on $\Y$.

By \eqref{eq:dividedpower}, each surviving divided power is the fundamental class $[\BG\GL_m]$, up to the fixed sign convention that is not important for us. We add the vertices in order from the original vertex outwards according to the formula for $D_n$ in \Cref{lem:signlemma}. At each step, the previous representation has zero space at the new vertex, while the new quotient is supported entirely at that vertex. The subrepresentation condition forces the reverse map $b_r$ to vanish. The inward map $a_r$ is arbitrary, and becomes injective on the chosen stable open.

After quotienting by changes of basis at the new vertex, this is the Grassmann bundle selecting a subspace of the preceding space. Iterating gives the full flag bundle $p$. All the $b_r$ vanish, which is exactly $\eta=0$ in the flag model (\Cref{subsection:flag-model}). Since the two factors have disjoint vertex supports at each step, there are no common-vertex off-diagonal matrix blocks: the pullback appearing in the definition of the cohomological Hall induction on this locus is the smooth pullback for the Grassmann bundle. The evaluation map is its closed inclusion in the target. The composite is therefore $j_*p^!$, followed by stable restriction to the smooth stable stack $\Y$. This is precisely the definition of the flag map \eqref{equation:flag-map}.
\end{proof}

Together with \Cref{lem:signlemma}, this gives an independent proof that the flag map \Cref{equation:flag-map} is injective on affinized BPS, and this could replace the use of the HLRV formula to prove the tautological generation theorem, as explained in \Cref{subsection:dimension-argument} after the proof of \Cref{thm:main}.

\subsection{The two constructions of Weyl group actions}

We have constructed two a priori different Weyl group actions on $\HO^*(Y,\BoQ)$: the Lie-theoretic action, obtained from the integrable $\mathfrak{s}$-representation through the Tits representatives (\Cref{lem:polynomial}), and the geometric action, defined by using a smooth family of quiver varieties (\Cref{sec:family}). Although these constructions have different origins, they are closely related.

Indeed, \Cref{prop:Heckeflag} and \Cref{thm:signintro} show that their sign-isotypical subspaces coincide. More precisely, both are identified with the image of the same BPS space $\Fg_{\Pi_Q}$ under the same flag correspondence $\Phi$ \eqref{equation:flag-map}, up to the scalar normalisation of the corresponding maps. Thus, while we do not identify the two Weyl group actions on the entire cohomology of the smooth quiver variety $\HO^*(Y,\BoQ)$, their sign-isotypical components agree as subspaces of the cohomology of the smooth target. This is precisely the identification needed for the tautological generation theorem (\Cref{thm:main}).

It is natural to expect that the two actions actually coincide. Let us explain how this comparison should work. Consider a simple reflection $s_r$, corresponding to the transposition of two consecutive steps of a flag. Writing $\CV_0,\ldots,\CV_n$ for the tautological bundles corresponding to a flag (i.e. to a leg of the extended quiver $Q^+$), the expected transformation of their classes in $K$-theory is
\[
[\CV_r']=[\CV_{r-1}]+[\CV_{r+1}]-[\CV_r],
\]
while the other tautological bundles remain unchanged. In particular, if
\[
x_r=c_1(\CV_{r-1})-c_1(\CV_r),\qquad x_{r+1}=c_1(\CV_r)-c_1(\CV_{r+1}),
\]
denote the Chern roots associated with two consecutive steps of the flag, then
\[
x_r'=c_1(\CV_{r-1})-c_1(\CV_r')=c_1(\CV_r)-c_1(\CV_{r+1})=x_{r+1},
\]
and similarly $x_{r+1}'=x_r$. Thus the expected geometric reflection acts on the tautological Chern roots by the corresponding transposition, in accordance with the standard Weyl group action.

To establish equality of the full actions, one would need to identify the geometric reflection or continuation operator with the Lie-theoretic Tits operator
\[
\dot{s}_r=\exp(E_r)\exp(-F_r)\exp(E_r)
\]
on the relevant weight-zero space. In particular, this requires a comparison of the geometric correspondences defining the two operators, beyond their effect on tautological classes. We do not pursue this comparison here, since the identification of the sign-isotypical subspaces is sufficient for our purposes. We expect to return to the comparison of the full actions in future work.

\section{Central deformations and generation for 2-Calabi--Yau categories}\label{sec:deformation}
All results and construction of this paper have natural extensions to deformed preprojective algebras \cite{CBH98}. In this section, we explain how they can be obtained by adapting the arguments of the undeformed case.

\subsection{The deformed preprojective algebra}
For $\lambda\in\BoC^{Q_0}$, we define following \cite{CBH98} the deformed preprojective algebra
\[
 \Pi_Q^\lambda=\BoC\overline Q\Big/\Big(\sum_{a\in Q_1}[a,a^*]-\sum_{i\in Q_0}\lambda_i e_i\Big).
\]
A representation of $\Pi_{Q}^{\lambda}$ of dimension $\dd$ can exist only if $\lambda\cdot\dd=0$, by taking traces. We let $\theta\in\mathbf{R}^{Q_0}$ be a King stability parameter such that $\theta\cdot \dd=0$. Under this condition, we define the stack of $\theta$-semistable $\dd$-dimensional representations of $\Pi_Q^{\lambda}$ by
\[
 \FM_{\Pi_Q^\lambda,\dd}^{\theta}:=[\mu^{-1}((\lambda_i\id_{V_i})_i)^\theta/G_{\dd}].
\]
Then $\mu-\lambda$ has values in $\mathfrak{gl}_{\dd}^{\,0}$, and we define the reduced virtual fundamental class
\[
 v_{\dd,\lambda}^{\theta}=[\FM_{\Pi_Q^\lambda,\dd}^{\theta}]^{\rd}
\]
by the same refined zero-section construction as in \Cref{sec:virtual}.

We extend $\lambda$ to a deformation parameter $\lambda^+\in \BoC^{Q_0^+}$ for $Q^+$ by putting zero at every new vertex in $Q_0^+\setminus Q_0$, and set
\[
 \Y_\lambda:=\FM_{\Pi_{Q^+}^{\lambda^+},\dd^+}^{\Theta}.
\]
We denote its coarse moduli space by $Y_\lambda$. These are the $\lambda$-deformed versions of the stacks $\Y$ and $Y$ considered in the undeformed case, see \eqref{equation:smoothstable}. The leg relations remain those of \eqref{eq:legrelations}, since the deformation parameter $\lambda^+$ vanishes on legs, and the original-vertex equation becomes
\[
 \mu(x)+\eta=\lambda.
\]
In particular the proof of injectivity for stability parameters in a positive-leg chamber (\Cref{lem:injectivity}), the classical and derived flag models (\Cref{subsection:flag-model} and \Cref{subsec:derivedflag}), and the zero-section calculation for defining virtual fundamental classes (\Cref{sec:virtual}) all apply without change.

We have the $\lambda$-deformed flag correspondence
\[
 \FM_{\Pi_Q^{\lambda},\dd}  \xleftarrow{\ p\ }\mathfrak{F}^{\lambda} \xrightarrow{\ j\ }\Zleg^{\lambda}.
\]
Here, writing
\[
 \lambda_{\dd}=(\lambda_i\id_{V_i})_{i\in Q_0}, \qquad X_{\dd}^{\lambda}=\mu^{-1}(\lambda_{\dd}),
\]
we set
\[
 \mathfrak{F}^{\lambda}= [X_{\dd}^{\lambda}\times\CB/G_{\dd}], \qquad \CB=\prod_{i\in Q_0}\Flag(V_i).
\]
Thus $\mathfrak{F}^{\lambda}$ parametrises a $\Pi_Q^{\lambda}$-module $M$ of dimension vector $\dd$, together with a decreasing full flag
\[
 M_i=F_{i,0}\supset F_{i,1}\supset\cdots \supset F_{i,d_i}=0, \qquad \rk F_{i,r}=d_i-r,
\]
at each vertex. It is a smooth projective flag bundle of relative dimension $N=\sum_i\binom{d_i}{2}$ over $\FM_{\Pi_Q^{\lambda},\dd}$.

The stack $\Zleg^{\lambda}$ is the open substack of $\FM_{\Pi_{Q^+}^{\lambda^+},\dd^+}$ on which all inward leg maps are injective. The deformed version of \Cref{prop:flagmodel} identifies it with
\begin{equation}\label{equation:deformedZleg}
 \Zleg^{\lambda}\cong \bigl[ \{(x,F,\eta)\in V_{\dd}\times T^*\CB:             \mu(x)+\eta=\lambda_{\dd}\} /G_{\dd} \bigr],
\end{equation}
where the cotangent coordinate satisfies $\eta_i(F_{i,r})\subseteq F_{i,r+1}$. Under this identification,
\[
 p(x,F)=[x], \qquad j(x,F)=(x,F,0).
\]
Equivalently, $j$ places $F_{i,r}$ at the added vertex $(i,r)$, takes the inclusions as inward leg maps, and sets all reverse leg maps to zero. It is the closed immersion defined by the vanishing of the tautological section $\eta$ of the rank-$N$ vector bundle $\CE^{\lambda}$ on $\Zleg^{\lambda}$, obtained by pulling back $[T^*\CB/G_{\dd}]\to[\CB/G_{\dd}]$ along the morphism forgetting $x$ and $\eta$, as in \Cref{subsection:the-correspondence} in the undeformed case.

Let
\[
 \iota_\lambda:  \Y_\lambda =\FM_{\Pi_{Q^+}^{\lambda^+},\dd^+}^{\Theta} \hookrightarrow\Zleg^\lambda
\]
be the open immersion of the stable locus for the chosen generic positive-leg stability parameter $\Theta$ (\Cref{lem:chamber}). We define the \emph{$\lambda$-deformed flag map} by
\[ \Phi_\lambda :=\iota_\lambda^*j_*p^!: \HO_m^{\BM}(\FM_{\Pi_Q^\lambda,\dd},\BoQ) \longrightarrow \HO_{m+2N}^{\BM}(\Y_\lambda,\BoQ).
\]
Here $p^!$ is smooth pullback along the flag bundle, $j_*$ is proper pushforward along the closed immersion, and $\iota_\lambda^*$ is restriction to the stable open. The shift by $2N$ is twice the relative dimension of $p$. Since by \Cref{lem:legdimension}
\[
 (\dd^+,\dd^+)_{Q^+}=(\dd,\dd)_Q-2N,
\]
the map $\Phi_\lambda$ preserves the normalised CoHA grading. It is $R_{\dd}$-linear: the universal bundles at the original vertices pull back to the same bundles on $\mathfrak F^\lambda$, and the assertion follows from the projection formula as in \Cref{lem:Phi-linear}.

We now have the $\lambda$-deformed version of \Cref{prop:vfc}.
\begin{proposition}\label{prop:deformedgeometry}
For a generic positive-leg stability parameter $\Theta$ (as in \Cref{lem:chamber}) there is an $R_{\dd}$-linear flag map $\Phi_\lambda$ satisfying
\[
 \Phi_\lambda(v_{\dd,\lambda})=(-1)^N\Delta_{\dd}\cap[\Y_\lambda].
\]
It preserves the normalised CoHA degree.
\end{proposition}
\begin{proof}
Use $\mu-\lambda+\eta$ in the reduced derived flag model  (\Cref{prop:flagmodel}), see \eqref{equation:deformedZleg}. The independent linear equations at leg vertices given by the corresponding summand of the moment map are unchanged, so \Cref{lem:linear-elimination} still applies. The derived zero locus of $\eta$ is the flag bundle over the reduced derived fibre $\mu=\lambda$. The proof of \Cref{prop:vfc} gives the result.
\end{proof}

The critical presentation also persists: replace the tripled potential by its $\lambda$-deformed version
\[
 W_\lambda=W-\sum_i\lambda_i\omega_i.
\]
Its critical equations describe $\Pi_Q^\lambda[\omega]$. We can then define the coproduct on the Borel--Moore homology $\CA_{\Pi_Q^{\lambda}}^{\theta}$ exactly as in the undeformed case for any deformation parameter $\theta$.

\subsection{The smooth deformed target}
The family of quiver varieties considered in \Cref{subsection:smooth-GIT-family} already settles the smooth geometry in the deformed case.

\begin{proposition}\label{prop:deformedtarget}
Let $\lambda\in\BoC^{Q_0}$ and $\dd\in\BoN^{Q_0}$ be such that $\lambda\cdot\dd=0$. The cohomology of the smooth stack $\Y_\lambda$ has a specialized $W$-action fixing the scalar class $u_{\dd}$. Its tautological map $\kappa_{\lambda}\colon P\rightarrow\HO^*(\Y_{\lambda})$ is surjective and $W$-equivariant, and
\[
 \HO^*(\Y_\lambda)^{\sgn}=\kappa_\lambda(\Delta_{\dd}R_{\dd}).
\]
The graded character of this sign isotypic component is the same as at $\lambda=0$ and is given by \Cref{cor:signcharacter}.
\end{proposition}
\begin{proof}
Put $\sigma_{i,r}=\lambda_i$ for every $r$. Under \eqref{eq:eigenparameters}, this is the central parameter equal to $\lambda$ at the original vertices and zero at the added vertices. Thus the fibre $\mathcal X_\sigma$ of $f$ at $\sigma$ (see \eqref{eq:diagram}) is $Y_\lambda$. Every element of the Weyl group $W$ fixes this parameter, so the equivariant intersection-complex construction in \Cref{subsection:springer-action} restricts to an action on its cohomology $\HO^*(Y_{\lambda},\BoQ)$. By \Cref{prop:constant,prop:compare}, the comparison with the cohomology of $Y$ is $W$-equivariant.

The flag quotient line bundles $\mathcal{L}_{i,r}$ (\Cref{subsection:equivariance-tautological-classes}) extend over the whole family $\mathcal{X}$ and restrict to the tautological line bundles at $\sigma$. Consequently the comparison isomorphism between $\HO^*(\Y)$ and $\HO^*(\Y_{\lambda})$ identifies $\kappa_0$ with $\kappa_\lambda$, that is there is a commutative diagram
\[\begin{tikzcd}
	P & {\HO^*(\Y)} \\
	& {\HO^*(\Y_{\lambda})}
	\arrow["\kappa", from=1-1, to=1-2]
	\arrow["{\kappa_{\lambda}}"', from=1-1, to=2-2]
	\arrow["\simeq", from=1-2, to=2-2]
\end{tikzcd}\]
Surjectivity and equivariance follow from \Cref{prop:Kirwan,prop:Wequiv}. The scalar class coming from the $\mathrm{B}\mathbf{G}_{\mathrm{m}}$-gerbe is unchanged since it is $W$-invariant and defined over the whole family $\mathcal{X}$. Now we can apply \Cref{lem:alternating}. The comparison preserves all cohomological degrees and the dimension of the smooth fibres.
\end{proof}

\begin{remark}
This argument does not require a generic central parameter for the original dimension vector. If $\dd=m\dd_0$, then $\lambda\cdot\dd=0$ forces $\lambda\cdot\dd_0=0$, so such a parameter could not remove all smaller zero-trace dimension vectors. The smooth family $\mathcal{X}$ coming from the extended quiver $Q^+$. For other intersection-cohomology constructions of quiver Weyl actions, with generic-orbit hypotheses, see \cite{Letellier}.
\end{remark}

\subsection{BPS cohomology of the central deformation}
We keep the virtual normalisations used above so that for a smooth Artin stack $\FM$, $\BoQ_{\FM}^{\vir}\coloneqq \BoQ_{\FM}[\dim\FM]$. For the potential
$W_\lambda=W-\sum_i\lambda_i\omega_i$, put
\[
 \mathcal K_{\lambda,\dd}^{\theta}:=\JH_{\dd,*}^{\theta} {}^p\varphi_{\Tr(W_\lambda)_{\dd}^{\theta}}\BoQ_{\FM_{\widetilde Q,\dd}^{\theta}}^{\vir},\qquad\calBPS_{\widetilde Q,W_\lambda,\dd}^{\theta}:={}^p\mathcal H^1(\mathcal K_{\lambda,\dd}^{\theta}).
\]
As before in the undeformed situation (\Cref{equation:embeddingBPS}), the canonical truncation morphism and dimensional reduction give an embedding
\[
 \Fg_{\Pi_Q^\lambda,\dd}^{\theta}:=\HO^*\!\left(\CM_{\widetilde Q,\dd}^{\theta},\calBPS_{\widetilde Q,W_\lambda,\dd}^{\theta}[-1]\right) \lhook\joinrel\longrightarrow\CA_{\Pi_Q^\lambda,\dd}^{\theta}.
\]
Its affinization $\widehat{\Fg}_{\Pi_Q^{\lambda},\dd}^{\theta}$ is the image of $\Fg_{\Pi_Q^\lambda,\dd}^{\theta}\otimes\BoQ[u]$ under the determinant-class action, exactly as in \Cref{thm:PBW}. We omit $\theta$ from the notation when $\theta=0$.

The dimensional-reduction comparison needed here is already contained in \cite[\S3.3, equation~(3.10) and Proposition~3.5]{DavBF}. The relevant case is the linear deformation, namely $n=1$ and parameter $-\lambda$ in the notation of \emph{loc.cit.}.

\begin{lemma}\label{lem:deformedBPScharacter}
For every nonzero $\dd$ with $\lambda\cdot\dd=0$, there is an isomorphism of graded vector spaces
\[
 \Fg_{\Pi_Q^\lambda,\dd}\cong\Fg_{\Pi_Q,\dd}.
\]
Consequently,
\begin{equation}\label{eq:deformedaffKac}
 \sum_k\dim\widehat{\Fg}_{\Pi_Q^\lambda,\dd}^{\,k}t^{k/2}=\frac{A_{Q,\dd}(t^{-1})}{1-t}.
\end{equation}
\end{lemma}
\begin{proof}
The linear case of \cite[equation~(3.10)]{DavBF} identifies the cohomological BPS invariants for $(\widetilde Q,W_\lambda)$ and $(\widetilde Q,W)$ whenever $\lambda\cdot\dd=0$. Dimensional reduction identifies their critical cohomology with the normalised Borel--Moore homology of the respective deformed preprojective stacks; see \cite[Proposition~3.5]{DavBF}. Applying the common shift $[-1]$ gives the stated comparison $\Fg_{\Pi_Q^{\lambda},\dd}\cong\Fg_{\Pi_Q,\dd}$. The BPS/Kac identity of \eqref{eq:BPSKac} and the determinant-class affinization give \eqref{eq:deformedaffKac}.
\end{proof}

One can see the role of the trace condition directly in this dimensional-reduction argument. Besides reducing along the loops $\omega_i$, one can reduce along the opposite arrows $a^*$. The remaining equations are
\[
 \omega_{t(a)}x_a=x_a\omega_{s(a)},\qquad a\in Q_1,
\]
and the remaining function is $-\sum_i\lambda_i\tr(\omega_i)$. The deformed dimensional-reduction theorem \cite[Theorem~1.3]{DavPadDR} applies with weight one on the $\omega_i$ and weight zero on $a,a^*$. On the common-eigenvalue support of the undeformed BPS sheaf, this function is $-(\lambda\cdot\dd)z$. It is zero when $\lambda\cdot\dd=0$. Since this function is linear, the vanishing cycle sheaf functor applies to the constant sheaf can be zero only if the function vanishes.

\subsection{Scalar support for the deformed potential}
We verify that the inclusion of the tautological subspace in the subspace of primitive elements (\Cref{prop:tautprimitive} in the undeformed situation) also holds in each deformed fibre. Write
\[
 f_{\lambda,\dd}(x,\omega):=\Tr(W_\lambda)_{\dd}(x,\omega)=\sum_i\tr\!\left(\omega_i(\mu_i(x)-\lambda_i\id_{V_i})\right)
\]
for the trace of the potential on the stack of $\dd$-dimensional representations of the tripled quiver $\widetilde{Q}$. The critical equations say that $x$ is a $\Pi_Q^\lambda$-module and that $\omega$ is an endomorphism of this module. In particular its generalised eigenspaces are submodules. Every summand of $x$ of dimension $\ee$ corresponding to a generalised eigenspace of $\omega$ satisfies $\lambda\cdot\ee=0$, by taking the trace of its preprojective relation. One can construct a coproduct on $\CA_{\Pi_Q^{\lambda}}$ using spectral factorisation of $\omega$ and restriction to disks in the $\mathbf{A}^1$-direction exactly as in \cite{DHKSV}.

\begin{lemma}\label{lem:deformedprimitives}
The spectral factorisation with respect to $\omega$ defines a graded cocommutative coproduct $\Delta_{\mathrm{fac},\lambda}$ on $\CA_{\Pi_Q^\lambda}$, compatible with the Hall product (so that $\CA_{\Pi_Q^{\lambda}}$ is a bialgebra), and
\[
 \Prim(\CA_{\Pi_Q^\lambda},\Delta_{\mathrm{fac},\lambda})
 =\widehat{\Fg}_{\Pi_Q^\lambda}.
\]
For every nonzero dimension vector $\dd\in\BoN^{Q_0}$ with $\lambda\cdot\dd=0$, one has
\begin{equation}\label{eq:deformedprimitiveinclusion}
 \HO^*(\FM_{\Pi_Q^\lambda,\dd},\BoQ)\cap v_{\dd,\lambda} \subseteq\widehat{\Fg}_{\Pi_Q^\lambda,\dd}.
\end{equation}
\end{lemma}
\begin{proof}
The proof is similar to that of \Cref{prop:primitives} and \Cref{prop:tautprimitive}, to which we refer for additional details. The construction of the coproduct is in a more general context is explained in \cite{DHKSV}. The small-diagonal support makes the BPS classes primitive, i.e. $\Fg_{\Pi_Q^{\lambda}}\subset \Prim(\CA_{\Pi_Q^{\lambda}},\Delta_{\mathrm{fac},\lambda})$. Their determinant-class multiples are primitive as well. Indeed, the classes $\delta_{\ee}=\sum_i c_1(\CV_i)$ are additive under direct sum, and the resulting operator satisfies
\[
 \Delta_{\mathrm{fac},\lambda}(\delta\cdot a)=(\delta\otimes1+1\otimes\delta)                    \Delta_{\mathrm{fac},\lambda}(a),\qquad \delta\cdot1=0.
\]
This proves that $\widehat{\Fg}_{\Pi_Q^{\lambda},\dd}=\Fg_{\Pi_Q^{\lambda}}\otimes\BoQ[\delta_{\dd}]\subset\Prim(\CA_{\Pi_Q^{\lambda},\dd})$
On a fixed nonzero component, replacing $\delta_{\ee}$ by a normalised multiple does not change the affinized subspace. The symmetrised Hall multiplication is therefore a coalgebra morphism
\[
 \Sym(\widehat{\Fg}_{\Pi_Q^\lambda})      \longrightarrow\CA_{\Pi_Q^\lambda},
\]
where the generators on the left are primitive. Cohomological integrality makes this an isomorphism of graded vector spaces, hence of coalgebras. Over $\BoQ$, the primitives of this symmetric coalgebra $\Sym(\widehat{\Fg}_{\Pi_Q^\lambda})$ are exactly its generators. This proves the first assertion.

For the second assertion, split
\[
 \omega_i=z\id_{V_i}+\eta_i,\qquad z=\frac{\sum_i\tr(\omega_i)}{|\dd|},\qquad \sum_i\tr(\eta_i)=0.
\]
The trace condition gives
\begin{equation}\label{eq:deformedscalarfunction}
 f_{\lambda,\dd}(x,\omega)=\sum_i\tr\!\left(\eta_i                  (\mu_i(x)-\lambda_i\id_{V_i})\right).
\end{equation}
Thus the scalar coordinate is a free factor and the remaining function pairs the reduced equation bundle with its dual. The zero-section Thom calculation of \Cref{thm:scalar}, applied to the section $\mu-\lambda$, identifies scalar-supported critical cohomology over a disk with $\HO^*(\FM_{\Pi_Q^\lambda,\dd},\BoQ)$, with the same degree shift. Forgetting supports and applying dimensional reduction sends $a$ to $a\cap v_{\dd,\lambda}$. This is \cite[Corollary~2.8 and Proposition~2.13]{DHKSV} applied to our situation. See \Cref{thm:scalar} for the detailed explanation in the undeformed case.

A scalar endomorphism has no decomposition into two nonzero spectral summands in disjoint disks. Its supported cohomology therefore has primitive image. Applying the first assertion gives \eqref{eq:deformedprimitiveinclusion}.
\end{proof}

\subsection{Stability restriction and tautological generation}
\begin{lemma}\label{lem:deformedstability}
For $\lambda\cdot\dd=\theta\cdot\dd=0$, open restriction induces an isomorphism
\[
 r_\theta:\widehat{\Fg}_{\Pi_Q^\lambda,\dd}\xrightarrow{\sim}\widehat{\Fg}_{\Pi_Q^\lambda,\dd}^{\theta}.
\]
For $a\in R_{\dd}$, it satisfies
\[
 r_\theta(a\cap v_{\dd,\lambda})=a|_{\FM_{\Pi_Q^\lambda,\dd}^{\theta}} \cap v_{\dd,\lambda}^{\theta}.
\]
\end{lemma}
\begin{proof}
This is the $\lambda$-deformed version of \Cref{cor:Todaaff} and \Cref{lem:restrictionvirtual}. The tripled quiver is symmetric. Therefore, Toda's restriction theorem applies to the potential $W_\lambda$: \cite[Lemma~4.10 and Remark~5.6]{Toda}, or \cite[Proposition~10.1]{DHSM}, see \Cref{thm:Toda}. It identifies the first perverse truncations through the natural restriction morphism. Dimensional reduction therefore gives the restriction isomorphism on BPS cohomology for deformed preprojective algebras.

As noted earlier, semistability in the critical locus agrees with semistability of the underlying $\Pi_Q^\lambda$-module, since a representation in the critical locus of $f_{\lambda,\dd}$ is a pair of a representation of $\Pi_Q^{\lambda}$ and an endomorphism of that representation. If that module is unstable, its maximal destabilising subobject is preserved by every endomorphism, in particular by $\omega$. Conversely, every submodule preserved by $\omega$ is a submodule of the underlying module. This also justifies dimensional reduction on the indicated semistable opens, since a priori the semistable open for the triple quiver is bigger than the pullback of the semistable open from the double quiver. We refer to \Cref{subspace:BPS-subspace} for the detailed explanation in the undeformed situation.

The determinant line bundles over $\FM_{\Pi_Q^{\lambda}}$ restrict to the corresponding line bundles on the semistable opens $\FM_{\Pi_Q^{\lambda}}^{\theta}$, so the BPS comparison extends to affinized BPS. The last formula follows from the standard compatibility of refined zero-section pullback and cap products with open restriction.
\end{proof}

\begin{theorem}\label{thm:deformed}\label{conj:deformed}
Let $Q$ be a finite quiver and $\dd\ne0$ a dimension vector for $Q$. For every deformation parameter $\lambda\in\BoC^{Q_0}$ and King stabilility parameter $\theta\in\BoQ^{Q_0}$ satisfying $\lambda\cdot\dd=\theta\cdot\dd=0$, one has
\[
 \HO^*(\BG G_{\dd})\cap[\FM_{\Pi_Q^\lambda,\dd}^{\theta}]^{\rd} =\widehat{\Fg}_{\Pi_Q^\lambda,\dd}^{\theta}.
\]
\end{theorem}
\begin{proof}
First take $\theta=0$ and put
\[
 T_{\dd,\lambda}:=R_{\dd}\cap v_{\dd,\lambda},\qquad S_{\dd,\lambda}:= \HO^*(\Y_\lambda)^{\sgn}\cap[\Y_\lambda],
\]
where the second space is endowed with the normalised CoHA grading, that is it is a graded subspace of $\HO^{\mathrm{BM}}_*(\Y_{\lambda},\BoQ^{\vir})$. By \Cref{lem:deformedprimitives},
\[
 T_{\dd,\lambda}\subseteq\widehat{\Fg}_{\Pi_Q^\lambda,\dd}.
\]
The formula of \Cref{prop:deformedgeometry} and the tautological-generation statement of \Cref{prop:deformedtarget} give a degree-preserving surjection
\[
 \Phi_\lambda:T_{\dd,\lambda}\twoheadrightarrow S_{\dd,\lambda}.
\]
Indeed, for $a\in R_{\dd}$,
\[
 \Phi_\lambda(a\cap v_{\dd,\lambda}) =(-1)^N\kappa_\lambda(a\Delta_{\dd})\cap[\Y_\lambda],
\]
and these classes span the sign summand. By \Cref{prop:deformedtarget,cor:signcharacter},
\[
 \sum_k\dim S_{\dd,\lambda}^{\,k}t^{k/2}=\frac{A_{Q,\dd}(t^{-1})}{1-t}.
\]
This is also the character of $\widehat{\Fg}_{\Pi_Q^\lambda,\dd}$, by \Cref{lem:deformedBPScharacter}. All the graded pieces are finite-dimensional. For each $k$,
\[
 \dim S_{\dd,\lambda}^{\,k}\leq\dim T_{\dd,\lambda}^{\,k} \leq\dim\widehat{\Fg}_{\Pi_Q^\lambda,\dd}^{\,k} =\dim S_{\dd,\lambda}^{\,k}.
\]
Hence $T_{\dd,\lambda}=\widehat{\Fg}_{\Pi_Q^\lambda,\dd}$.

Finally, apply \Cref{lem:deformedstability}. Its compatibility with the reduced virtual class and tautological cap products gives
\[
 R_{\dd}\cap v_{\dd,\lambda}^{\theta}=r_\theta(T_{\dd,\lambda}) =\widehat{\Fg}_{\Pi_Q^\lambda,\dd}^{\theta}.
\]
\end{proof}

\subsection{A generation question for $2$-Calabi--Yau categories}
\label{subsection:generation-2CY}
The first inclusion (\Cref{prop:tautprimitive} and \eqref{eq:deformedprimitiveinclusion}) in the proof of the tautological generation theorem (\Cref{thm:main}) has a wider scope. For a moduli stack of objects in a $2$-Calabi--Yau Abelian category with the appropriate geometric assumptions (see \cite{DHSM}, \cite{bu2025cohomology}, in particular, the moduli stack has a good moduli space), there is an inclusion
\[
 \HO^*(\FM,\BoQ)\cap[\FM]^{\rd}\subseteq\Prim(\HO_*^{\BM}(\FM)).
\]
Here $[\FM]^{\rd}$ is the reduced virtual fundamental class of $\FM$. For sheaves on symplectic surfaces, this is the starting point of the tautological-generation theorem in \cite{DHKSV}.

The opposite inclusion, i.e. that any primitive elements of the cohomological Hall algebra is a cohomology multiple of the reduced fundamental class, needs additional geometry. In the present case of quiver varieties, it is given by extending the quiver by legs and imposing a generic stability condition, providing a smooth moduli space with tautological generation, and a sign character formula for a symmetric group action. It is natural to ask for analogues of this construction for other moduli categories. For moduli spaces of sheaves on symplectic surfaces as in \cite{DHKSV}, the tautological generation of affinized BPS cohomology (the generalisation of Markman's theorem to singular moduli spaces) is proved by using Hecke modifications at points of the surface.

Even if the notion of tautological classes is not well-defined for the moduli stack $\FM$ under consideration, one can still ask whether there is an equality
\[
 \HO^*(\FM,\BoQ)\cap[\FM]^{\rd}=\Prim(\HO_*^{\BM}(\FM,\BoQ)),
\]
that is using the whole of the cohomology of the stack $\HO^*(\FM,\BoQ)$ instead of just the a priori smaller subalgebra generated by tautological classes.

We know that this equality holds for moduli stacks of sheaves on symplectic surfaces by \cite{DHKSV} and of representations of preprojective algebras of quivers by the present paper since in both cases, the cohomology of the stack contains tautological classes.

The second author discussed this general generation question with Tasuki Kinjo during a workshop at KIAS, Seoul in June 2025 and this seems a worthy future direction of research.

\enlargethispage{\baselineskip}
\begingroup
\small
\setstretch{1.3}
\bibliographystyle{amsalpha}
\bibliography{tautological-generation-references}
\endgroup
\end{document}